\documentclass{amsart}
\usepackage{amssymb}
\usepackage{amsmath,amsthm}
\usepackage{afterpage}
\usepackage{graphicx}
\title[Elliptic surfaces without 1-handles]{Elliptic surfaces without 1-handles}
\author[Kouichi Yasui]{Kouichi Yasui}
\date{February 22, 2008. Revised on July 6, 2008.}

\address{Department~of~Mathematics, Graduate~School~of~Science, Osaka~University, Toyonaka, Osaka 560-0043, Japan}
\email{kyasui@cr.math.sci.osaka-u.ac.jp}
\thanks{The author is partially supported by JSPS Research Fellowships for Young Scientists.}

\newtheorem{theorem}{Theorem}[section]
\newtheorem{proposition}[theorem]{Proposition}
\newtheorem{lemma}[theorem]{Lemma}
\newtheorem{corollary}[theorem]{Corollary}
\theoremstyle{definition}

\newtheorem{remark}[theorem]{Remark}

\newtheorem{ack}{Acknowledgement}

\begin{document}

\begin{abstract}
Harer-Kas-Kirby conjectured that every handle decomposition of the elliptic surface $E(1)_{2,3}$ requires both $1$- and $3$-handles. 
We prove that the elliptic surface $E(n)_{p,q}$ has a handle decomposition without 1-handles for $n\geq 1$ and $(p,q)=(2,3),(2,5),(3,4),(4,5)$. 
\end{abstract}

\maketitle

\section{Introduction}
It is not known whether or not the 4-sphere $S^4$ and the complex projective plane $\mathbf{CP}^2$ admit an exotic smooth structure. If such a structure exists, then each handle decomposition of it has at least either a $1$- or $3$-handle (cf.~\cite{Y2}). On the contrary, many simply connected closed topological $4$-manifolds are known to 
admit infinitely many different smooth structures which have neither $1$- nor $3$-handles in their handle decompositions (cf.~Gompf-Stipsicz \cite{GS}).

Problem 4.18 in Kirby's problem list \cite{Ki} is the following: ``Does every simply connected, closed $4$-manifold 
have a handlebody decomposition without $1$-handles? Without $1$- and $3$-handles?'' It is not known whether or not the simply connected elliptic surface $E(n)_{p,q}$ 
($n\geq 1$, $p,q\geq 2$, $\gcd (p,q)=1$) 
admits a handle decomposition without $1$-handles. 
In particular, Harer, Kas and Kirby conjectured in \cite{HKK} that 
every handle decomposition of $E(1)_{2,3}$ requires both $1$- and $3$-handles. Gompf \cite{G} notes the following: it is a good conjecture that $E(n)_{p,q}$ $(p,q\geq 2)$ has no handle decomposition without $1$- and $3$-handles.

In \cite{Y1} and \cite{Y2} we constructed a homotopy $E(1)_{2,3}$ which has the same Seiberg-Witten invariant as $E(1)_{2,3}$ and has a handle decomposition without $1$- and $3$-handles. 
Recently Akbulut \cite{Ak} proved that $E(1)_{2,3}$ has a handle decomposition without 1- and 3-handles, by using knot surgery on $E(1)$ and investigating a dual handle decomposition. He also proved that infinitely many different smooth structures on $\mathbf{C}\mathbf{P}^2\# 9\overline{\mathbf{C}\mathbf{P}^2}$ admit handle decompositions without 1-handles. 

In this paper, we prove the theorem below by improving our previous procedure (\cite{Y1}, \cite{Y2}). Our method is different from Akbulut. 
\begin{theorem}\label{th:1.1}
The elliptic surface $E(n)_{p,q}$ has a handle decomposition without $1$-handles, for $n\geq 1$ and $(p,q)=(2,3),(2,5),(3,4),(4,5)$. 
\end{theorem}
\begin{ack}
The author wishes to express his deep gratitude to his adviser, 
Hisaaki Endo, for encouragements and helpful comments. 
The author would like to thank Selman Akbulut, Motoo Tange and Yuichi Yamada for discussions. 
This work was partially done during the author's stay at Michigan State University. The author is grateful for their hospitality. Finally, the author would like to thank the referee for his/her useful suggestion (see the proof of Proposition~\ref{prop:E(n)_4}). 
\end{ack}
\section{Rational blow-down}
In this section we review the rational blow-down introduced 
by Fintushel-Stern \cite{FS1}. See also Gompf-Stipsicz \cite{GS}.

\begin{figure}[htbp]
\begin{center}
\includegraphics[width=3.5in]{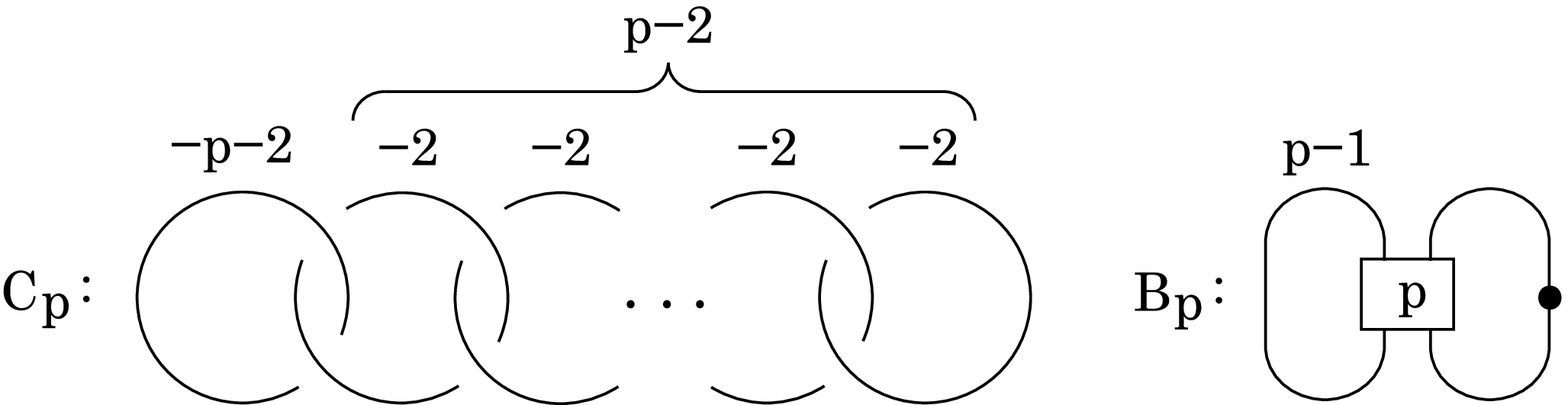}
\caption{}
\label{fig1}
\end{center}
\end{figure}

Let $C_p$ and $B_p$ ($p\geq 2$) be the smooth $4$-manifolds defined by Kirby diagrams in Figure~\ref{fig1}.
The boundary $\partial C_p$ of $C_p$ is diffeomorphic to the lens space $L(p^2,p-1)$ and to the boundary $\partial B_p$ of $B_p$. 

Suppose that $C_p$ embeds in a smooth $4$-manifold $X$. Let $X_{(p)}$ be a smooth $4$-manifold obtained from $X$ by removing $C_p$ and gluing $B_p$. The $4$-manifold $X_{(p)}$ is called the rational blow-down of $X$ along $C_p$. Note that $X_{(p)}$ is uniquely determined up to diffeomorphism by a fixed pair $(X,C_p)$. 
This operation has the following relation with the logarithmic transformation.
\begin{theorem}[{Fintushel-Stern \cite{FS1}, see also Gompf-Stipsicz \cite{GS}}]\label{th:2.1}
Suppose that a smooth $4$-manifold $X$ contains a cusp neighborhood, that is, a $0$-handle with a $2$-handle attached along a $0$-framed right trefoil knot. 
Let $X_{p}$ be the smooth $4$-manifold obtained from $X$ by 
performing a logarithmic transformation of multiplicity $p$ in the cusp neighborhood. 
Then there exists a copy of $C_p$ in 
$X\# (p-1)\overline{\mathbf{C}\mathbf{P}^2}$ such that the rational blow-down of 
$X\# (p-1)\overline{\mathbf{C}\mathbf{P}^2}$ along the copy of $C_p$ is diffeomorphic to $X_p$.
\end{theorem}
\section{Proof}
In this section we prove Theorem 1.1. We do not draw (whole) Kirby diagrams of elliptic surfaces. However, one can draw whole diagrams of elliptic surfaces without 1-handles, following the procedures in this section. 

Let $E(n)$ be the simply connected elliptic surface with Euler characteristic $12n$ and with no multiple fibers, and $E(n)_{p_1,\dots,p_k}$ the elliptic surface obtained from $E(n)$ by performing logarithmic transformations of multiplicities $p_1,\dots,p_k$.
\begin{proposition}\label{prop:3.1}
For $n\geq 1$, the elliptic surface $E(n)_2$ has handle decompositions as in Figure~$\ref{fig2}$ and $\ref{fig3}$. Each obvious cusp neighborhood in Figure~$\ref{fig2}$ and $\ref{fig3}$ is isotopic to a regular neighborhood of a cusp fiber of $E(n)_2$. 
\begin{figure}[ht]
\begin{center}
\begin{minipage}{.48\linewidth}
\begin{center}
\includegraphics[width=2.2in]{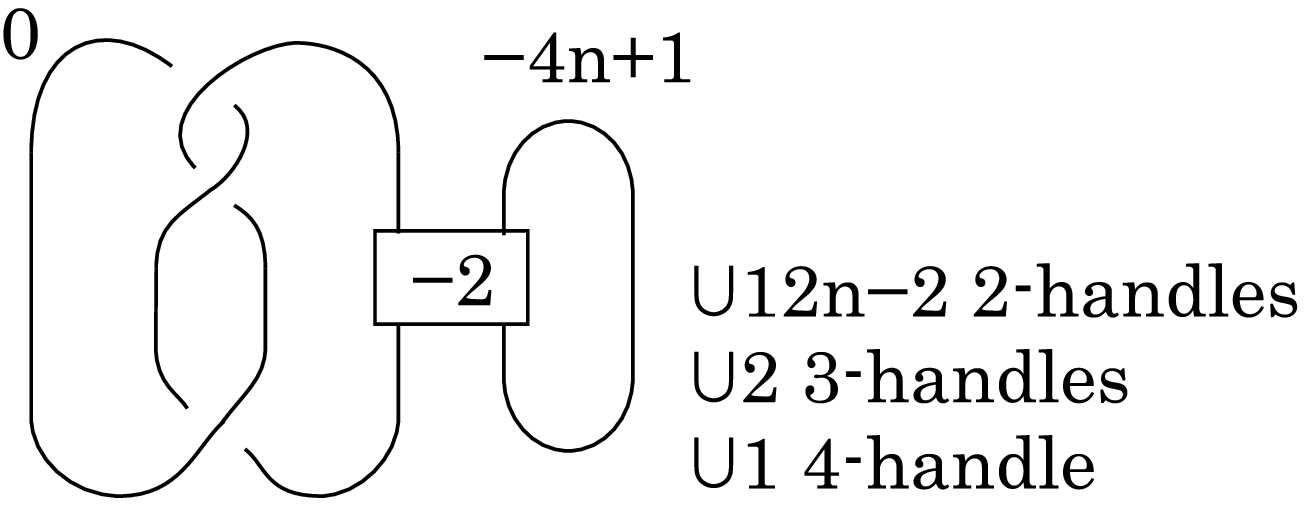}
\caption{$E(n)_2$}
\label{fig2}
\end{center}
\end{minipage}
\begin{minipage}{.48\linewidth}
\begin{center}
\includegraphics[width=2.2in]{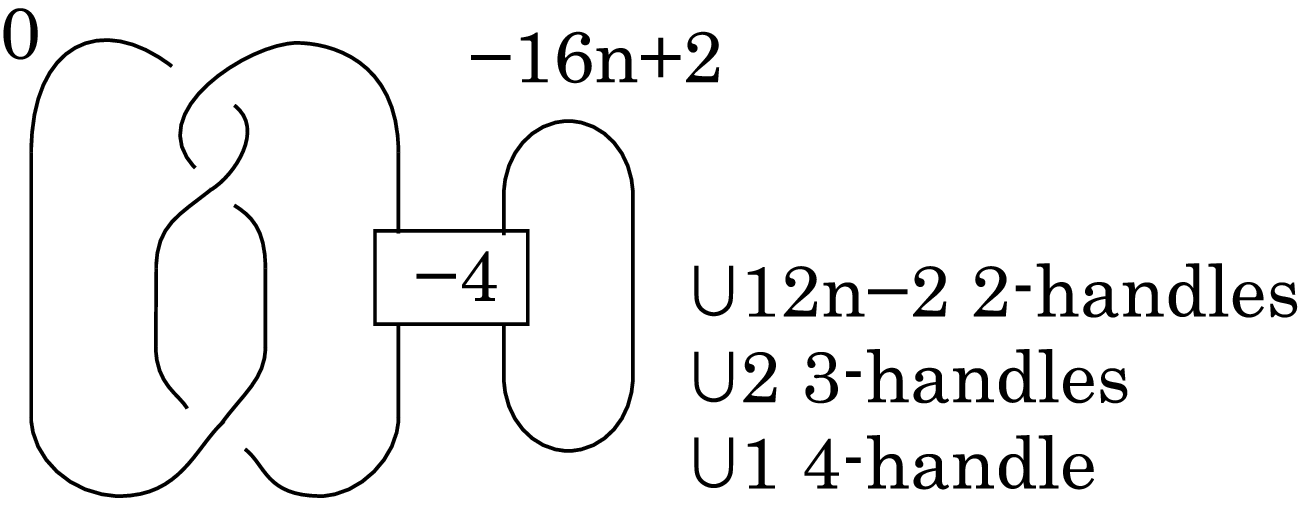}
\caption{$E(n)_2$}
\label{fig3}
\end{center}
\end{minipage}
\end{center}
\end{figure}
\end{proposition}
\begin{proof}
$E(n)_p$ admits a handle decomposition in Figure~$\ref{fig9}$ (see Gompf-Stipsicz \cite[page $315 \sim 316$]{GS} and Harer-Kas-Kirby \cite{HKK}). The obvious cusp neighborhood in Figure~$\ref{fig9}$ is isotopic to a regular neighborhood of a cusp fiber of $E(n)_p$ (see \cite{GS} and \cite{HKK}). Figure~$\ref{fig10}$ is the $p=2$ case of Figure~$\ref{fig9}$. Note that we do not draw $6n-1$ 2-handles in Figure~$\ref{fig10}$. 
We change Figure~$\ref{fig10}$ into Figure~$\ref{fig2}$ and $\ref{fig3}$ without sliding the cusp neighborhood over any handles, as follows. In Figure~$\ref{fig10}$, we slide $-4n+2$ framed knot over vertical $-1$ framed knots as shown in Figure~$\ref{fig10}\sim \ref{fig13}$. Note that $\frac{k}{2}$ in the boxes denotes $k$ right half-twists. By repeating handle slides similar to Figure~$\ref{fig11}\sim \ref{fig13}$, we obtain Figure~$\ref{fig14}$. An isotopy gives Figure~$\ref{fig15}$. By canceling $1$-handles, we get Figure~$\ref{fig2}$. 

In Figure~$\ref{fig15}$, we slide a vertical $-1$ framed knot over $-4n+1$ framed knot. We get Figure~$\ref{fig16}$. We slide $-4n$ framed knot over $-4n+1$ framed knot as shown in Figure~$\ref{fig17}$. Sliding $-16n+3$ framed knot over a vertical $-1$ framed knot gives Figure~$\ref{fig18}$. By repeating handle slides similar to Figure~$\ref{fig11}\sim \ref{fig13}$, we obtain Figure~$\ref{fig19}$. An isotopy gives Figure~$\ref{fig20}$. By canceling $1$-handles, we get Figure~$\ref{fig3}$. 
\end{proof}
\begin{proposition}\label{prop:3.2}
For $n\geq 1$, the elliptic surface $E(n)_3$ admits a handle decomposition as in Figure~$\ref{fig4}$. The obvious cusp neighborhood in Figure~$\ref{fig4}$ is isotopic to a regular neighborhood of a cusp fiber of $E(n)_3$. 
\begin{figure}[ht!]
\begin{center}
\includegraphics[width=2.2in]{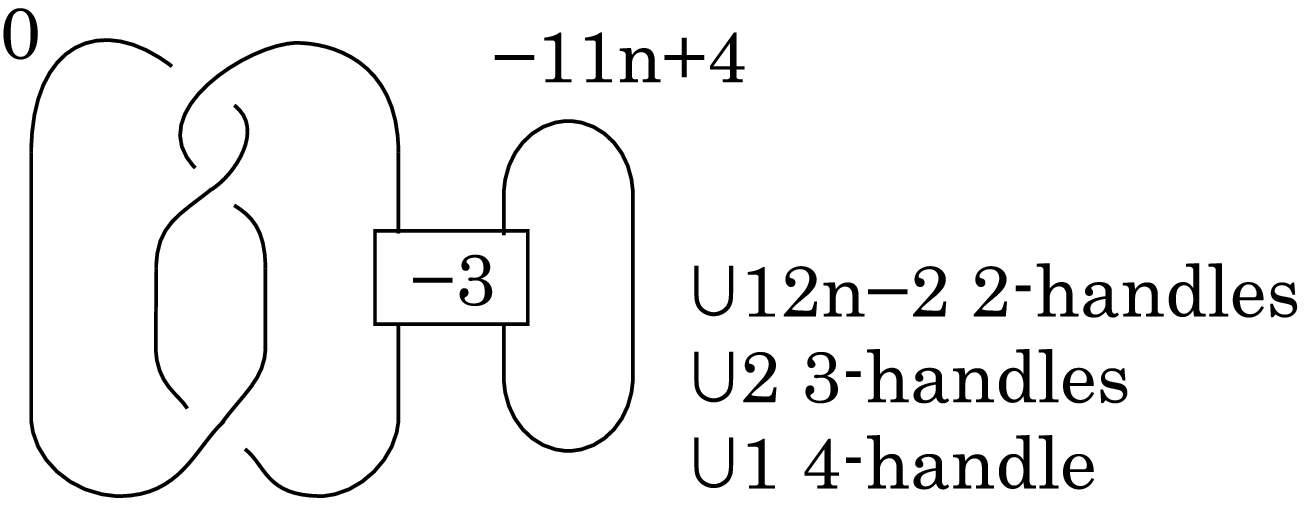}
\caption{$E(n)_3$}
\label{fig4}
\end{center}
\end{figure}
\end{proposition}
\begin{proof}
Figure~$\ref{fig21}$ is the $p=3$ case of Figure~$\ref{fig9}$. 
We change Figure~$\ref{fig21}$ into Figure~$\ref{fig4}$ without sliding the cusp neighborhood over any handles. In Figure~$\ref{fig21}$, we slide $-9n+3$ framed knot over vertical $-1$ framed knots as shown in Figure~$\ref{fig21}\sim \ref{fig28}$. An isotopy gives Figure~\ref{fig29}. We get Figure~\ref{fig30} by sliding $-9n+2$ framed knot over a vertical $-1$ framed knot. In the $n\geq 2$ case, we obtain Figure~\ref{fig31} by repeating handle slides similar to Figure~$\ref{fig21}\sim \ref{fig30}$. We get the $n=1$ case of Figure~\ref{fig31} by an isotopy in the $n=1$ case of Figure~$\ref{fig21}$. Handle slides similar to Figure~$\ref{fig21}\sim \ref{fig24}$ give Figure~\ref{fig32}. An isotopy gives Figure~\ref{fig33}. By canceling $1$-handles, we get Figure~$\ref{fig4}$.
\end{proof}
\begin{proposition}\label{prop:E(n)_4}
For $n\geq 1$, the elliptic surface $E(n)_4$ admits a handle decomposition as in Figure~$\ref{fig5}$. The obvious cusp neighborhood in Figure~$\ref{fig5}$ is isotopic to a regular neighborhood of a cusp fiber of $E(n)_4$. 
\begin{figure}[ht!]
\begin{center}
\includegraphics[width=2.2in]{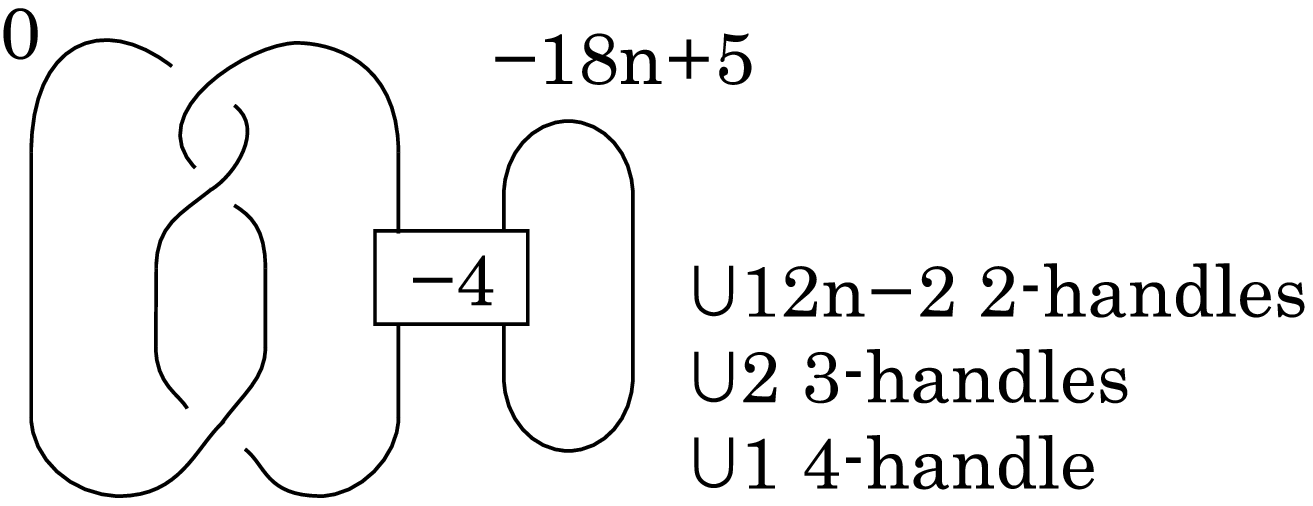}
\caption{$E(n)_4$}
\label{fig5}
\end{center}
\end{figure}
\end{proposition}
\begin{proof}
In Figure~\ref{fig9} of $E(n)_p$, we repeat handle slides shown in Figure~\ref{fig34}. We then get the diagram of $E(n)_p$ in Figure~\ref{fig35}. (This diagram is a key of our proof for $n\geq 2$. The way to construct this diagram is suggested by the referee.) Note that we did not slide the cusp neighborhood in Figure~\ref{fig9} over any handles. Figure~$\ref{fig36}$ is the $p=4$ case of Figure~$\ref{fig35}$. We change Figure~$\ref{fig36}$ into Figure~$\ref{fig5}$ without sliding the cusp neighborhood over any handles, as follows. 

The $n\geq 2$ case.  We slide handles as shown in Figure~$\ref{fig36}\sim \ref{fig46}$. We then get Figure~\ref{fig47}. Isotopies give Figure~\ref{fig48} and \ref{fig49}. We have Figure~\ref{fig53} by handle slide as shown in Figure~$\ref{fig49}\sim \ref{fig52}$. By repeating handle slides similar to Figure~$\ref{fig36}\sim \ref{fig53}$, we obtain Figure~\ref{fig54}. We slide handles similarly to Figure~$\ref{fig36}\sim \ref{fig43}$. We then get Figure~\ref{fig55}. An isotopy gives Figure~\ref{fig56}. By cancelling 1-handles, we have Figure~\ref{fig5}.

The $n=1$ case. Figure~$\ref{fig54}$ is isotopic to the $n=1$ case of Figure~\ref{fig36}. We slide handles similarly to Figure~$\ref{fig36}\sim \ref{fig43}$. We then get Figure~\ref{fig55}. An isotopy gives Figure~\ref{fig56}. By cancelling 1-handles, we have Figure~\ref{fig5}.
\end{proof}
\begin{lemma}\label{without-handle}Suppose that a simply connected closed smooth $4$-manifold $X$ has a handle decomposition 
as in Figure~$\ref{fig6}$. Here $q$ is an arbitrary integer. $h_2$ and $h_3$ are arbitrary non-negative integers. 
Let $X_{(p)}$ be the rational blow-down of $X$ along the copy of $C_p$ in Figure~$\ref{fig6}$. 
Then $X_{(p)}$ admits a handle decomposition 
\begin{equation*}
X_{(p)}=\text{one $0$-handle} \cup \text{$(h_2+1)$ $2$-handles} \cup \text{$h_3$ $3$-handles} \cup \text{one $4$-handle}.
\end{equation*}
In particular $X_{(p)}$ admits a handle decomposition without $1$-handles.
\begin{figure}[h!]
\begin{center}
\includegraphics[width=3.4in]{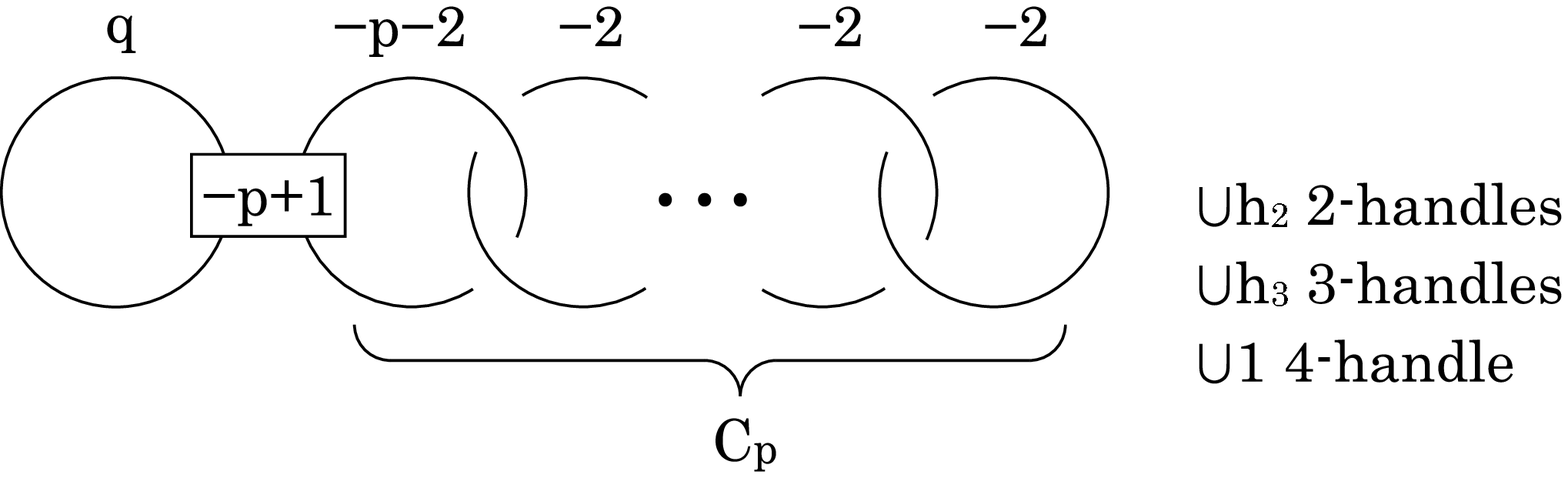}
\caption{Handle decomposition of $X$}
\label{fig6}
\end{center}
\end{figure}
\end{lemma}
\begin{proof}Draw a Kirby diagram of $X_{(p)}$, following the procedure introduced in Gompf-Stipsicz~\cite[Section 8.5]{GS} (see also \cite[page 516 Solution of Exercise 8.5.1.(a)]{GS}). Then we have a handle decomposition of $X_{(p)}$ as in Figure~\ref{fig7}. We easily get a meridian of the unique dotted circle by a handle slide. Thus we can cancel the $1$-handle/$2$-handle pair. 
\begin{figure}[ht]
\begin{center}
\includegraphics[width=2.2in]{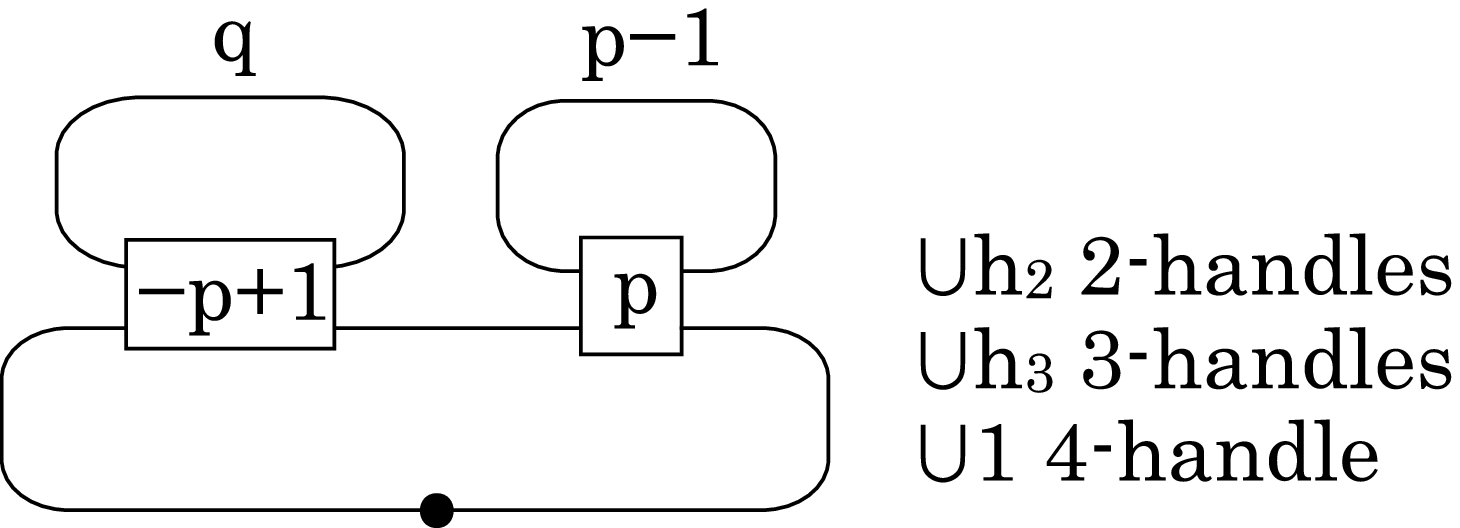}
\caption{Handle decomposition of $X_{(p)}$}
\label{fig7}
\end{center}
\end{figure}
\end{proof}
\begin{corollary}\label{cor:key}
Suppose that a simply connected closed smooth $4$-manifold $X$ has a handle decomposition as in Figure~$\ref{fig8}$. Here $q$ is an arbitrary integer. $h_2$ and $h_3$ are arbitrary non-negative integers. Let $X_{p}$ be the smooth $4$-manifold obtained from $X$ by performing a logarithmic transformation of multiplicity $p$ in the obvious cusp neighborhood in Figure~$\ref{fig8}$. 
Then $X_p$ admits a handle decomposition 
\begin{equation*}
X_p=\text{one $0$-handle} \cup \text{$(h_2+2)$ $2$-handles} \cup \text{$h_3$ $3$-handles} \cup \text{one $4$-handle}.
\end{equation*}
In particular $X_p$ admits a handle decomposition without $1$-handles.
\begin{figure}[ht]
\begin{center}
\includegraphics[width=2.2in]{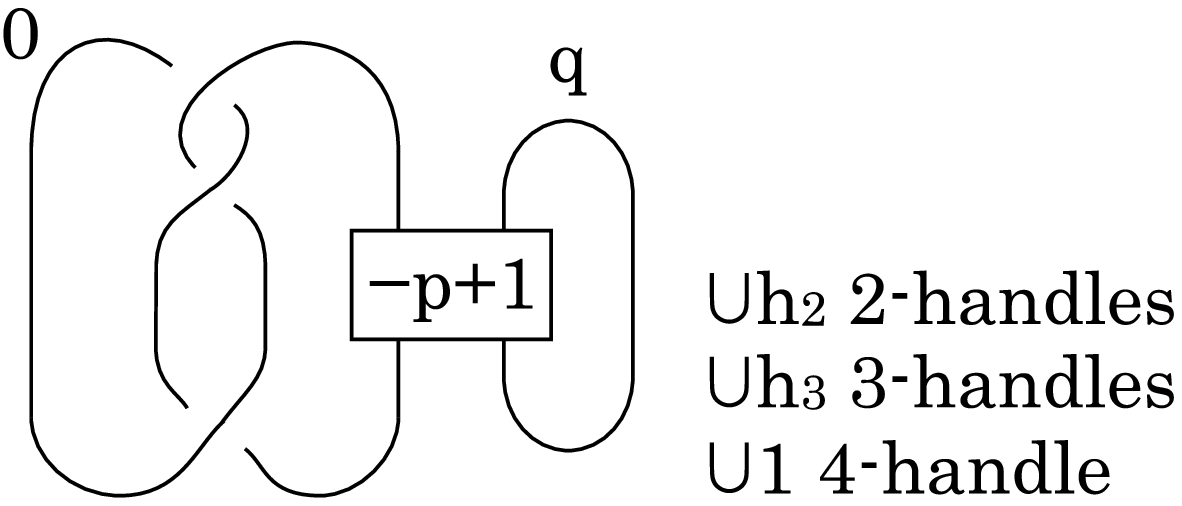}
\caption{Handle decomposition of $X$}
\label{fig8}
\end{center}
\end{figure}
\end{corollary}
\begin{proof}
Construct $C_p$ from Figure $\ref{fig8}$, following the procedure given by Fintushel-Stern~\cite[Example 1]{FS1} (and Gompf-Stipsicz~\cite[Section 8.5]{GS}). Then we have an embedding of $C_p$ into $X\# (p-1)\overline{\mathbf{C}\mathbf{P}^2}$ such that the rational blow-down of $X\# (p-1)\overline{\mathbf{C}\mathbf{P}^2}$ along $C_p$ is diffeomorphic to $X_p$. This embedding of $C_p$ clearly satisfies the assumption of Lemma~\ref{without-handle}. Therefore we get the required handle decomposition of $X_p$.
\end{proof}
\begin{remark}
One can prove Corollary~\ref{cor:key} without using rational blow-downs. Follow the procedure given by Gompf \cite[Section 4]{G}. 
\end{remark}
Propositions~\ref{prop:3.1}, \ref{prop:3.2} and \ref{prop:E(n)_4} together with Corollary~\ref{cor:key} clearly give the following main theorem:
\begin{theorem}
For $n\geq 1$ and $(p,q)=(2,3),(2,5),(3,4),(4,5)$, the elliptic surface $E(n)_{p,q}$ has a handle decomposition
\begin{equation*}
\text{one $0$-handle} \cup \text{$12n$ $2$-handles} \cup \text{two $3$-handles} \cup \text{one $4$-handle}.
\end{equation*}
\end{theorem}
\section{Further remarks}
\begin{figure}[p]
\begin{center}
\includegraphics[width=3.1in]{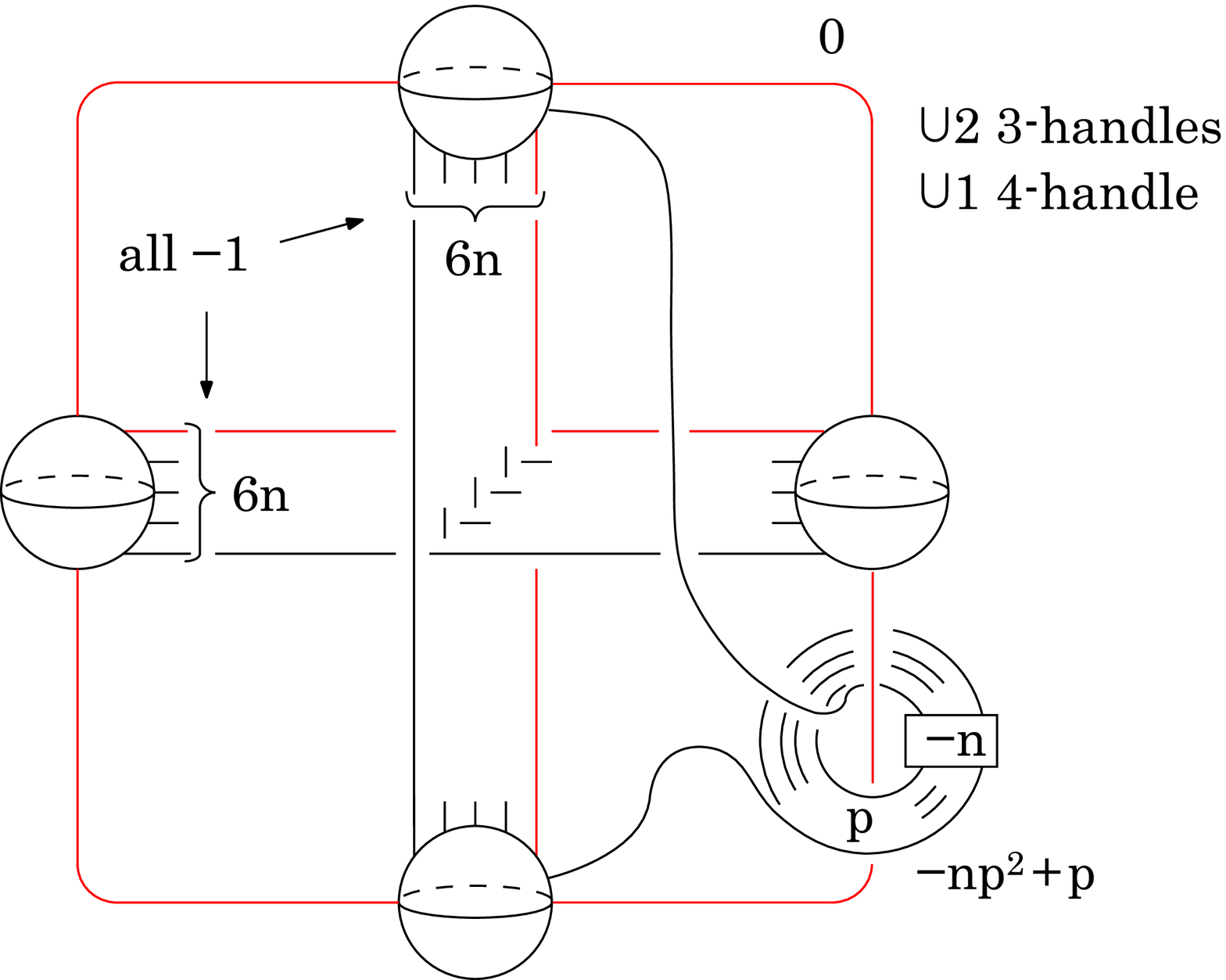}
\caption{$E(n)_p$}
\label{fig9}
\end{center}
\bigskip \medskip

\begin{center}
\includegraphics[width=3.1in]{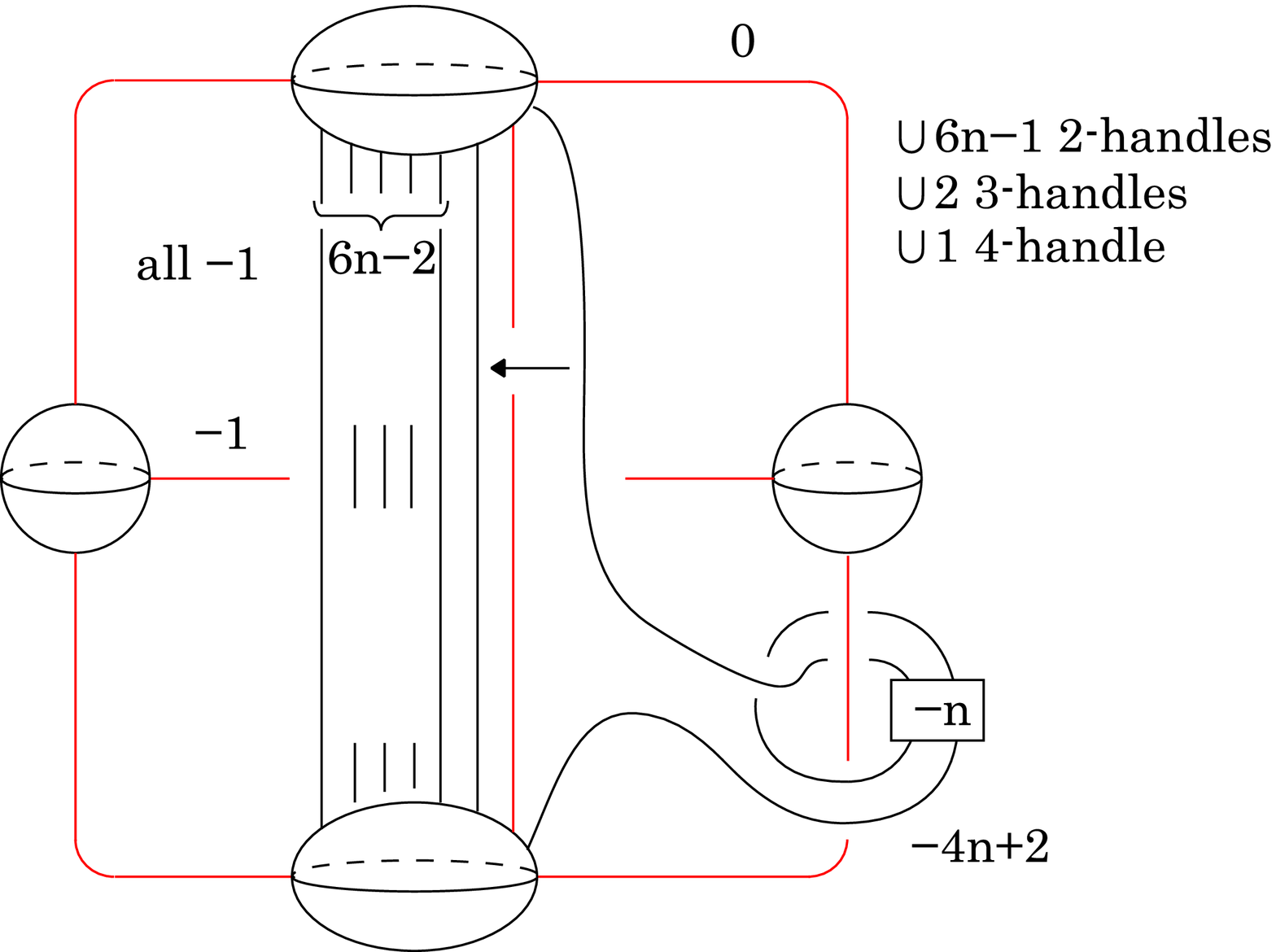}
\caption{$E(n)_2$}
\label{fig10}
\end{center}
\bigskip \medskip

\begin{center}
\includegraphics[width=3.1in]{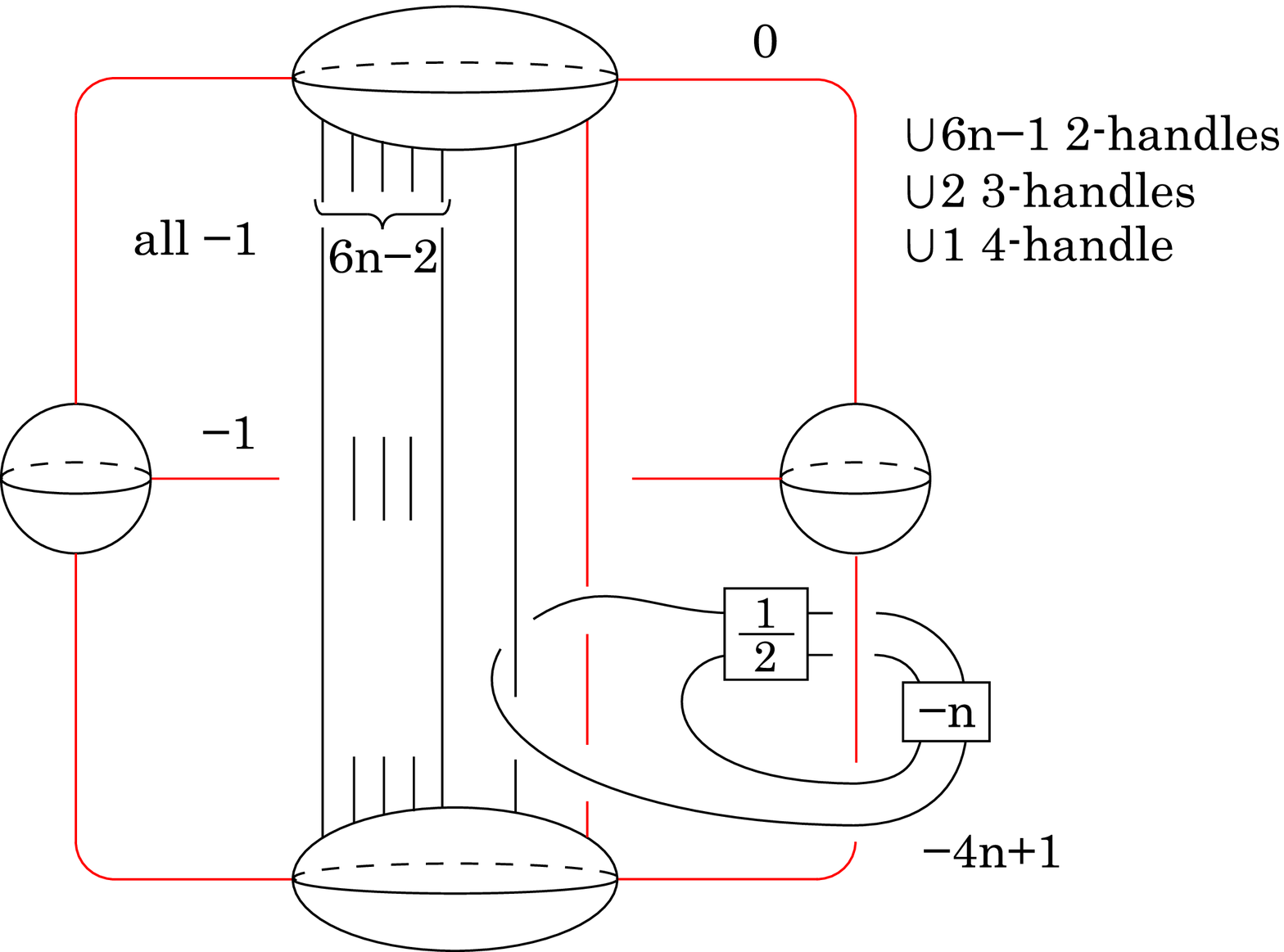}
\caption{$E(n)_2$}
\label{fig11}
\end{center}
\end{figure}

\begin{figure}[p]
\begin{center}
\includegraphics[width=3.1in]{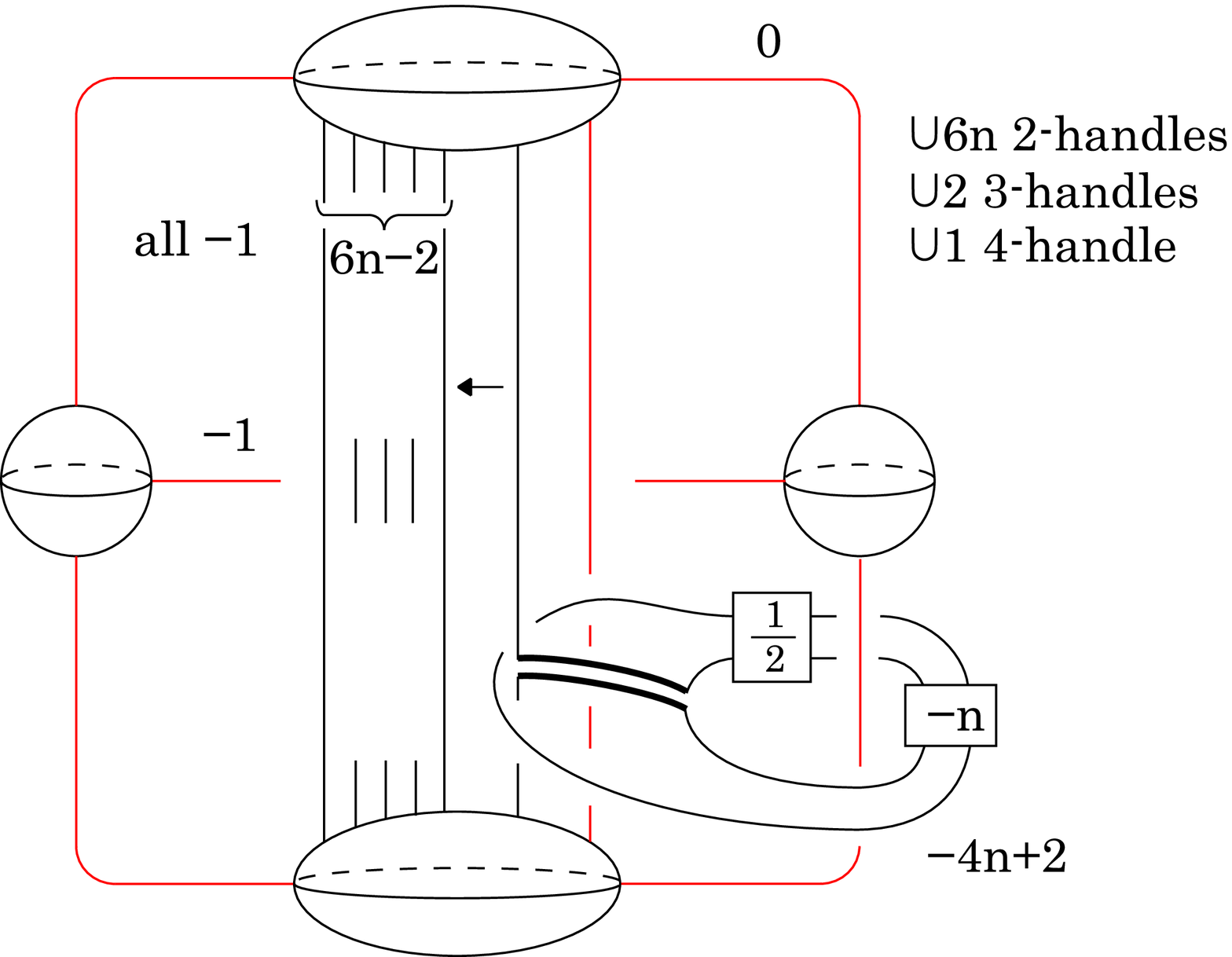}
\caption{$E(n)_2$}
\label{fig12}
\end{center}
\bigskip \medskip

\begin{center}
\includegraphics[width=3.1in]{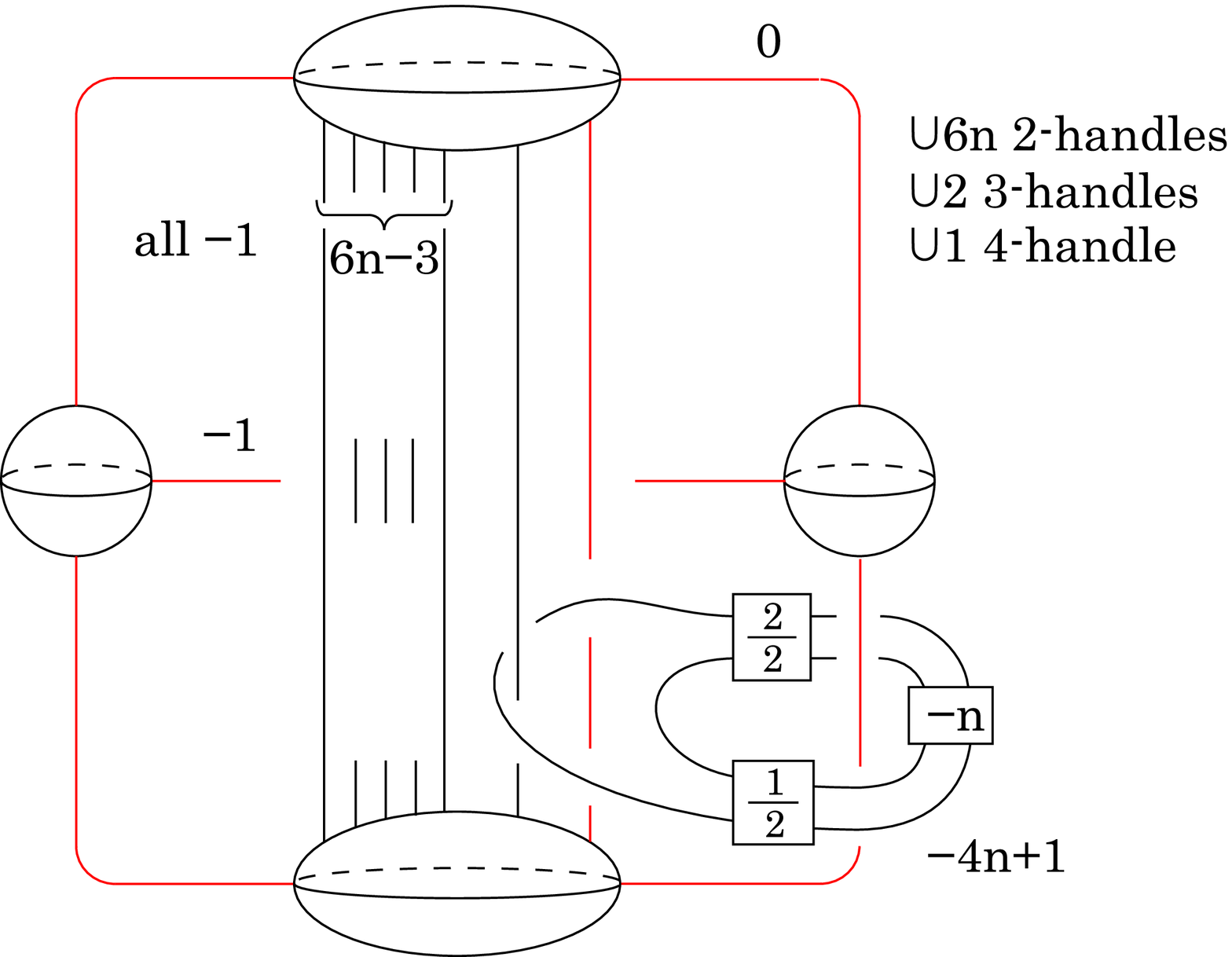}
\caption{$E(n)_2$}
\label{fig13}
\end{center}
\bigskip \medskip

\begin{center}
\includegraphics[width=3.1in]{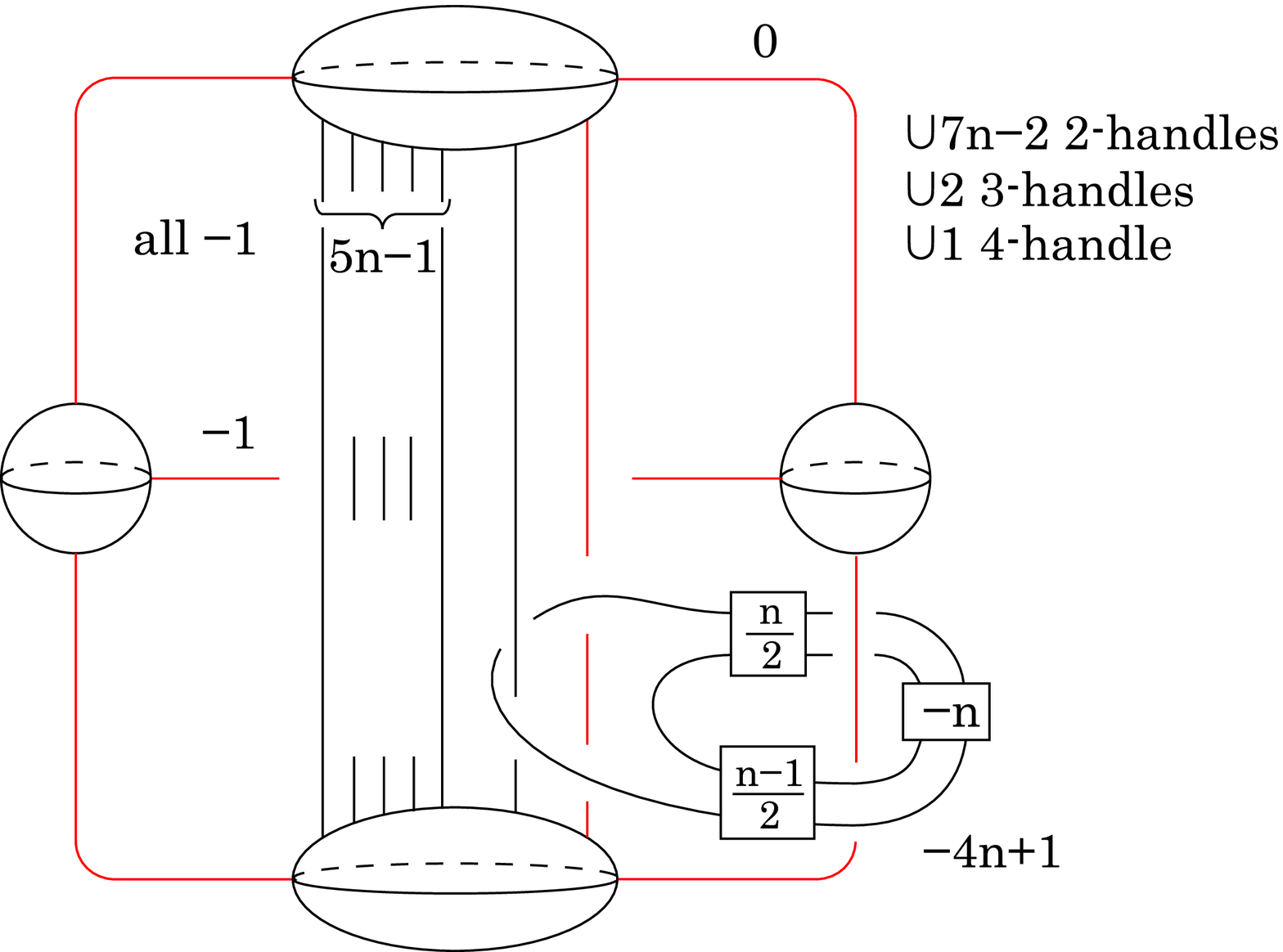}
\caption{$E(n)_2$}
\label{fig14}
\end{center}
\end{figure}

\begin{figure}[p]
\begin{center}
\includegraphics[width=3.1in]{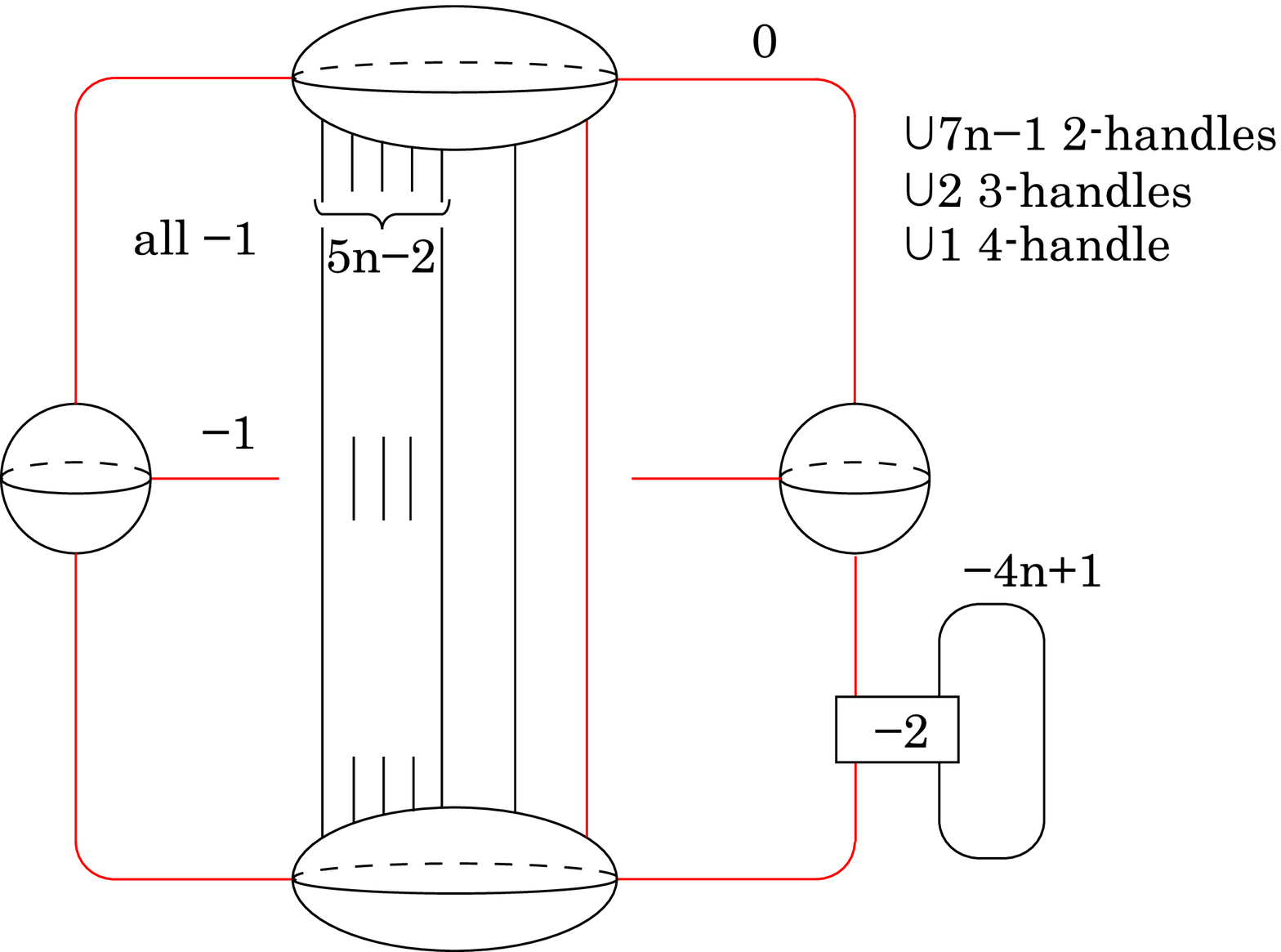}
\caption{$E(n)_2$}
\label{fig15}
\end{center}
\bigskip \medskip

\begin{center}
\includegraphics[width=3.1in]{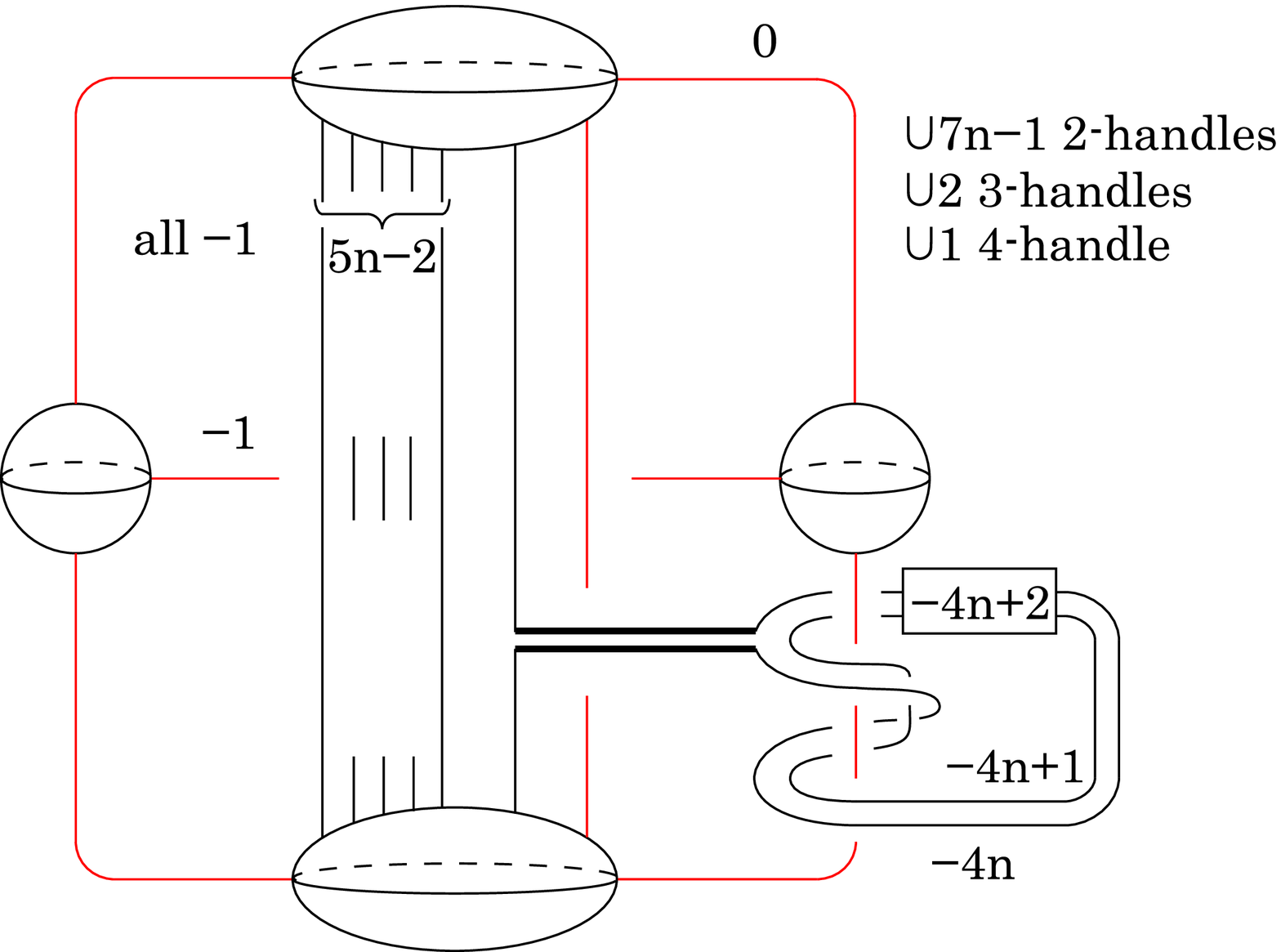}
\caption{$E(n)_2$}
\label{fig16}
\end{center}
\bigskip \medskip

\begin{center}
\includegraphics[width=3.1in]{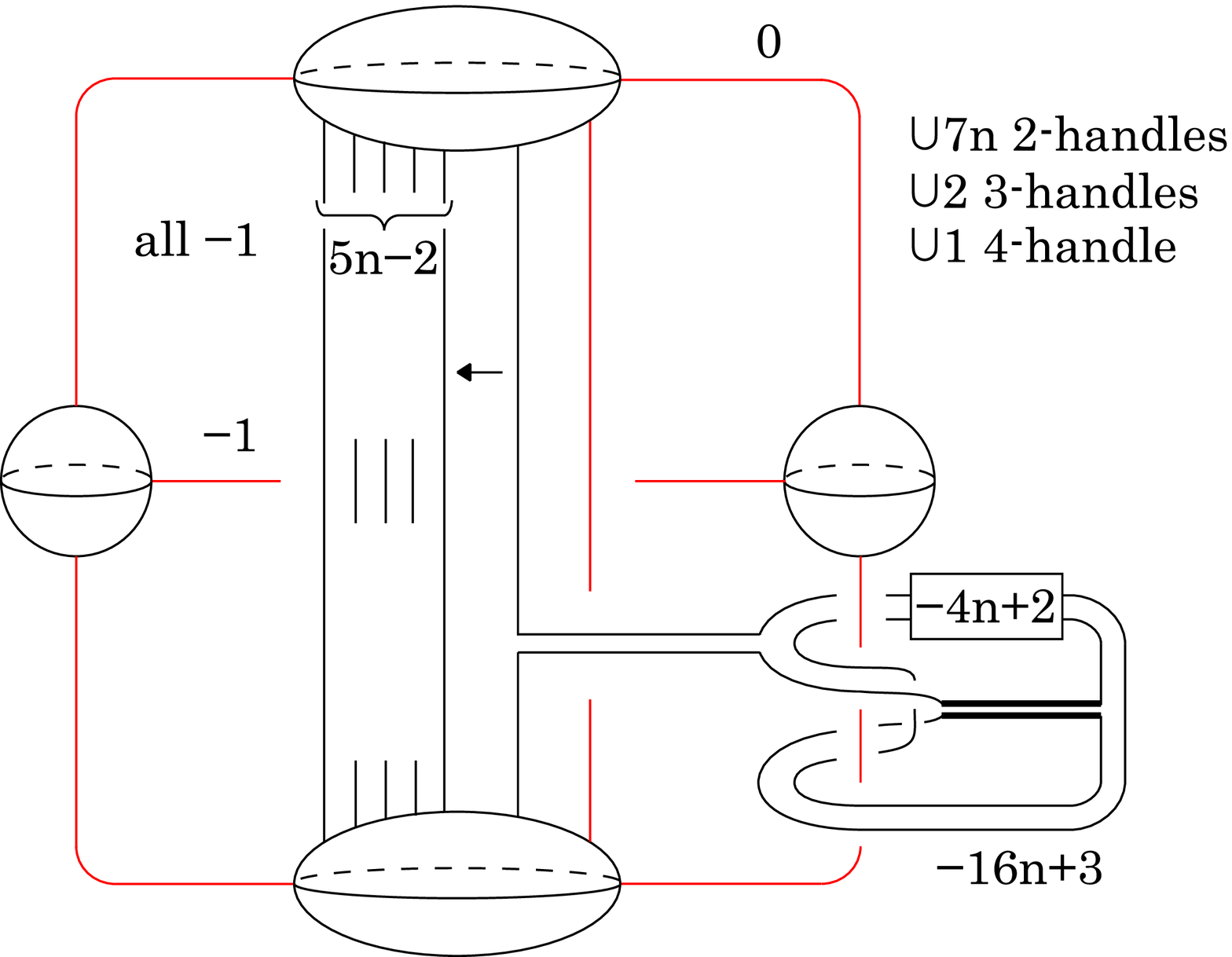}
\caption{$E(n)_2$}
\label{fig17}
\end{center}
\end{figure}

\begin{figure}[p]
\begin{center}
\includegraphics[width=3.1in]{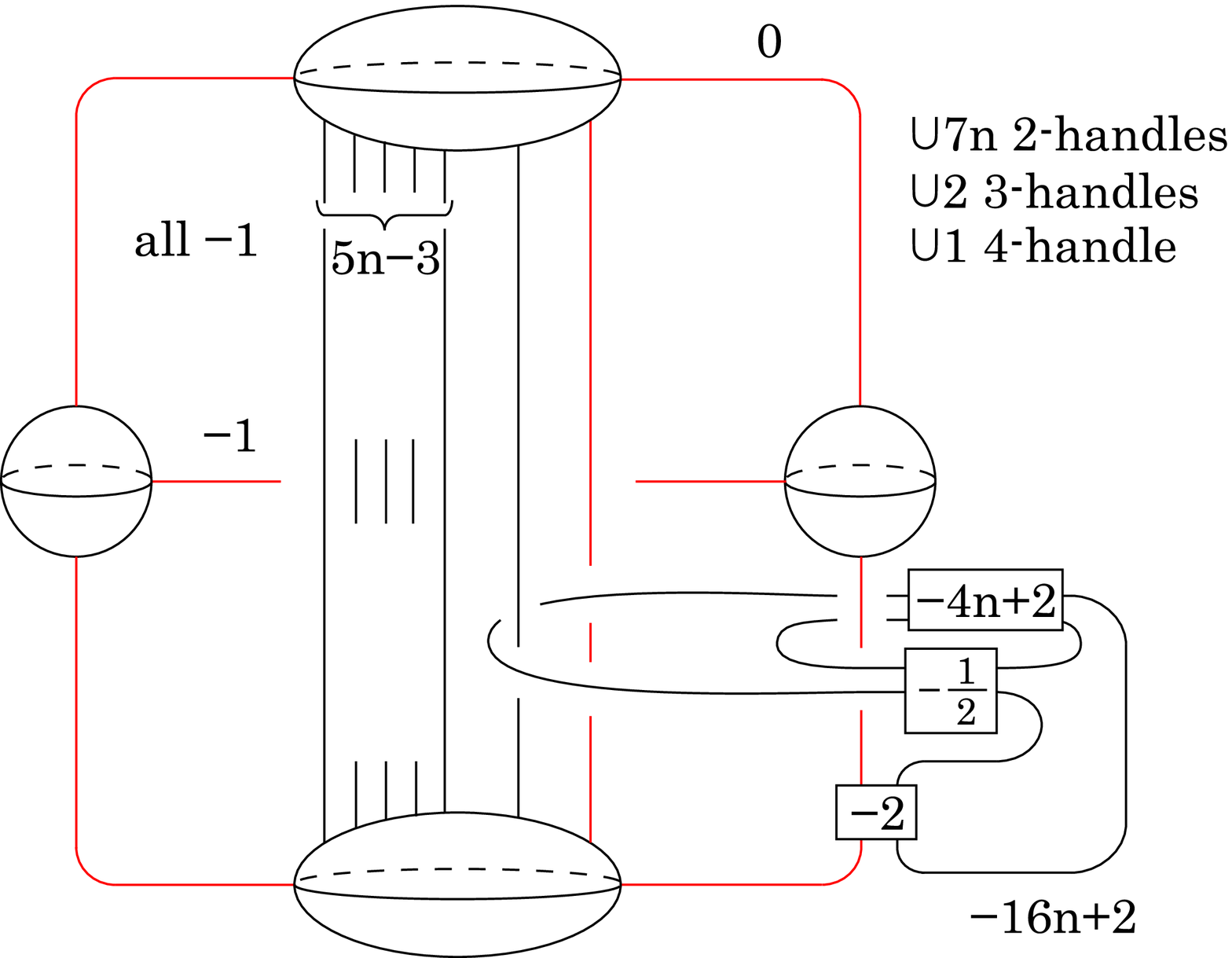}
\caption{$E(n)_2$}
\label{fig18}
\end{center}
\bigskip \medskip

\begin{center}
\includegraphics[width=3.1in]{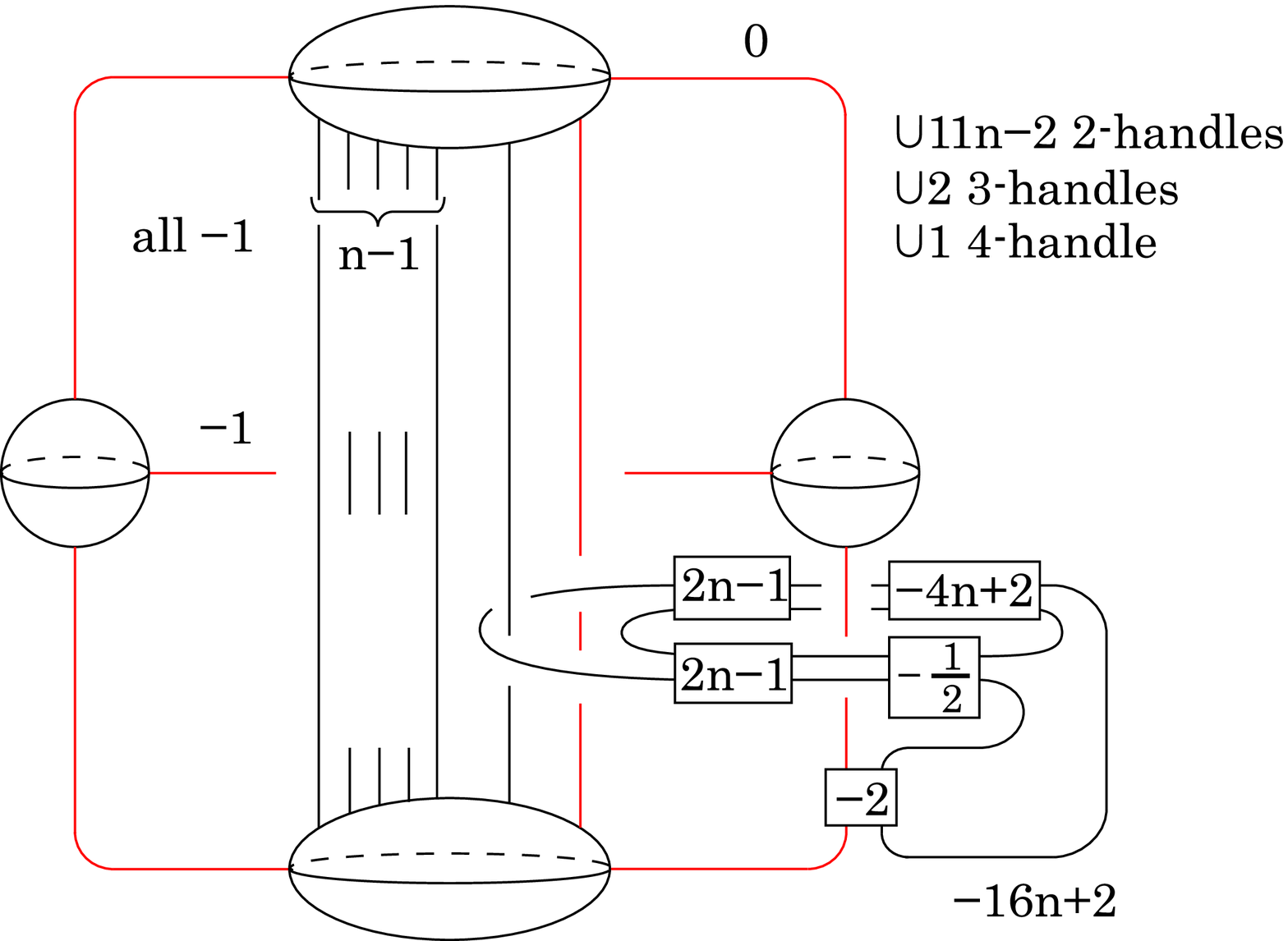}
\caption{$E(n)_2$}
\label{fig19}
\end{center}
\bigskip \medskip

\begin{center}
\includegraphics[width=3.1in]{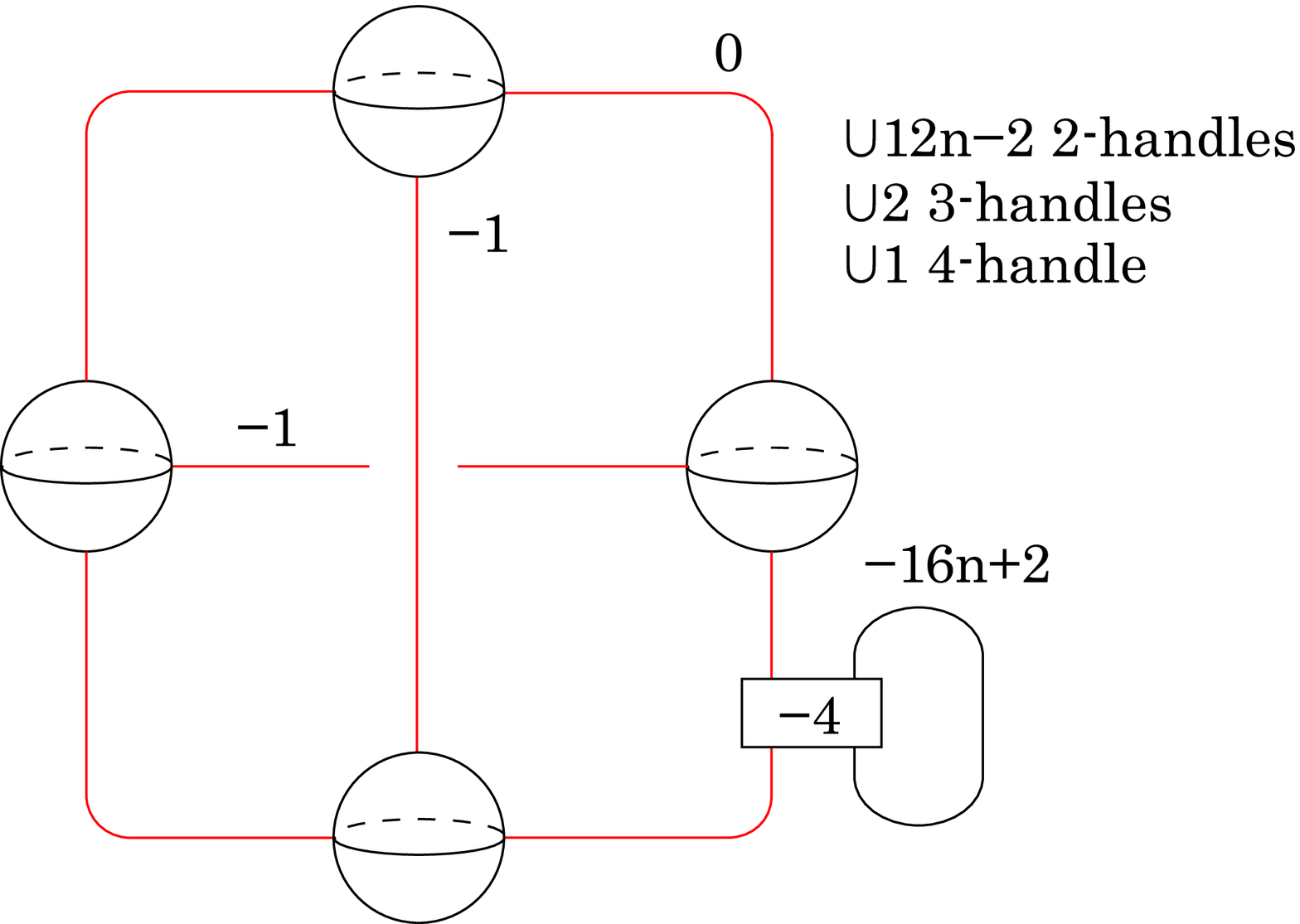}
\caption{$E(n)_2$}
\label{fig20}
\end{center}\end{figure}

\begin{figure}[p]
\begin{center}
\includegraphics[width=3.1in]{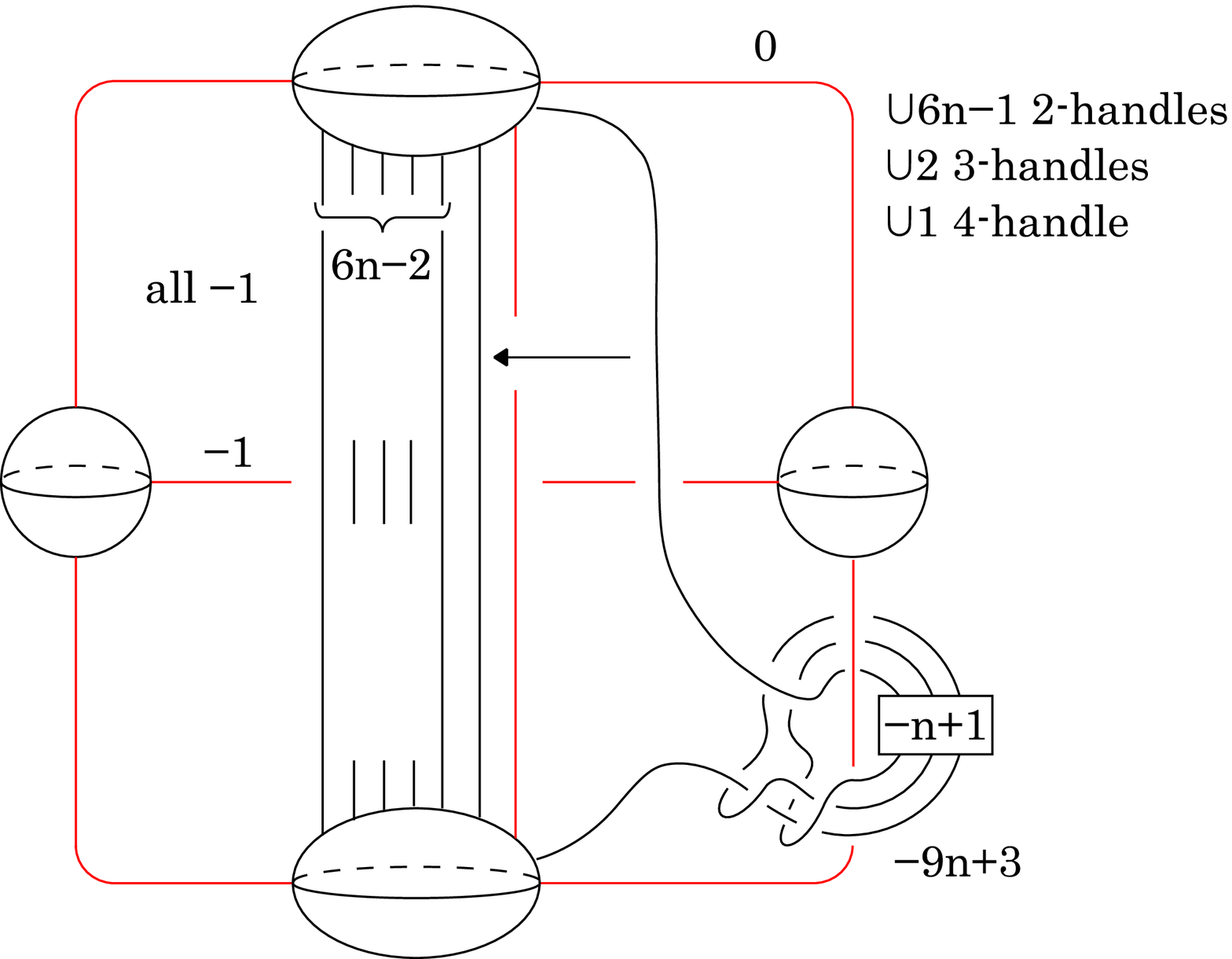}
\caption{$E(n)_3$}
\label{fig21}
\end{center}
\bigskip \medskip

\begin{center}
\includegraphics[width=3.1in]{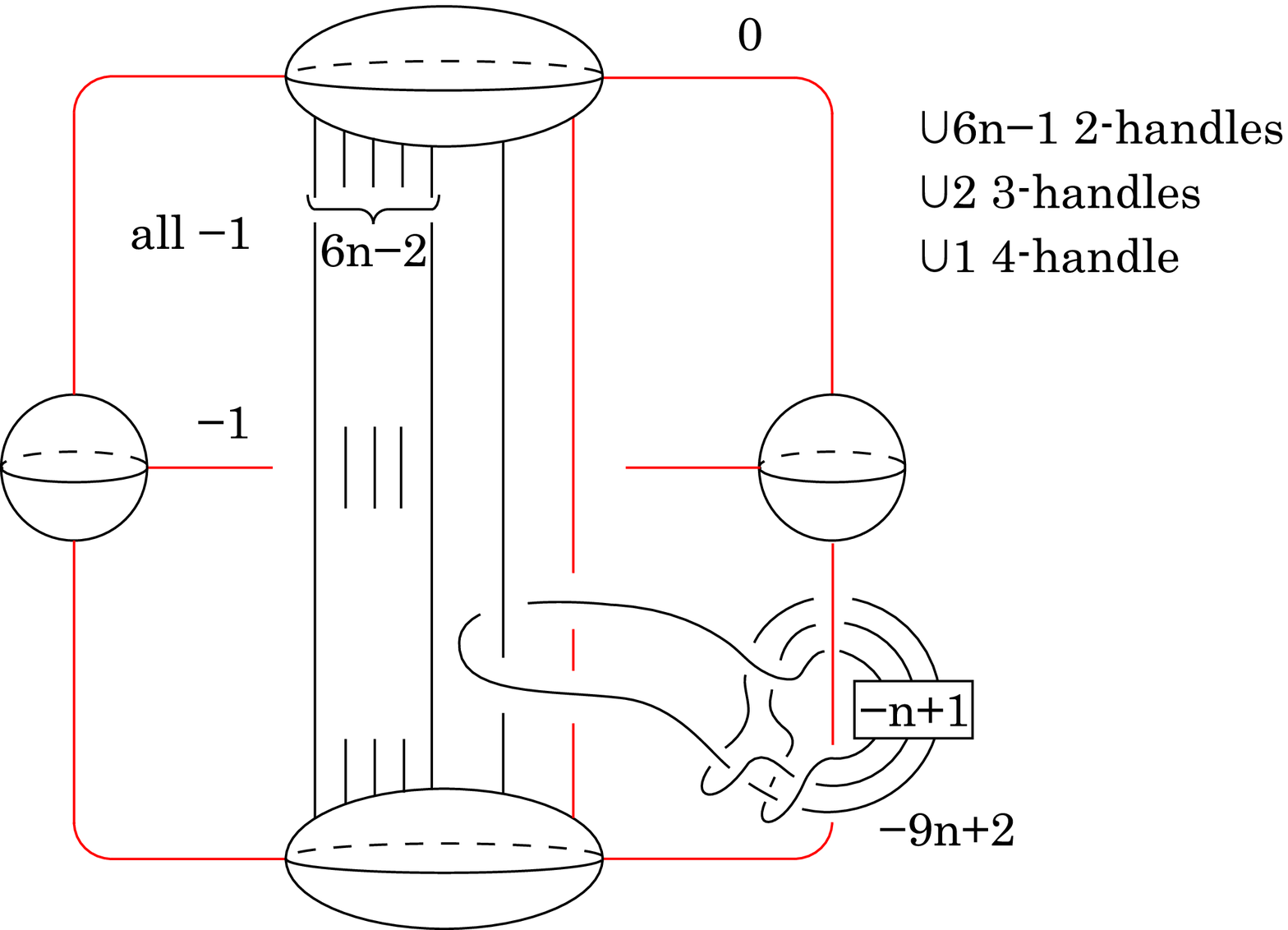}
\caption{$E(n)_3$}
\label{fig22}
\end{center}
\bigskip \medskip

\begin{center}
\includegraphics[width=3.1in]{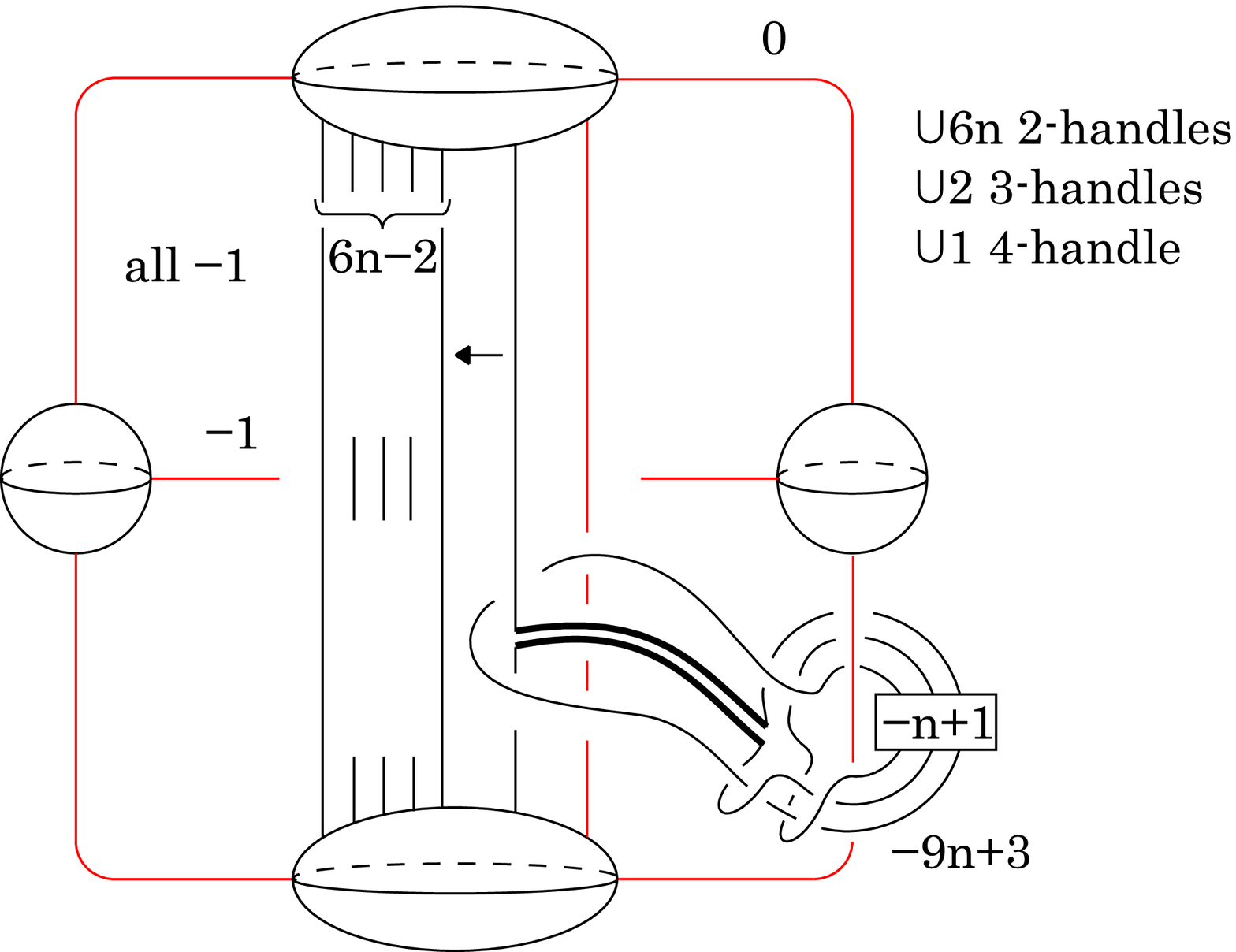}
\caption{$E(n)_3$}
\label{fig23}
\end{center}
\end{figure}

\begin{figure}[p]
\begin{center}
\includegraphics[width=3.1in]{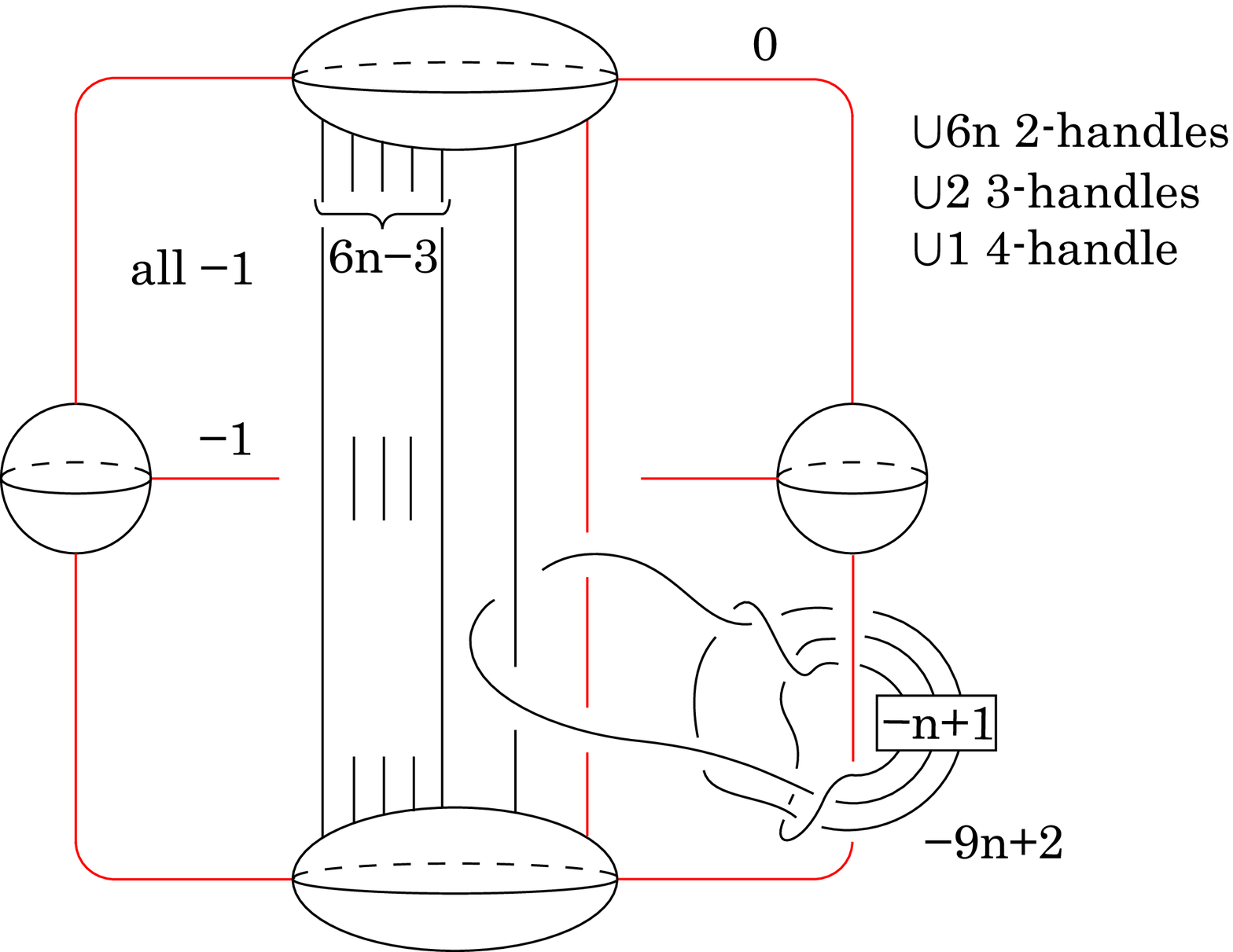}
\caption{$E(n)_3$}
\label{fig24}
\end{center}
\bigskip \medskip

\begin{center}
\includegraphics[width=3.1in]{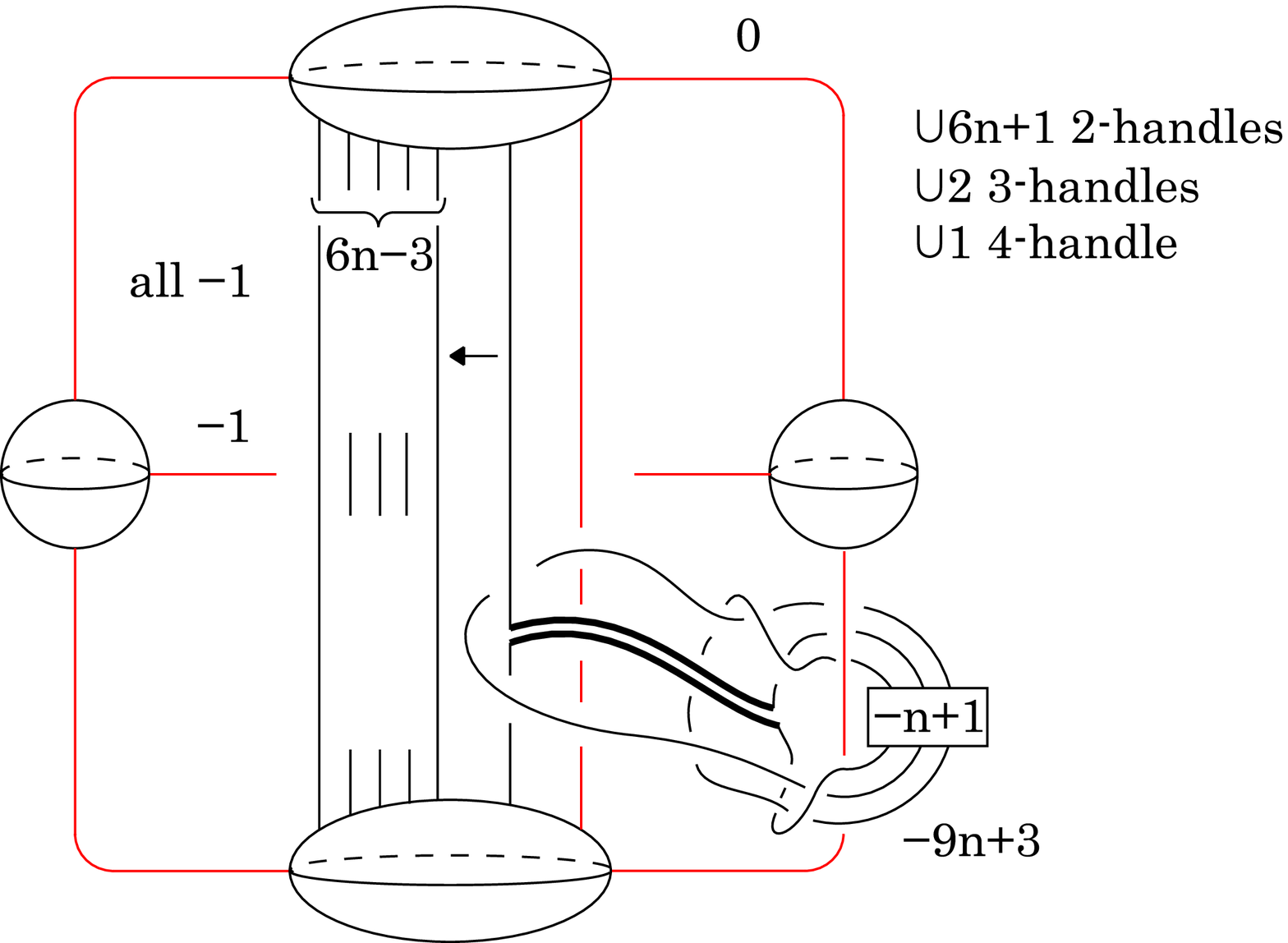}
\caption{$E(n)_3$}
\label{fig25}
\end{center}\bigskip \medskip

\begin{center}
\includegraphics[width=3.1in]{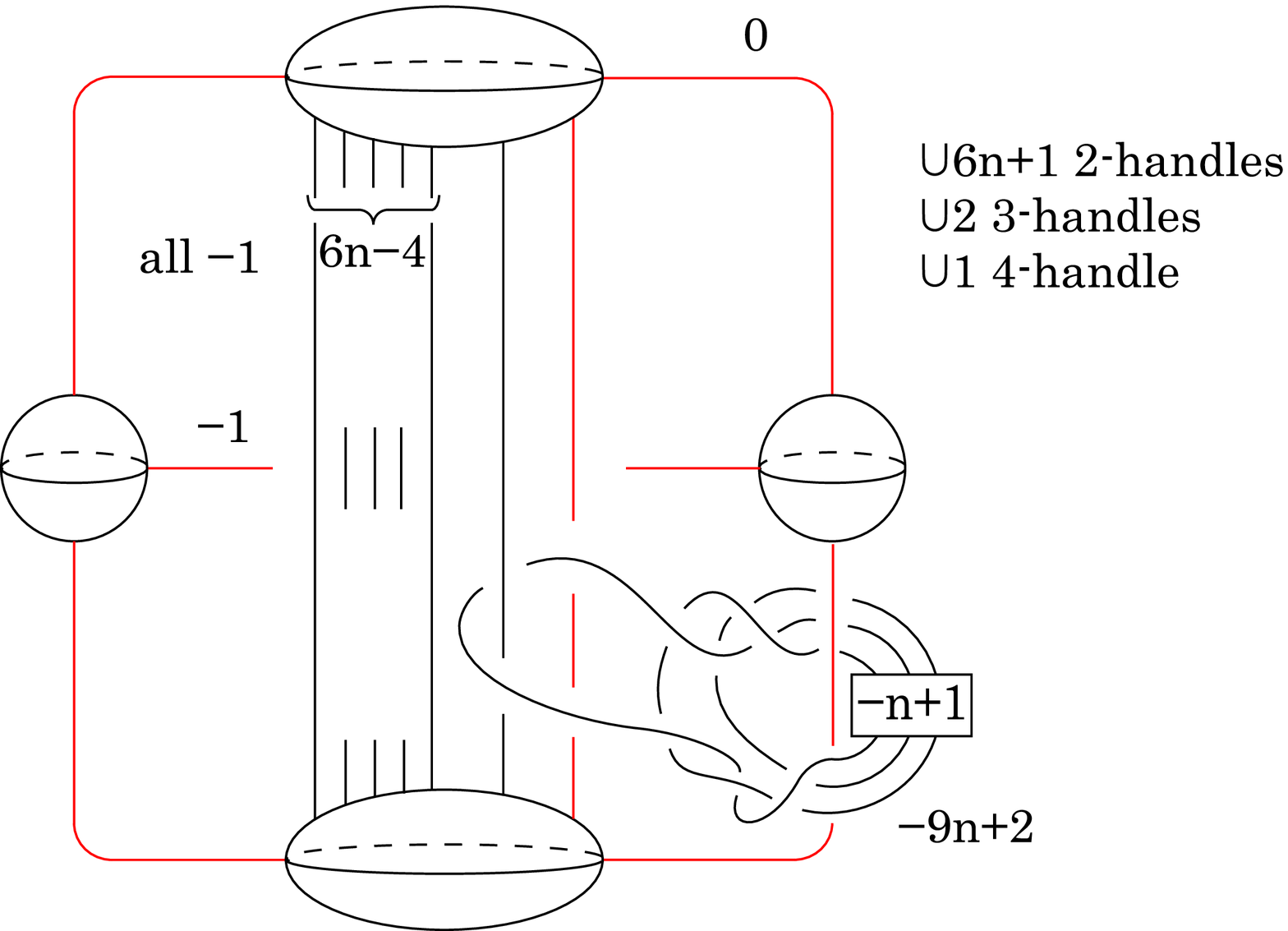}
\caption{$E(n)_3$}
\label{fig26}
\end{center}
\end{figure}

\begin{figure}[p]
\begin{center}
\includegraphics[width=3.1in]{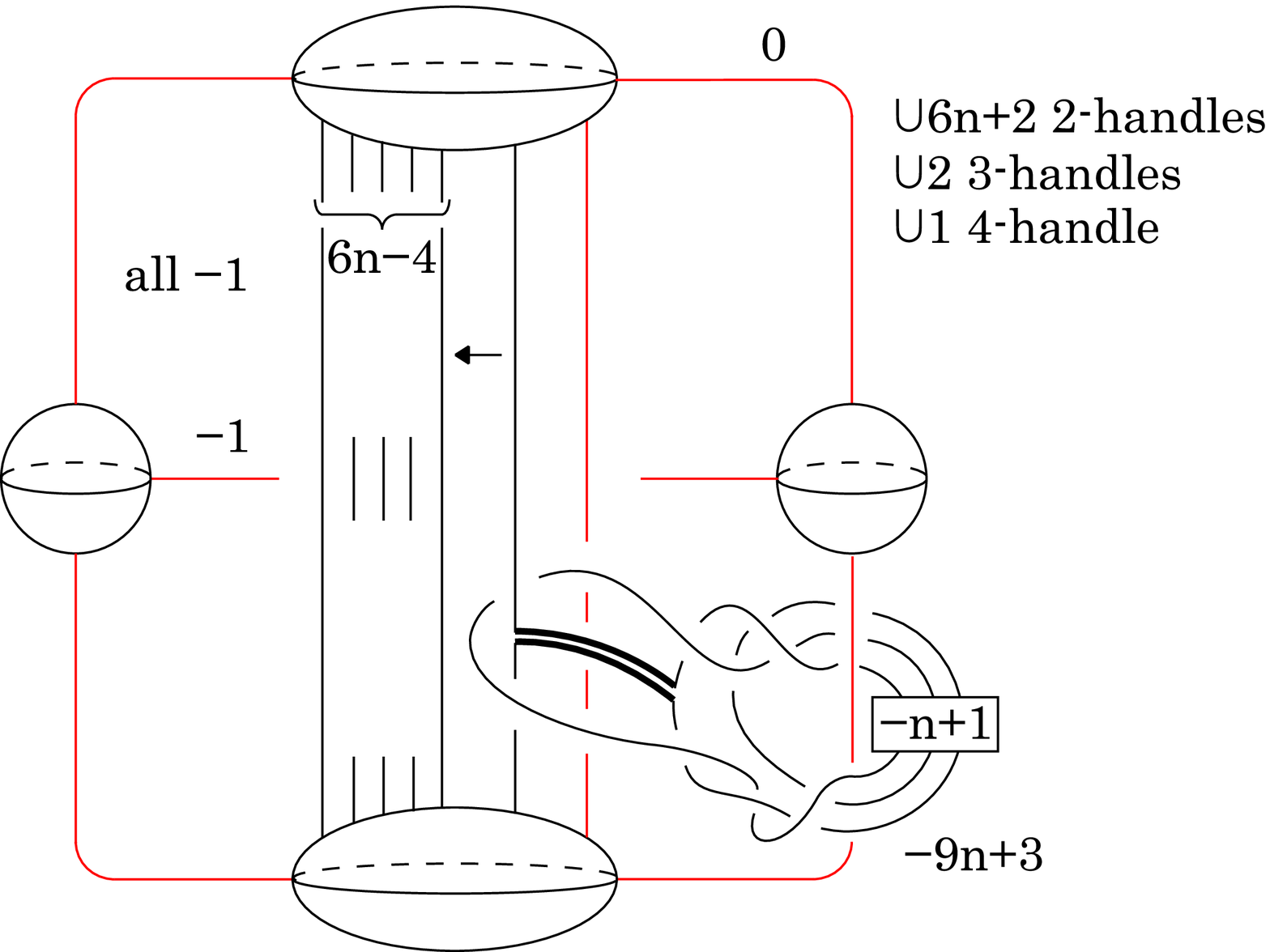}
\caption{$E(n)_3$}
\label{fig27}
\end{center}
\bigskip \medskip

\begin{center}
\includegraphics[width=3.1in]{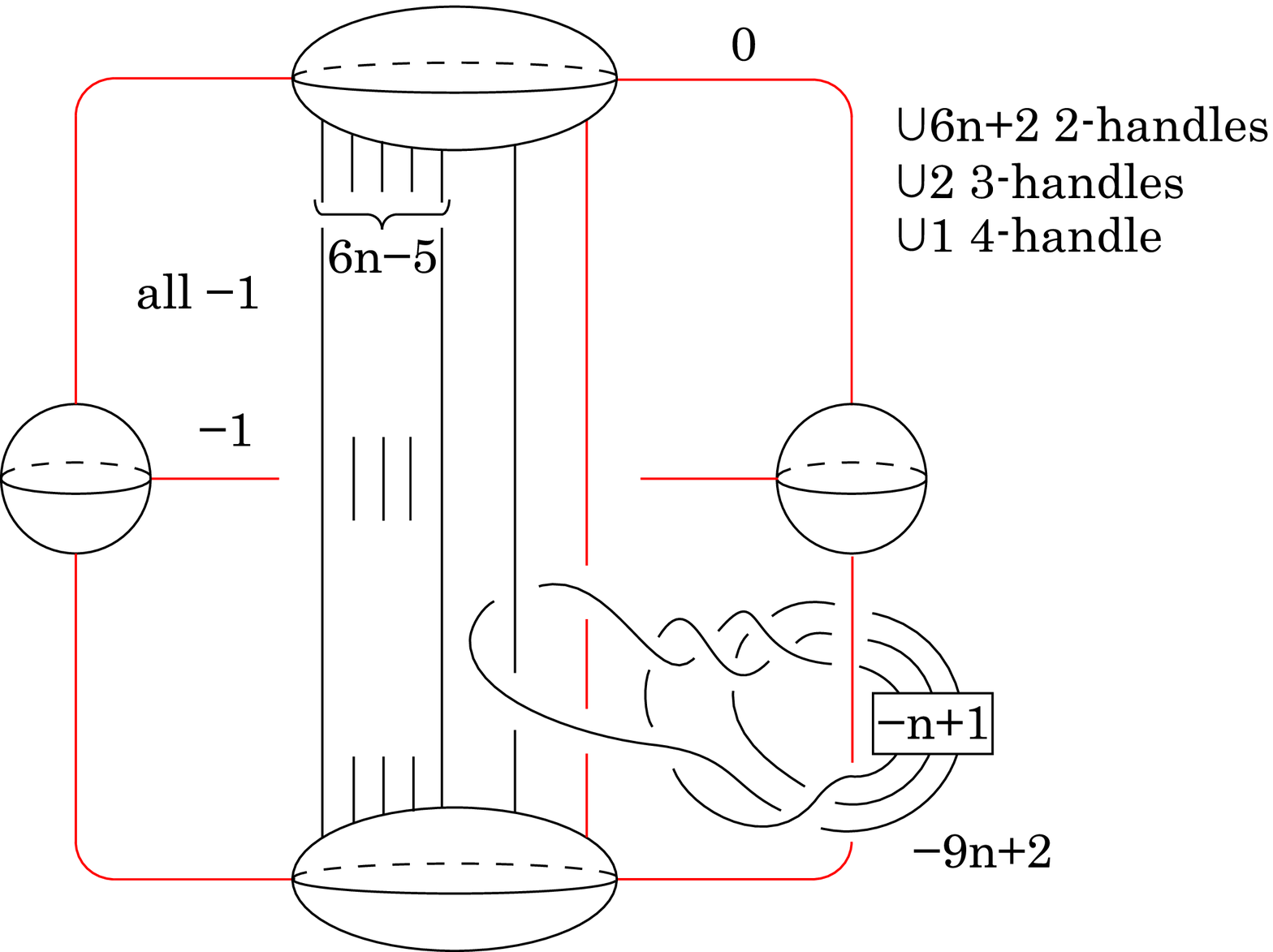}
\caption{$E(n)_3$}
\label{fig28}
\end{center}
\bigskip \medskip

\begin{center}
\includegraphics[width=3.1in]{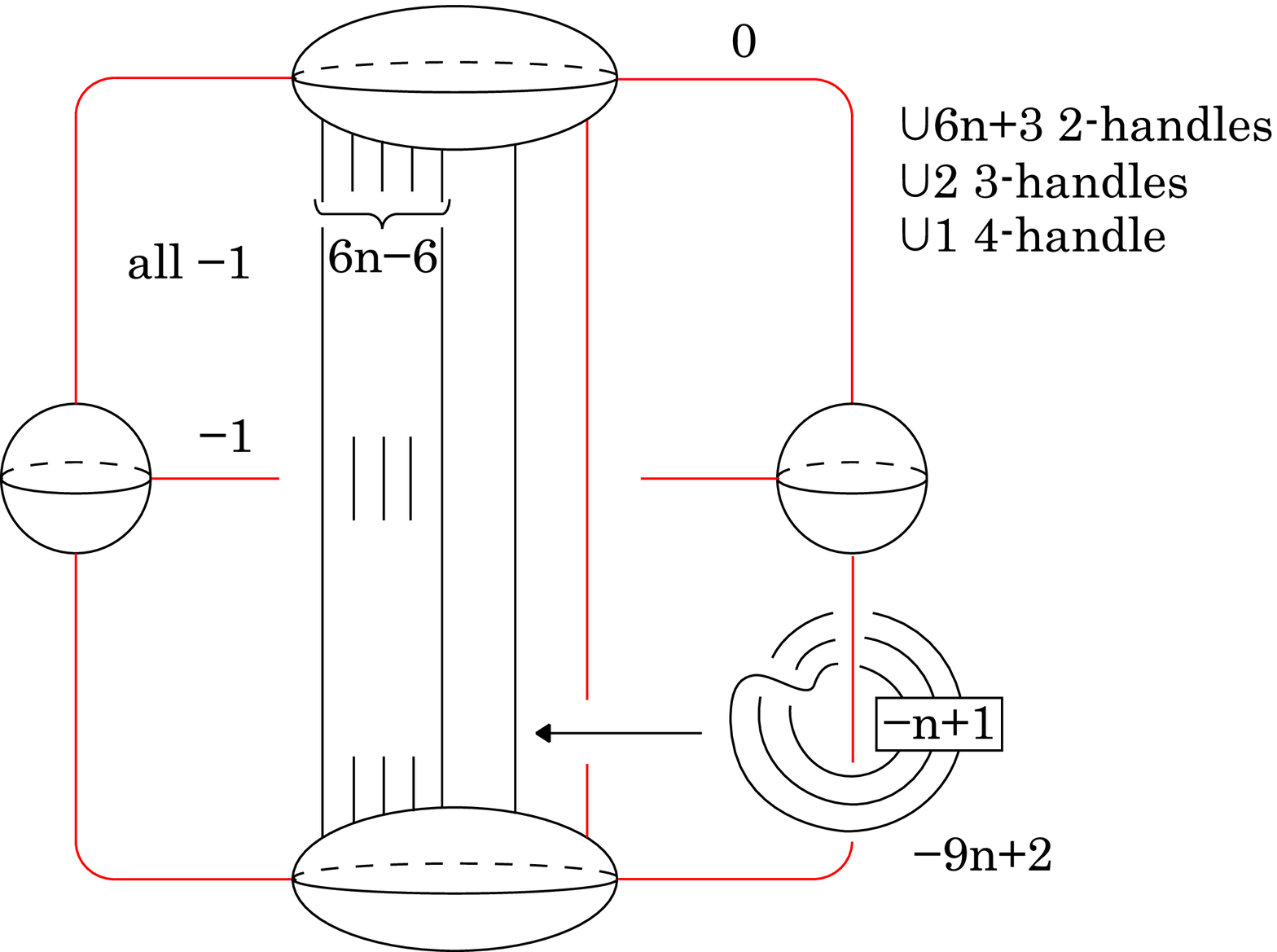}
\caption{$E(n)_3$}
\label{fig29}
\end{center}
\end{figure}

\begin{figure}[p]
\begin{center}
\includegraphics[width=3.1in]{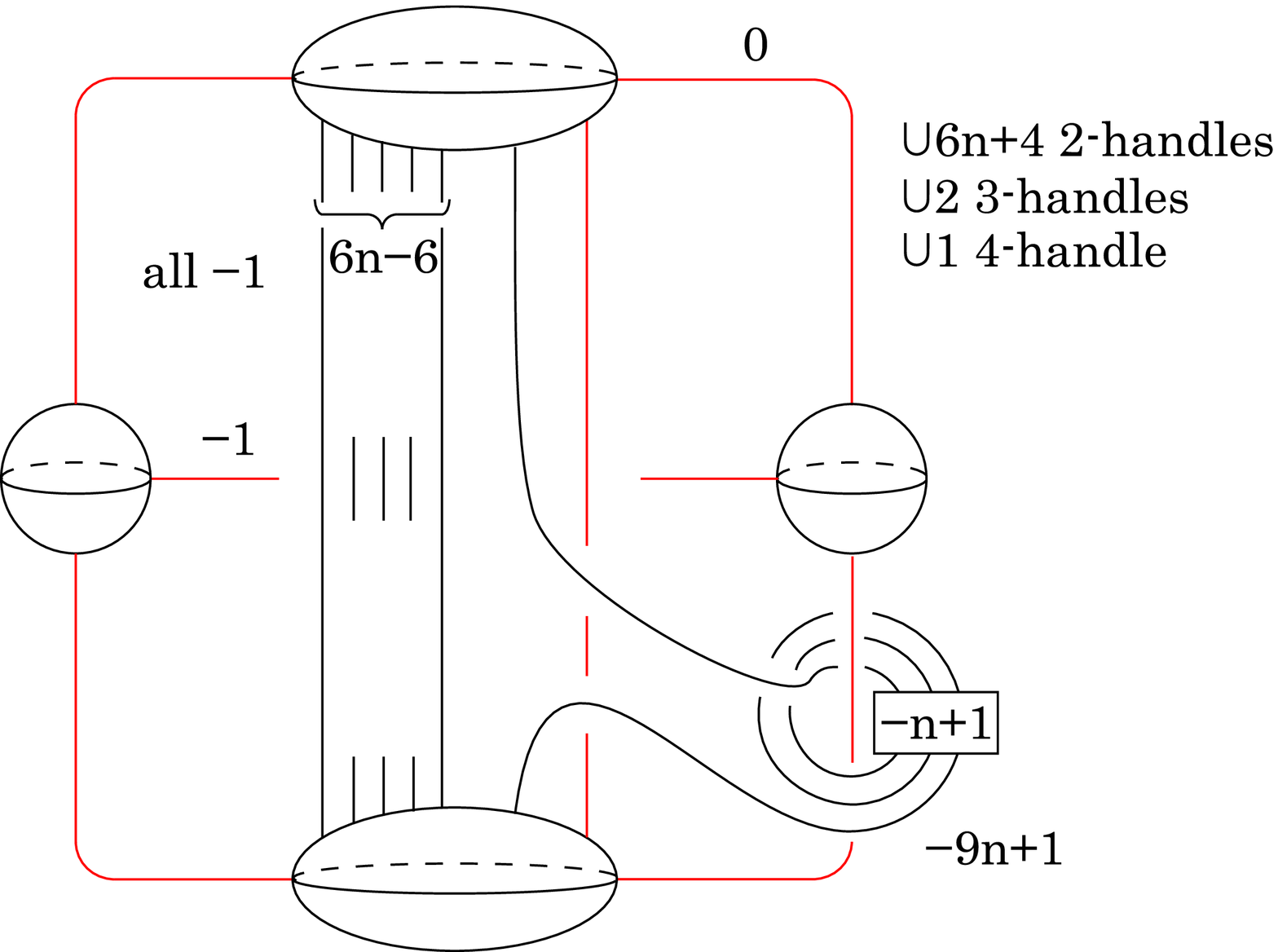}
\caption{$E(n)_3$}
\label{fig30}
\end{center}
\bigskip \medskip

\begin{center}
\includegraphics[width=3.1in]{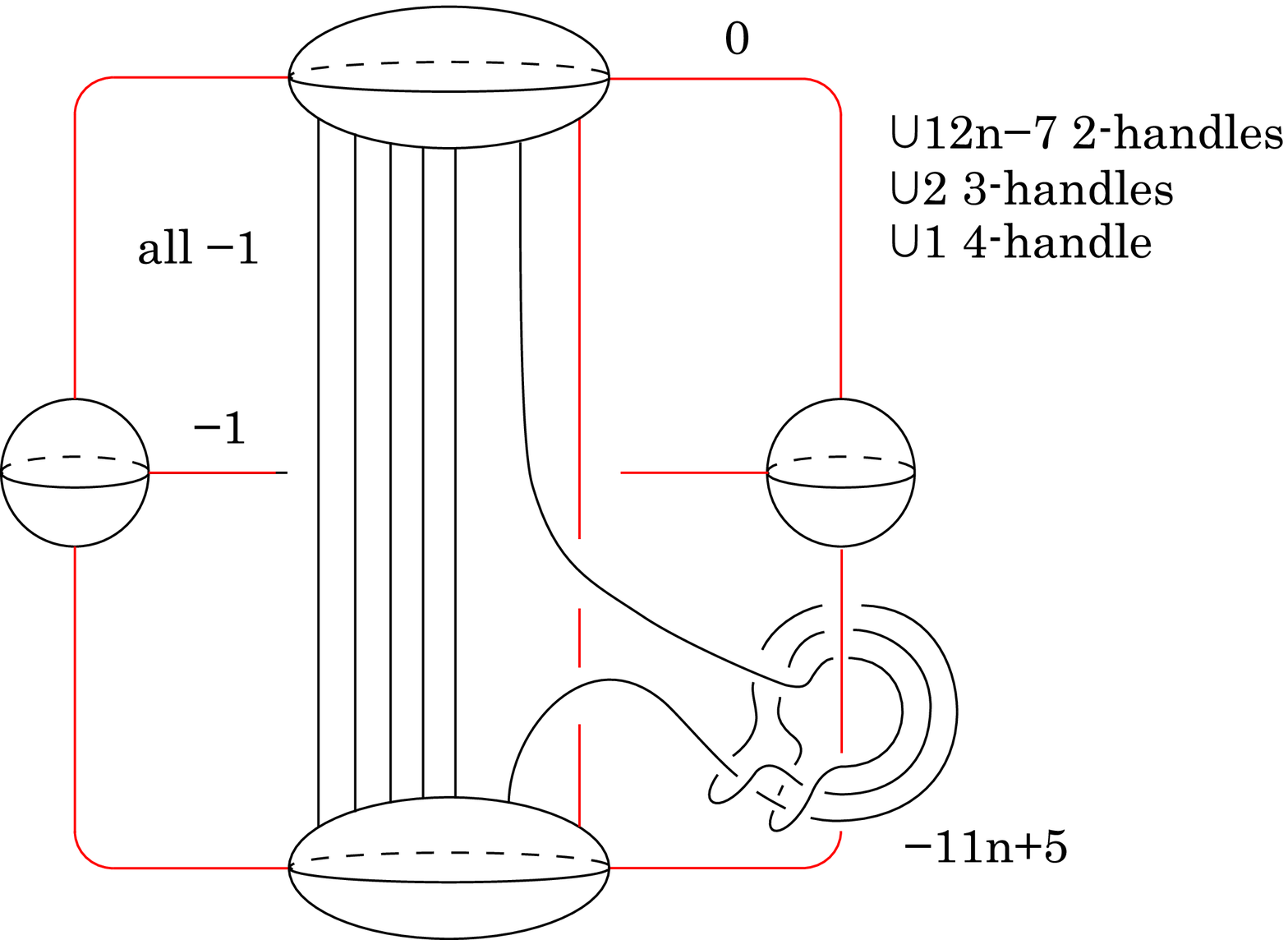}
\caption{$E(n)_3$}
\label{fig31}
\end{center}
\bigskip \medskip

\begin{center}
\includegraphics[width=3.1in]{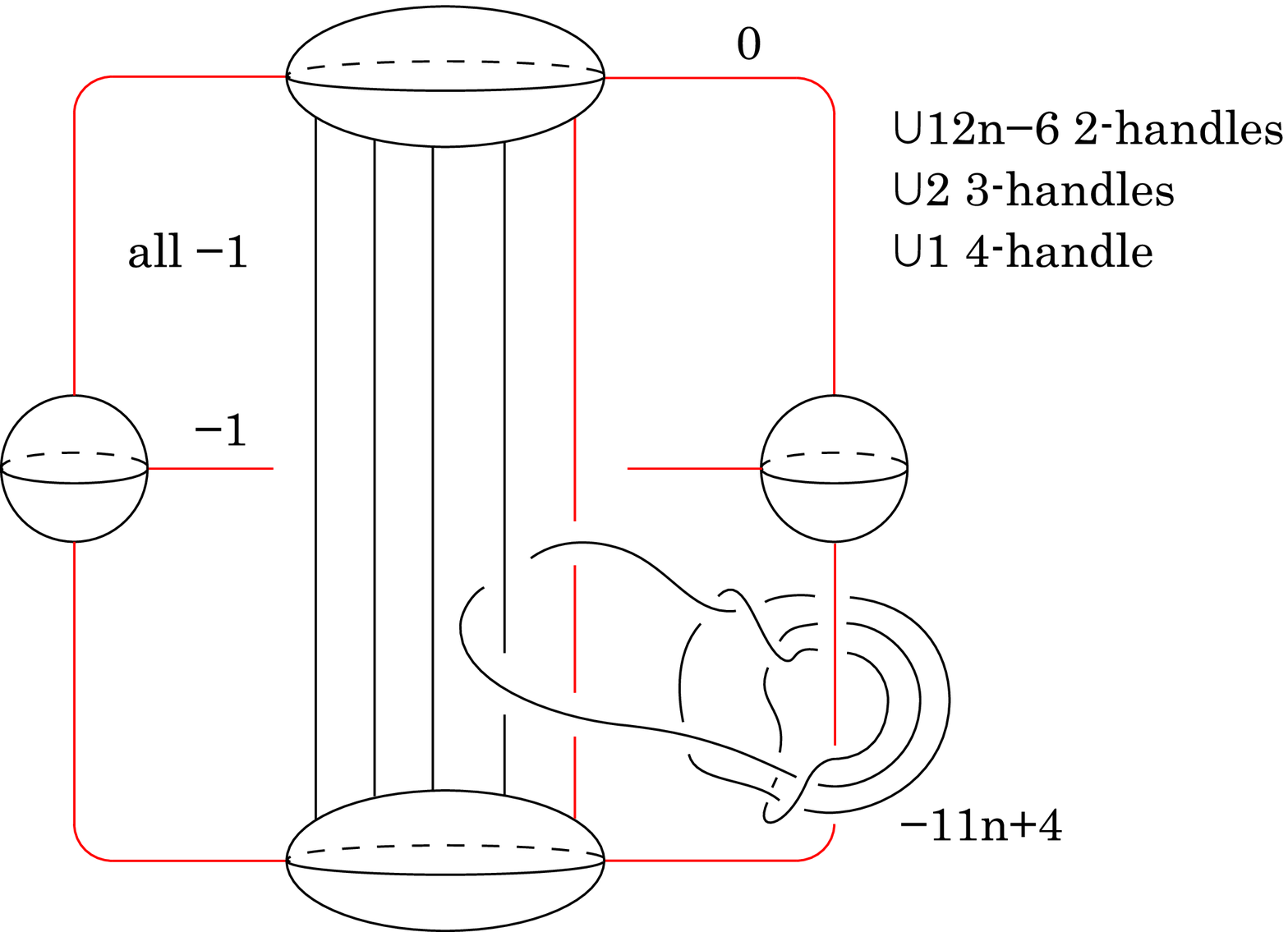}
\caption{$E(n)_3$}
\label{fig32}
\end{center}
\end{figure}

\begin{figure}[p]
\begin{center}
\includegraphics[width=3.1in]{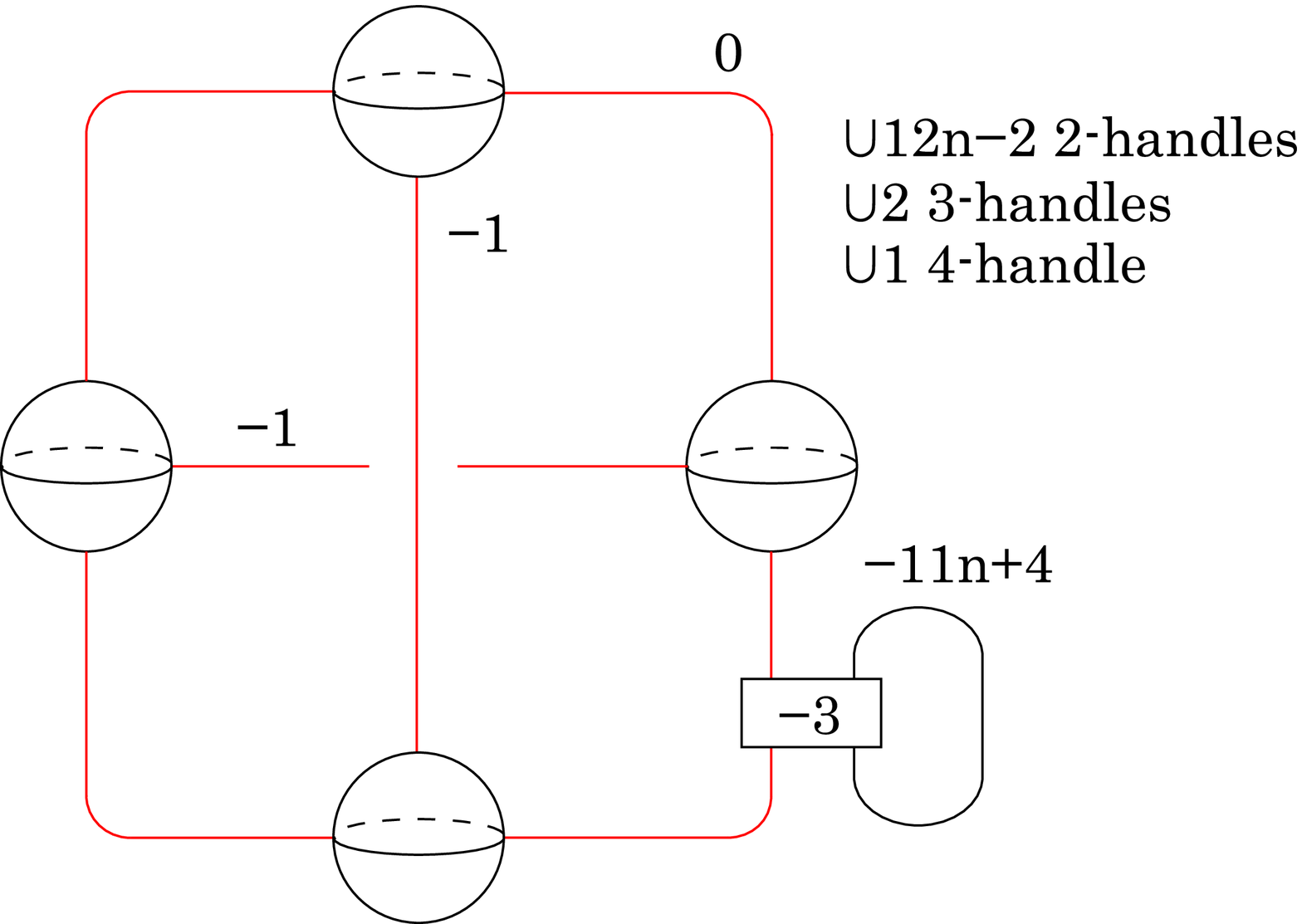}
\caption{$E(n)_3$}
\label{fig33}
\end{center}
\bigskip \medskip \bigskip \medskip \bigskip \medskip

\begin{center}
\includegraphics[width=3.6in]{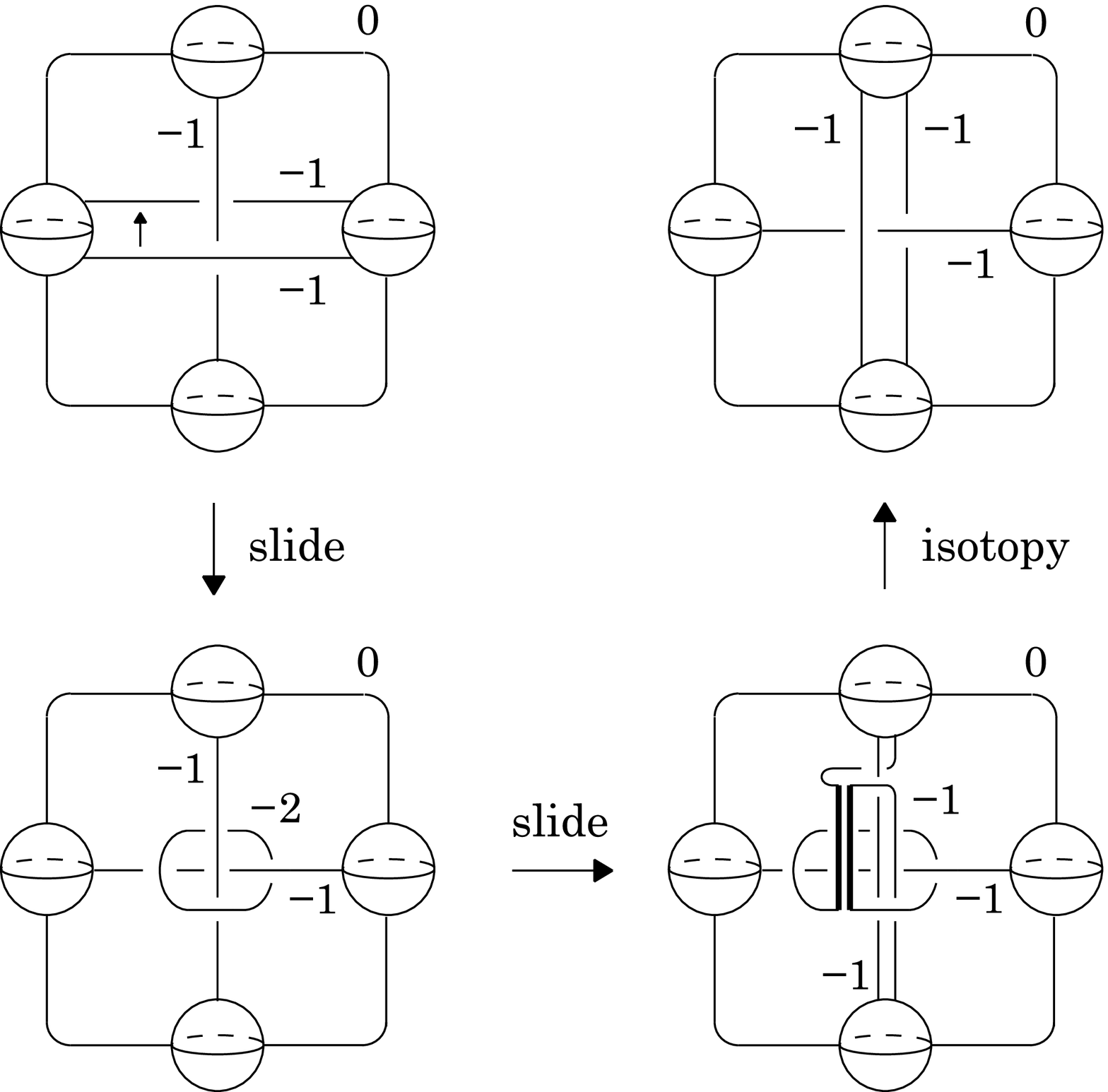}
\caption{handle slides}
\label{fig34}
\end{center}
\end{figure}

\begin{figure}[p]
\begin{center}
\includegraphics[width=3.1in]{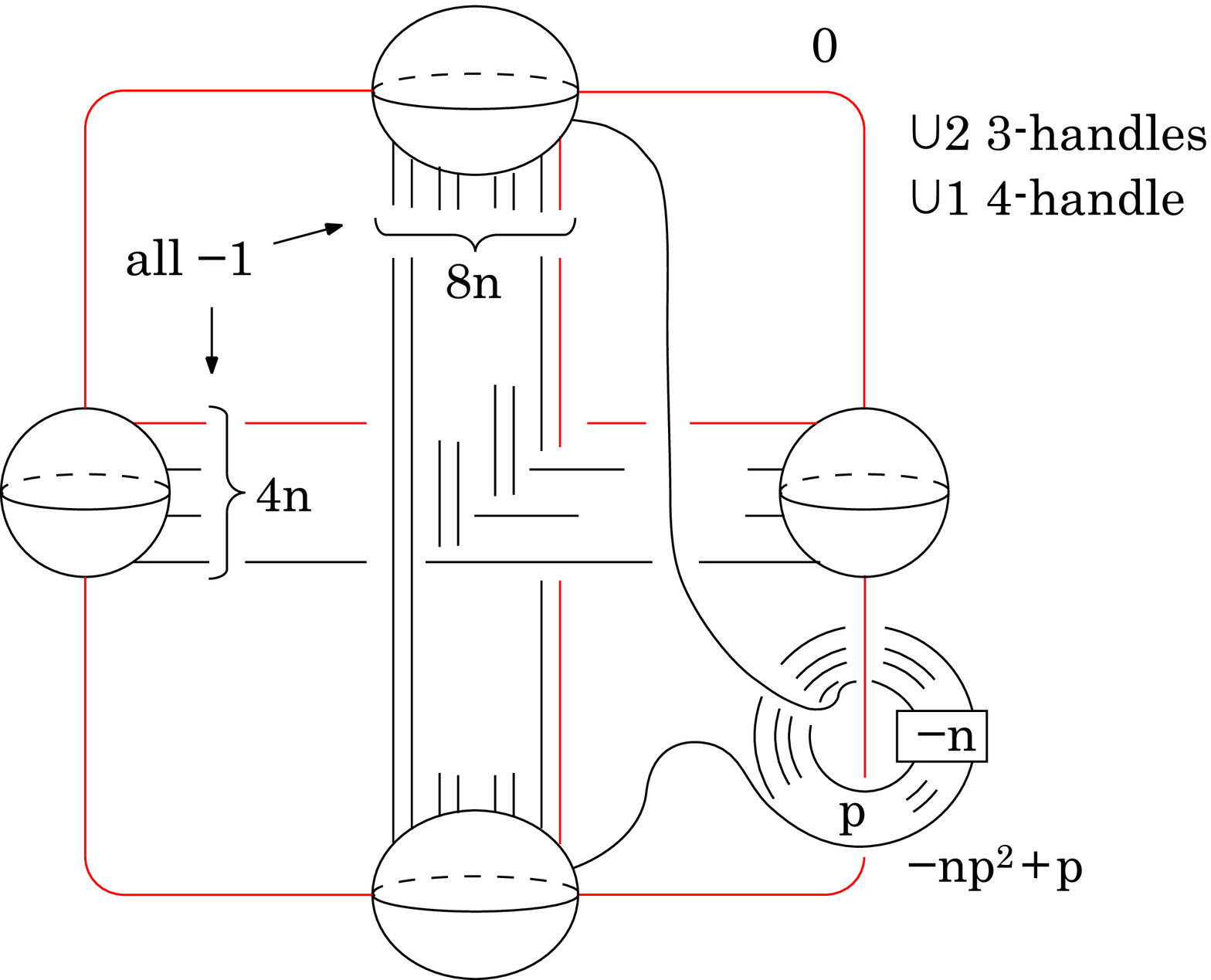}
\caption{$E(n)_p$}
\label{fig35}
\end{center}
\bigskip \medskip

%
\begin{center}
\includegraphics[width=3.2in]{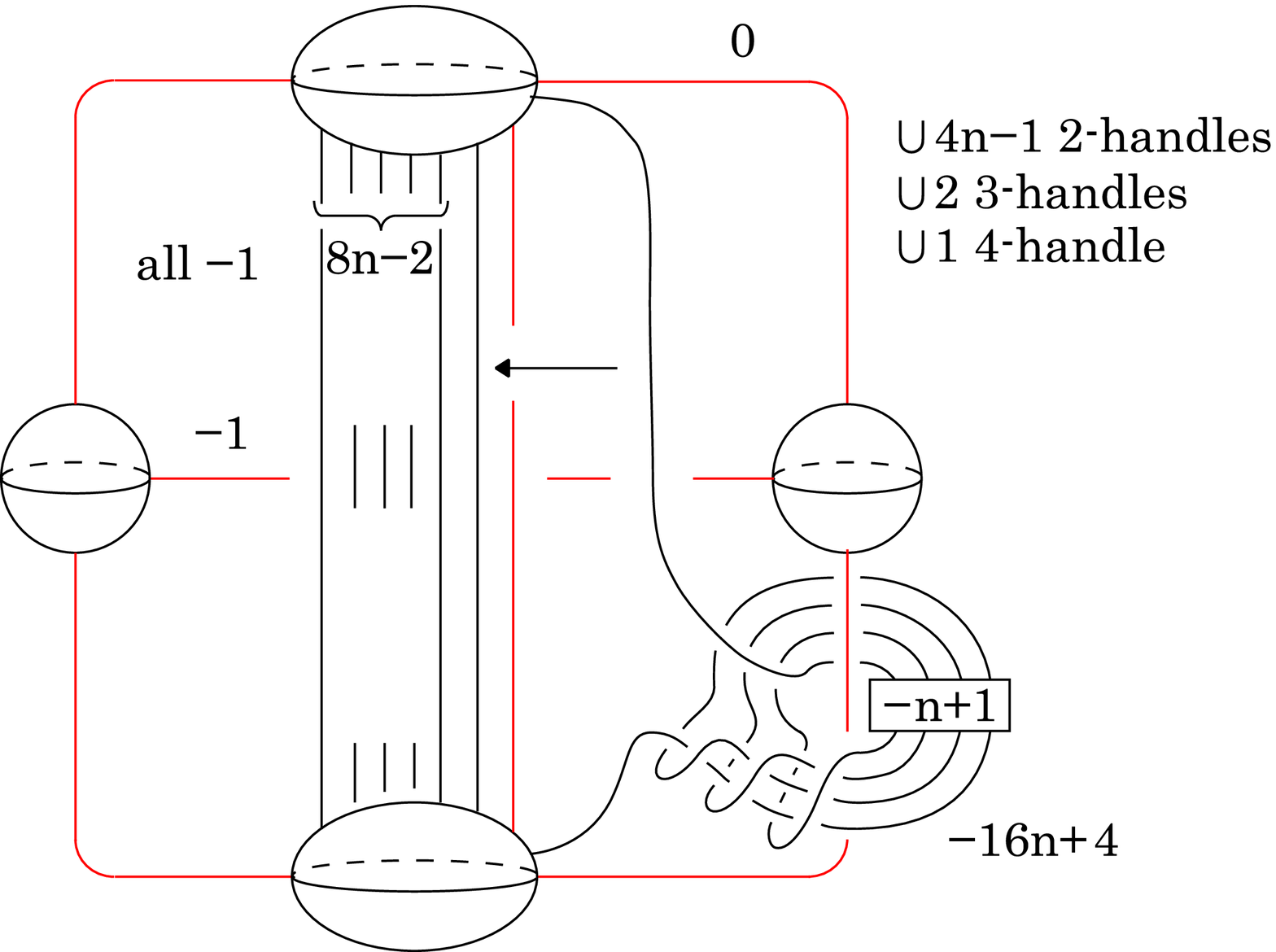}
\caption{$E(n)_4$}
\label{fig36}
\end{center}
\bigskip \medskip

\begin{center}
\includegraphics[width=3.1in]{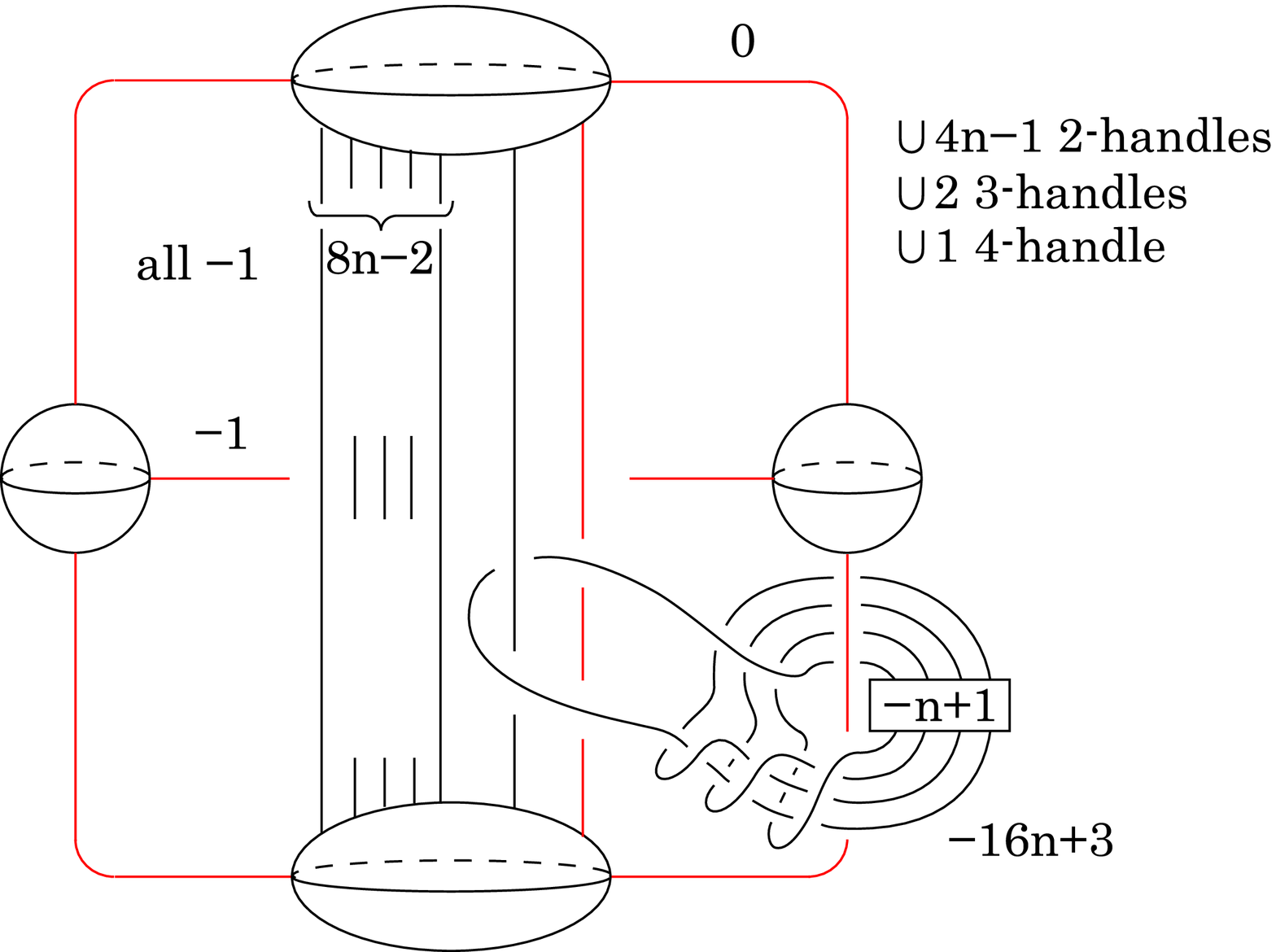}
\caption{$E(n)_4$}
\label{fig37}
\end{center}
\end{figure}

We finish this paper by making some remarks.
\begin{remark}A key of our proof of the main theorem is to elliminate extra twists of a 2-handle of $E(n)_p$ so that we can apply Corollary~\ref{cor:key}. To carry out the key, we used many vertical $-1$ framed $2$-handles of $E(n)_p$ in Figure~\ref{fig9} or \ref{fig35}. Perhaps, we may obtain more examples of elliptic surfaces without 1-handles by additionally using horizontal $2$-handles of $E(n)_p$ in Figure~\ref{fig9} or \ref{fig35}.  
\end{remark}
\begin{remark}
In \cite{Y1} and \cite{Y2}, we constructed a smooth $4$-manifold $E_3'$ which is homeomorphic to $E(1)_{2,3}$. The $4$-manifold $E_3'$ has the same Seiberg-Witten invariant as $E(1)_{2,3}$ and has a handle decomposition without $1$- and $3$-handles.  $E_3'$ is constructed from $\mathbf{C}\mathbf{P}^2\# 13\overline{\mathbf{C}\mathbf{P}^2}$ by rationally blowing down $C_5$. However, it is not known whether or not $E(1)_{2,3}$ can be obtained from $\mathbf{C}\mathbf{P}^2\# 13\overline{\mathbf{C}\mathbf{P}^2}$ by rationally blowing down $C_5$. 
We do not know whether or not manifolds in \cite{Y2} are diffeomorphic to $E(1)_{2,q}$ $(q=3,5)$.
\end{remark}
\begin{remark}
It seems more interesting to investigate handle decompositions of exotic $4$-manifolds with small Euler characteristics, because there exist no exotic $S^4$ and no exotic $\mathbf{CP}^2$ which admit handle decompositions without 1- and 3-handles. 
In \cite{Y3}, we constructed exotic $\mathbf{CP}^2\# n\overline{\mathbf{C}\mathbf{P}^2}$ $(5\leq n\leq 9)$ which admit neither 1- nor 3-handles for $7\leq n\leq 9$. 
\end{remark}

%
%
%
\begin{figure}[p]
\begin{center}
\includegraphics[width=3.1in]{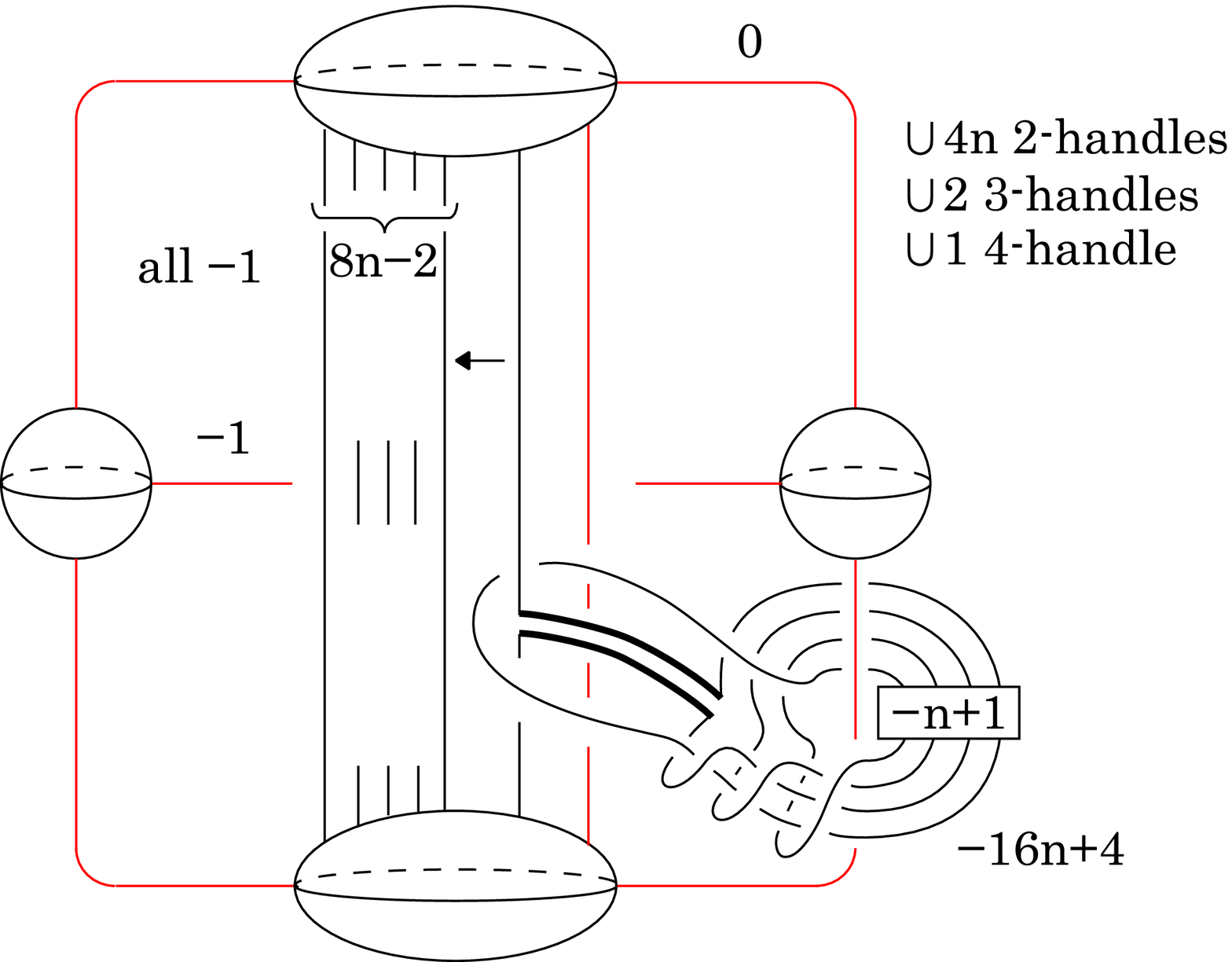}
\caption{$E(n)_4$}
\label{fig38}
\end{center}
\bigskip \medskip

\begin{center}
\includegraphics[width=3.1in]{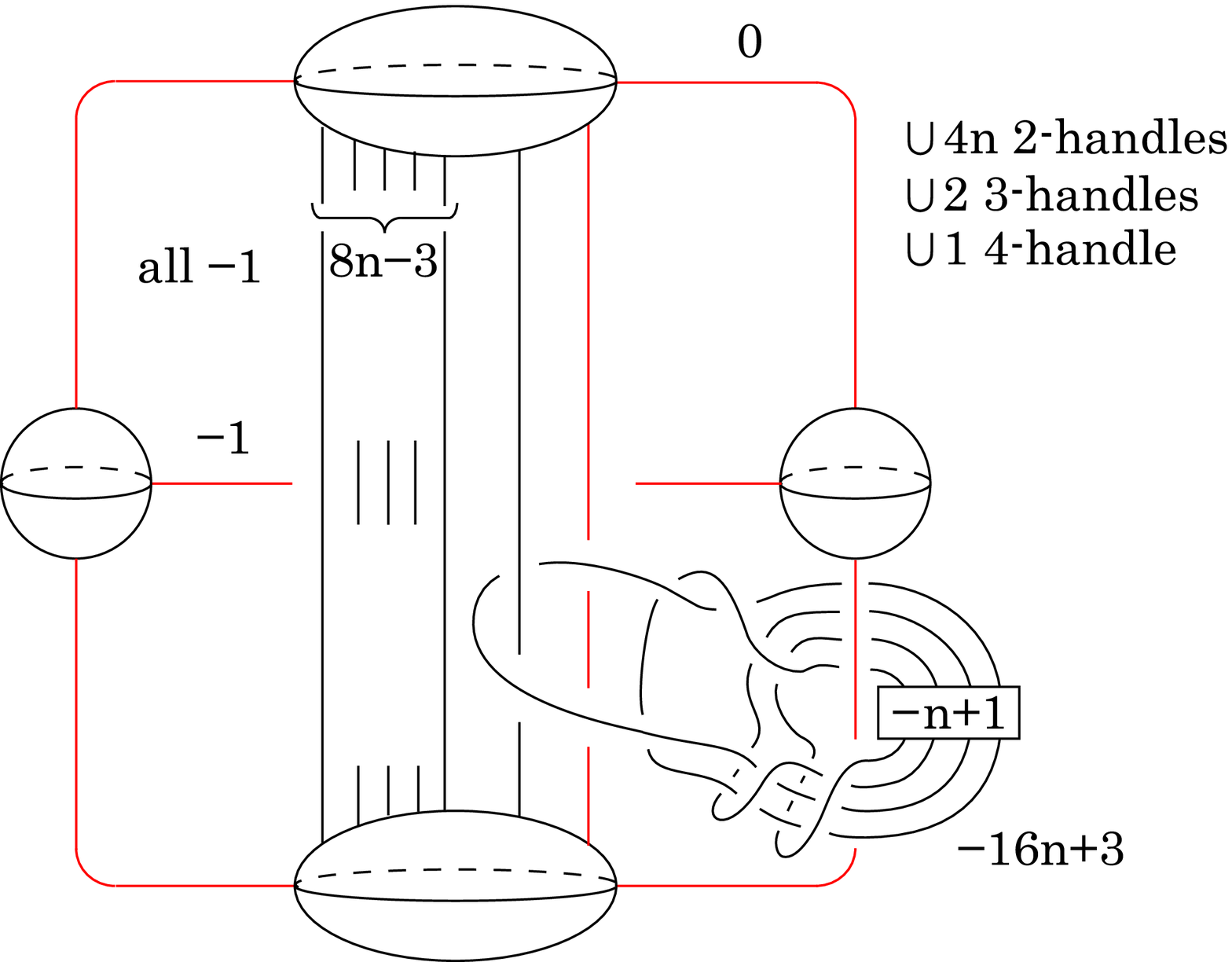}
\caption{$E(n)_4$}
\label{fig39}
\end{center}
\bigskip \medskip

\begin{center}
\includegraphics[width=3.1in]{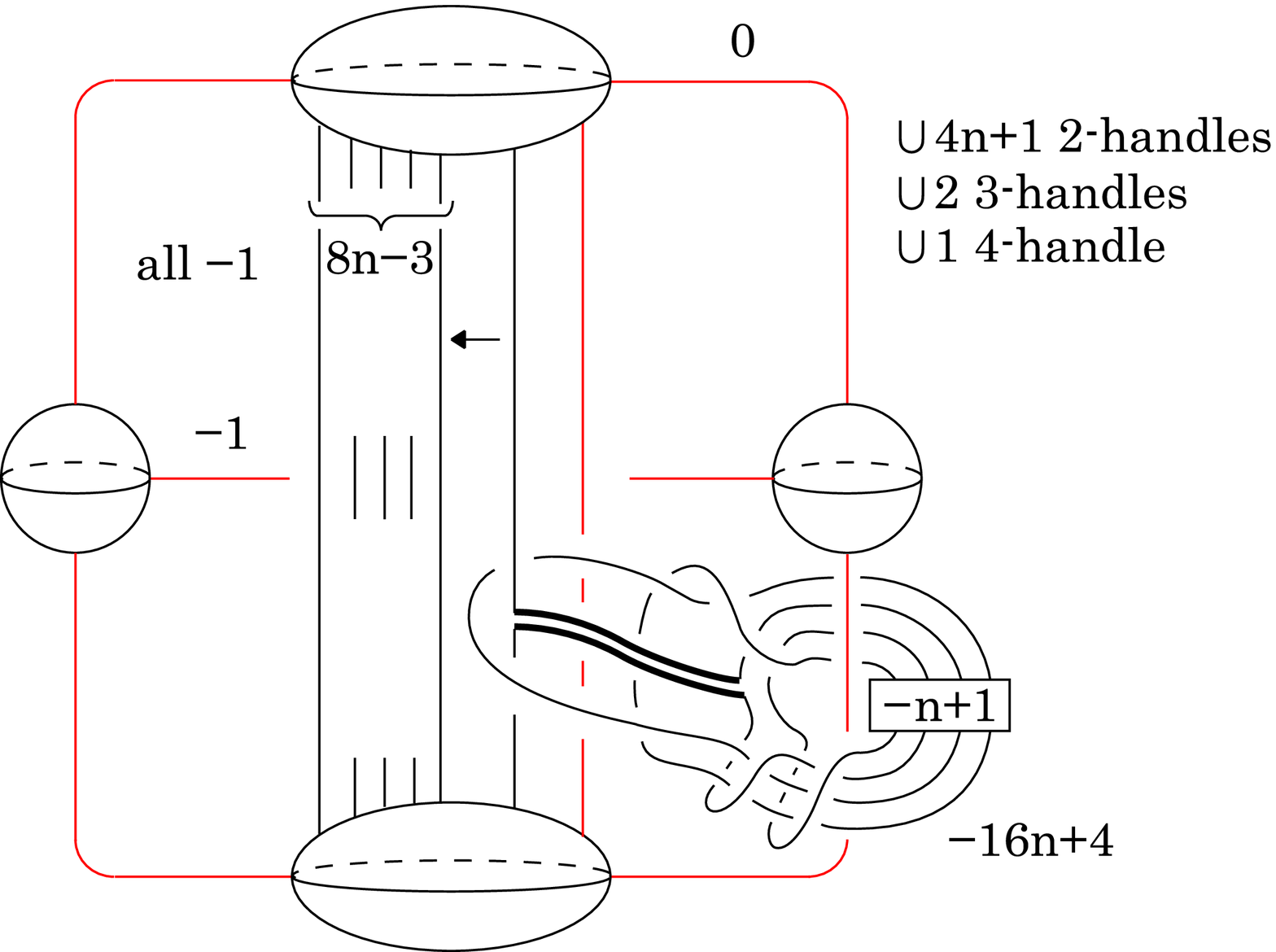}
\caption{$E(n)_4$}
\label{fig40}
\end{center}
\end{figure}
%
%
\begin{figure}[p]
\begin{center}
\includegraphics[width=3.1in]{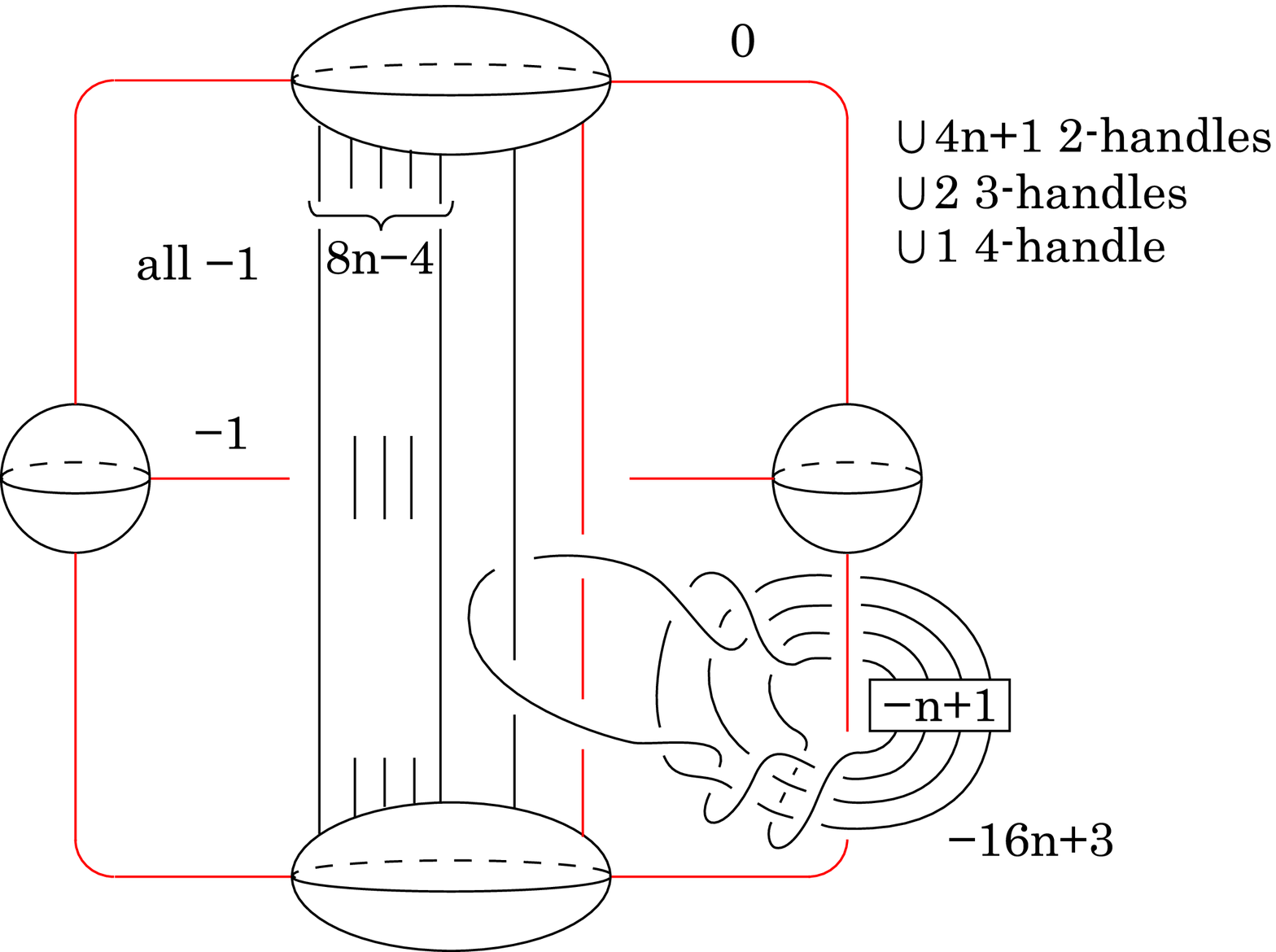}
\caption{$E(n)_4$}
\label{fig41}
\end{center}
\bigskip \medskip

\begin{center}
\includegraphics[width=3.1in]{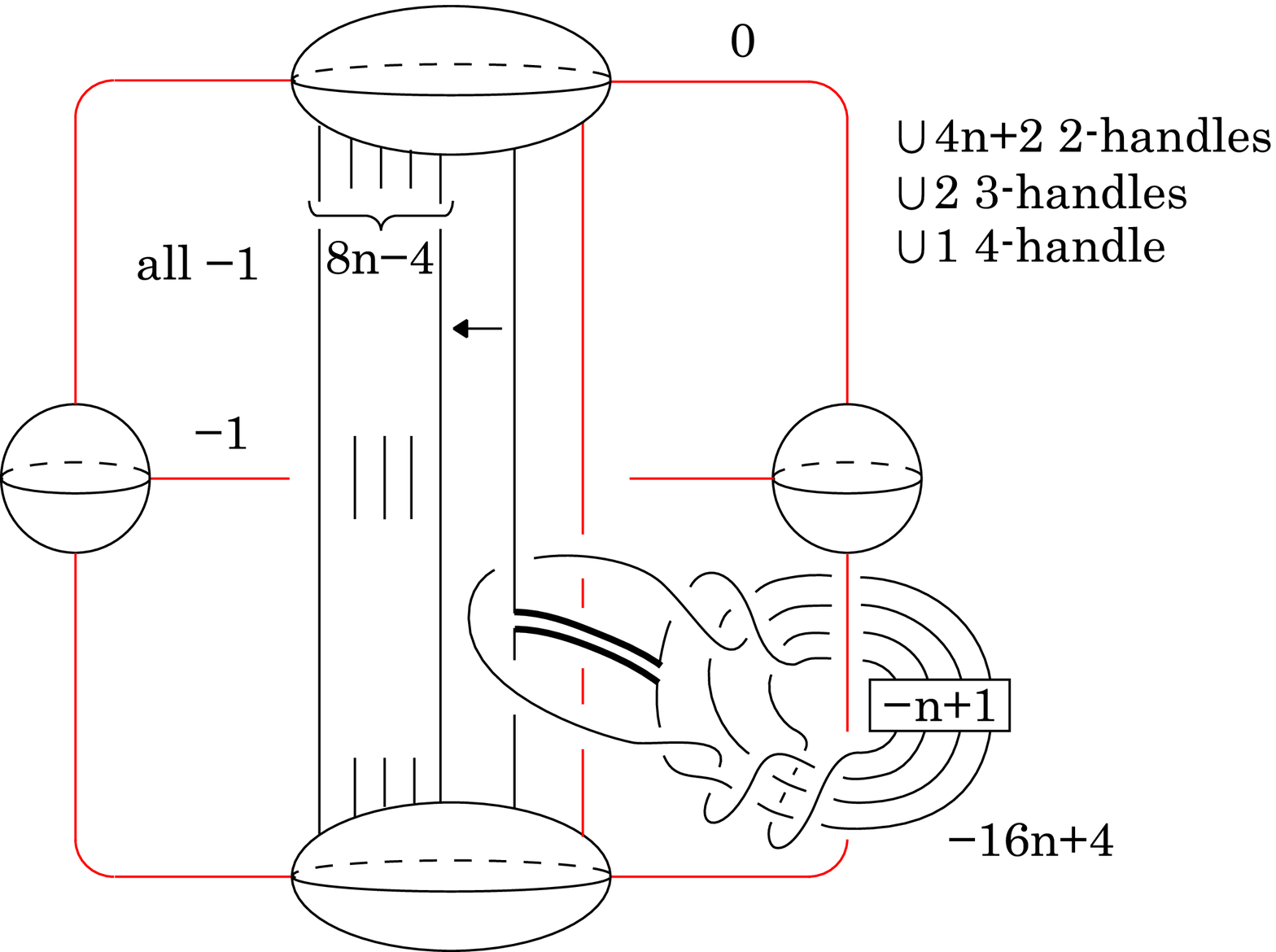}
\caption{$E(n)_4$}
\label{fig42}
\end{center}
\bigskip \medskip

\begin{center}
\includegraphics[width=3.1in]{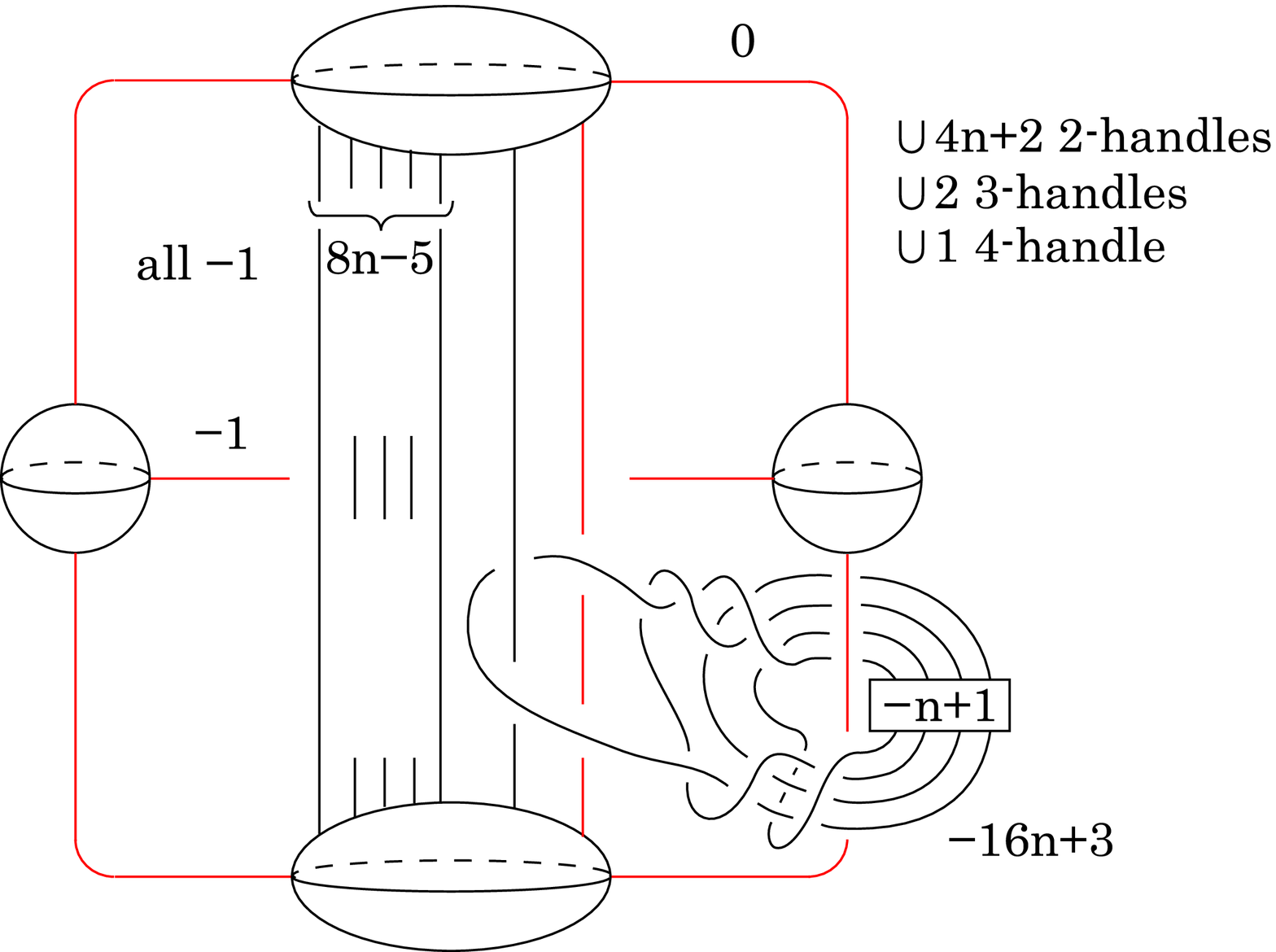}
\caption{$E(n)_4$}
\label{fig43}
\end{center}
\end{figure}
\begin{figure}[p]
\begin{center}
\includegraphics[width=3.1in]{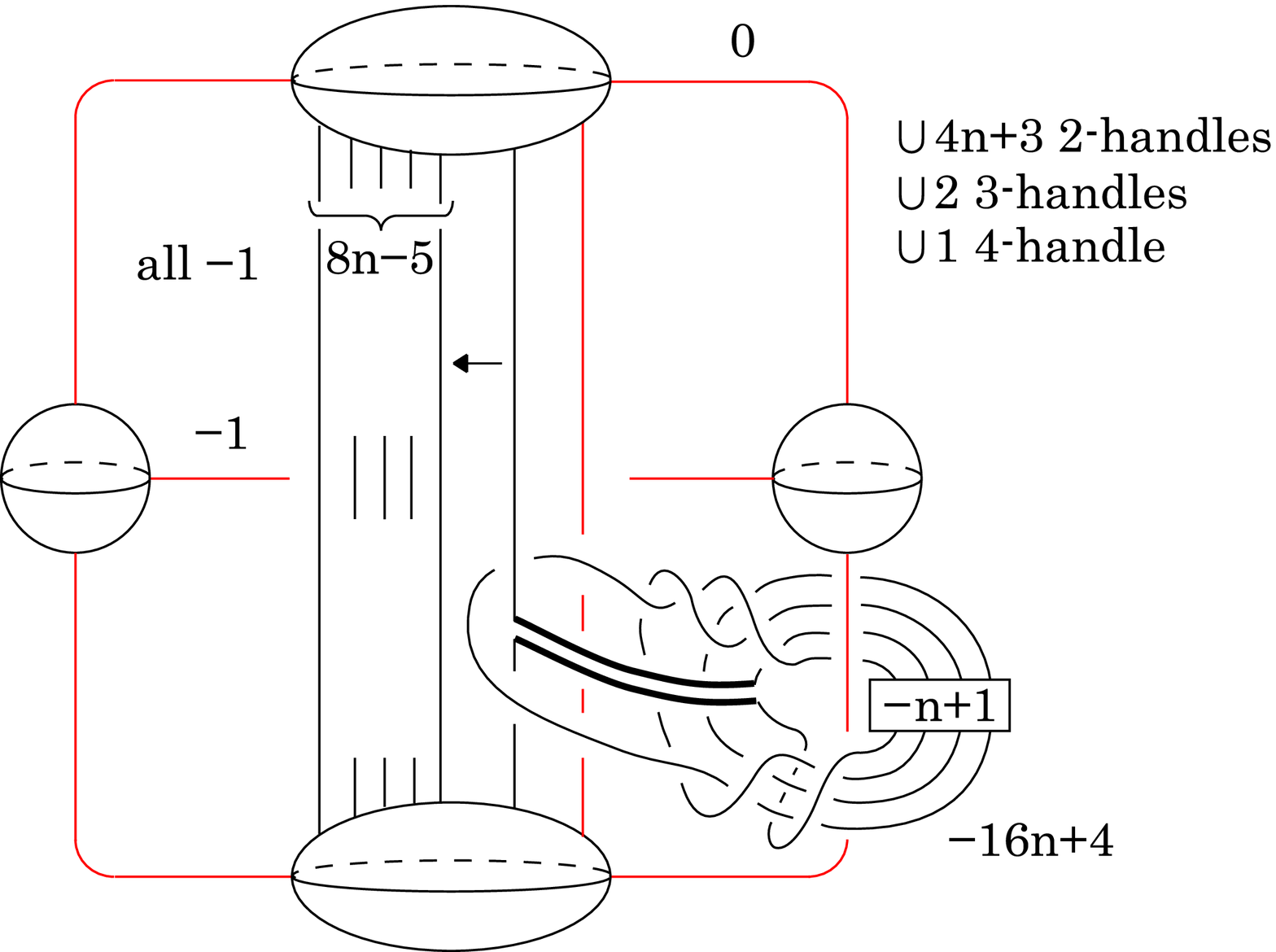}
\caption{$E(n)_4$}
\label{fig44}
\end{center}
\bigskip \medskip

\begin{center}
\includegraphics[width=3.1in]{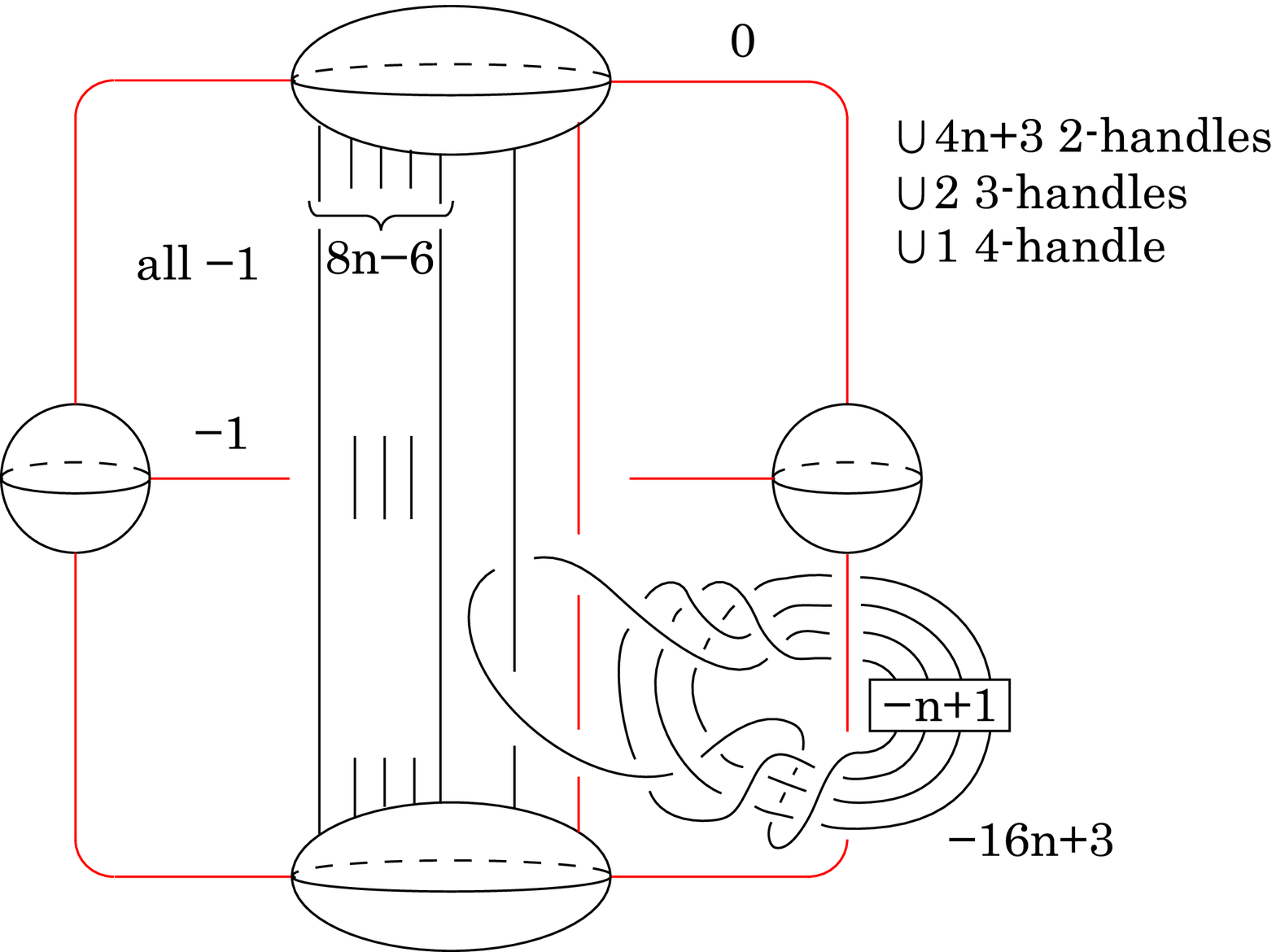}
\caption{$E(n)_4$}
\label{fig45}
\end{center}
\bigskip \medskip

\begin{center}
\includegraphics[width=3.1in]{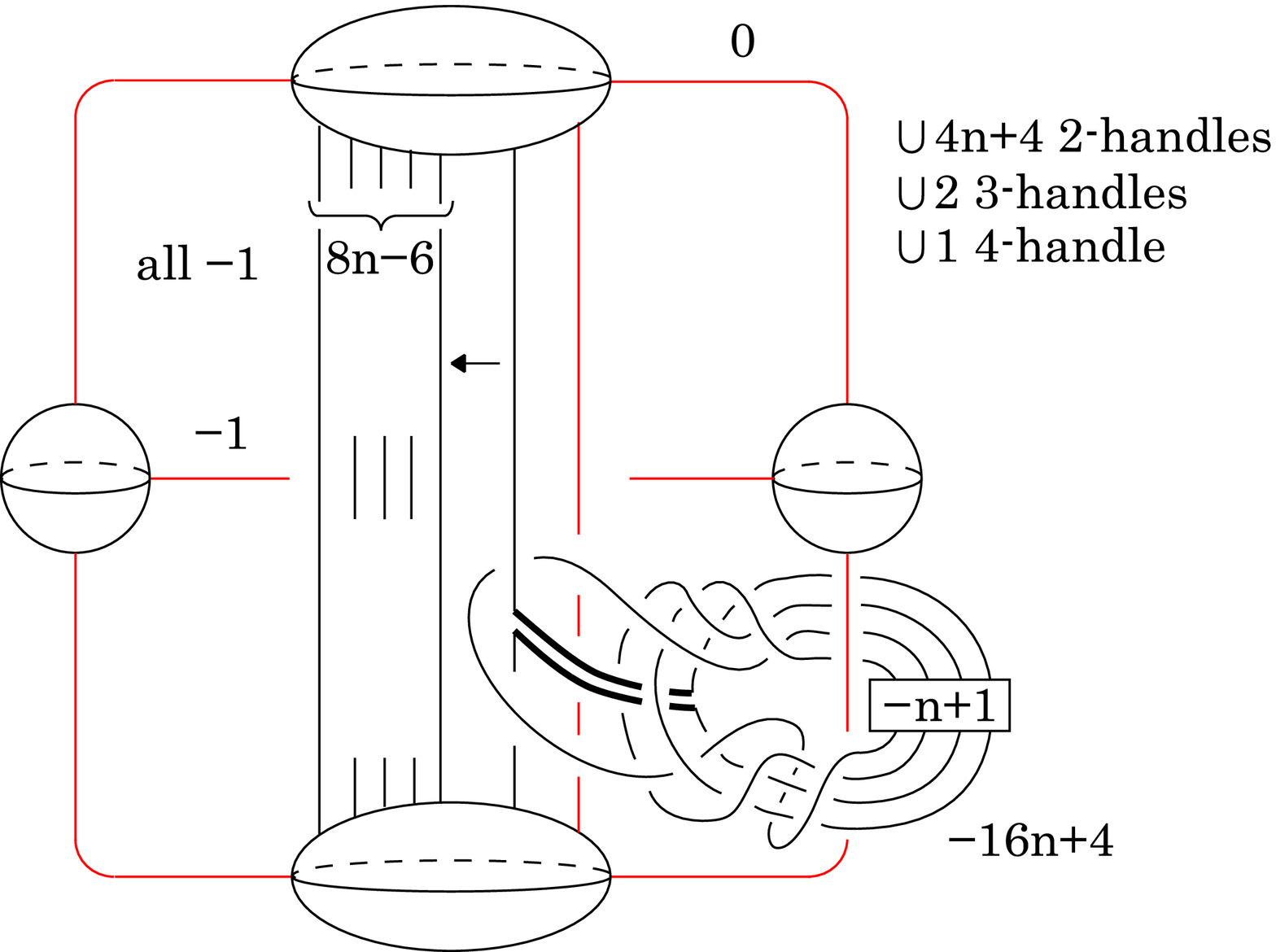}
\caption{$E(n)_4$}
\label{fig46}
\end{center}
\end{figure}
\begin{figure}[p]
\begin{center}
\includegraphics[width=3.1in]{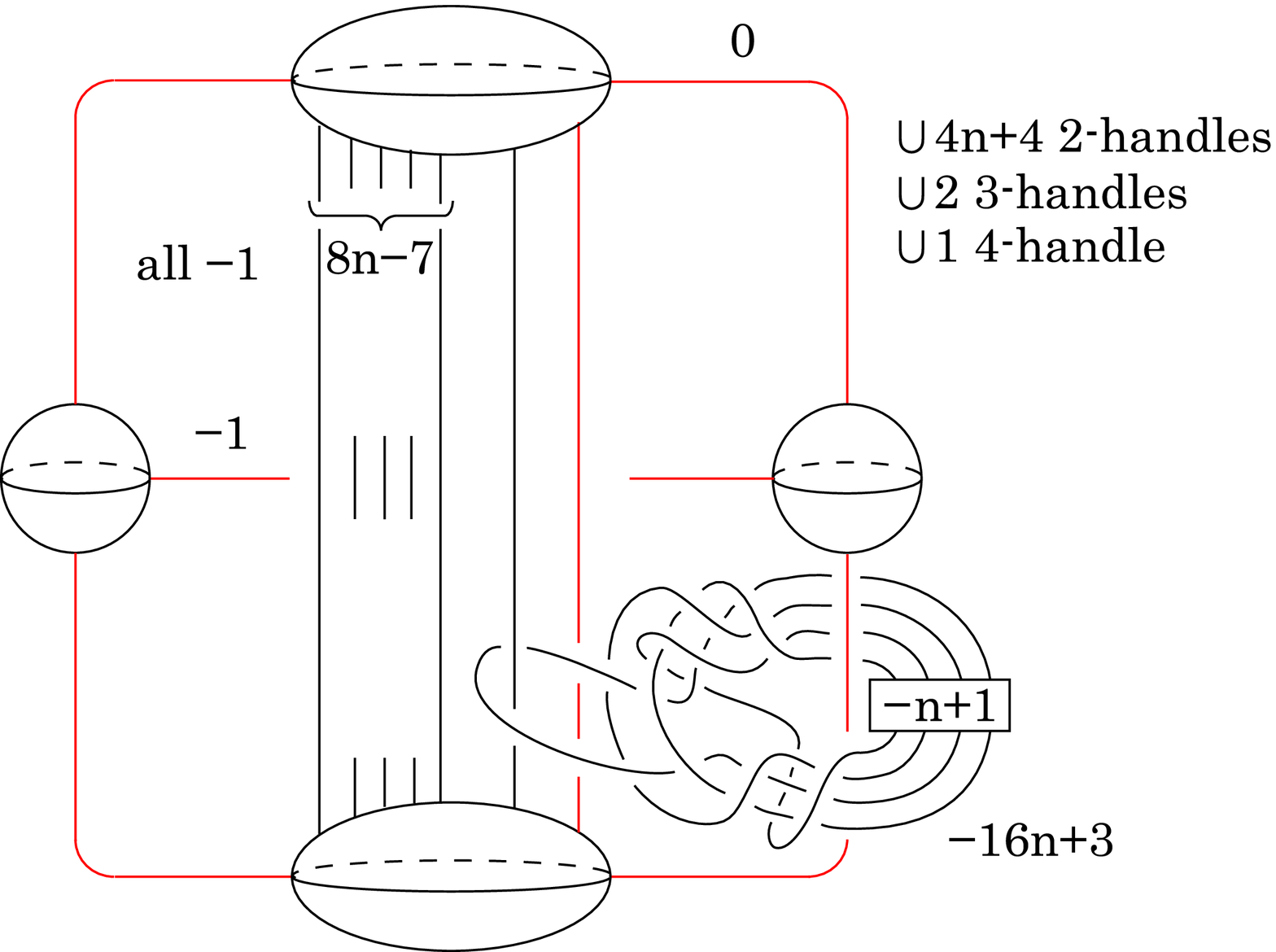}
\caption{$E(n)_4$}
\label{fig47}
\end{center}
\bigskip \medskip

\begin{center}
\includegraphics[width=3.1in]{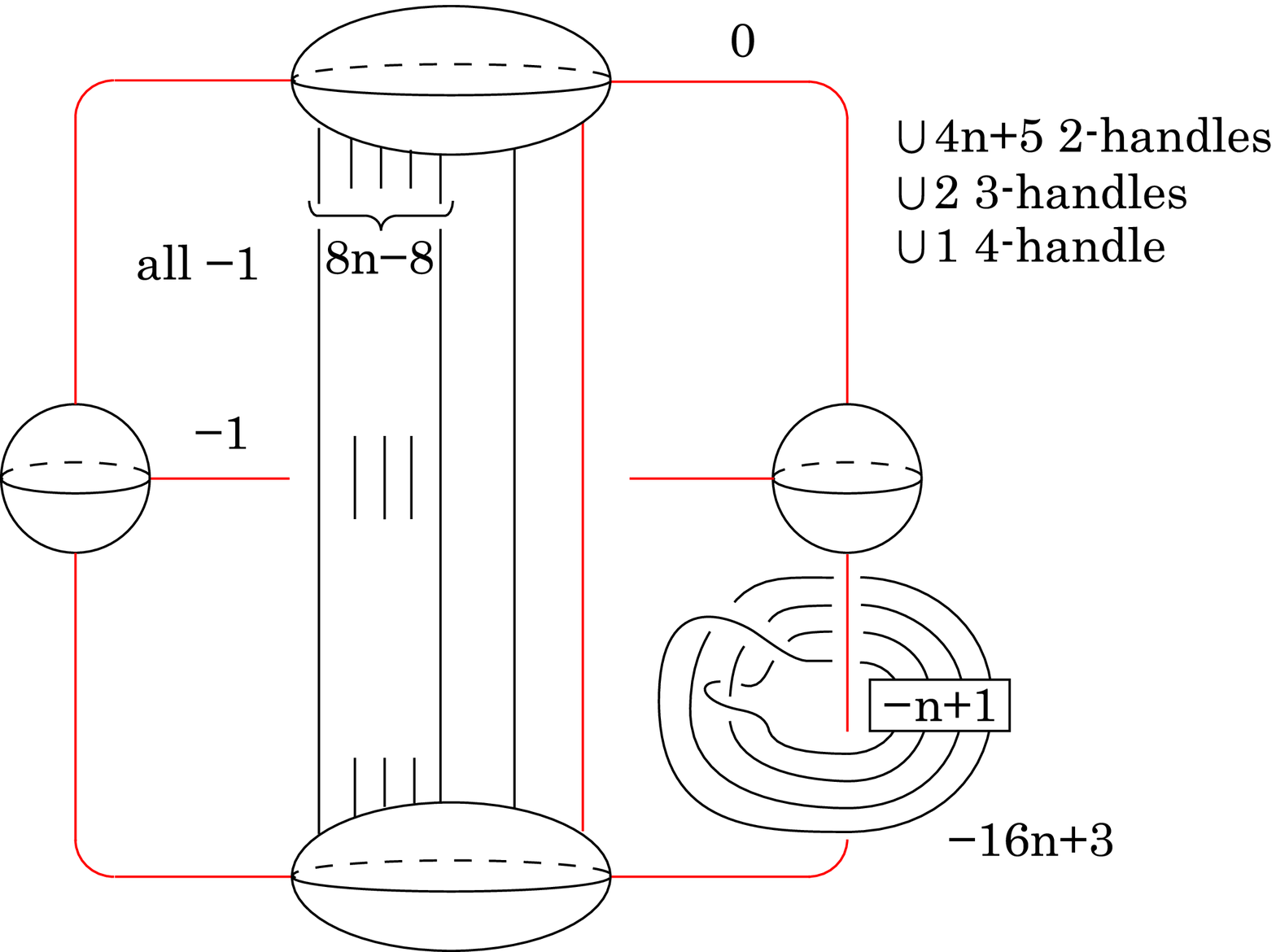}
\caption{$E(n)_4$}
\label{fig48}
\end{center}
\bigskip \medskip

\begin{center}
\includegraphics[width=3.1in]{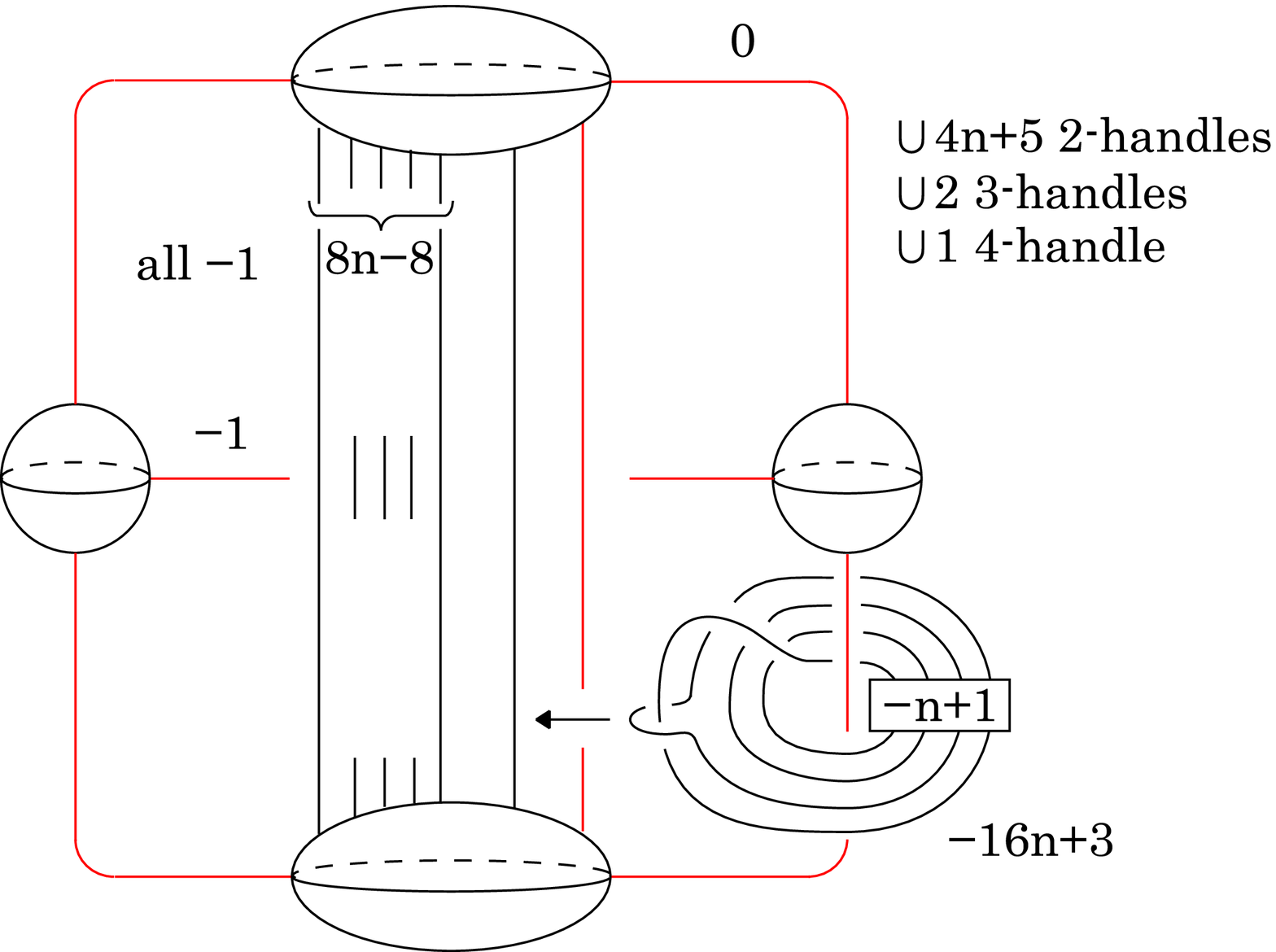}
\caption{$E(n)_4$}
\label{fig49}
\end{center}
\end{figure}
\begin{figure}[p]
\begin{center}
\includegraphics[width=3.1in]{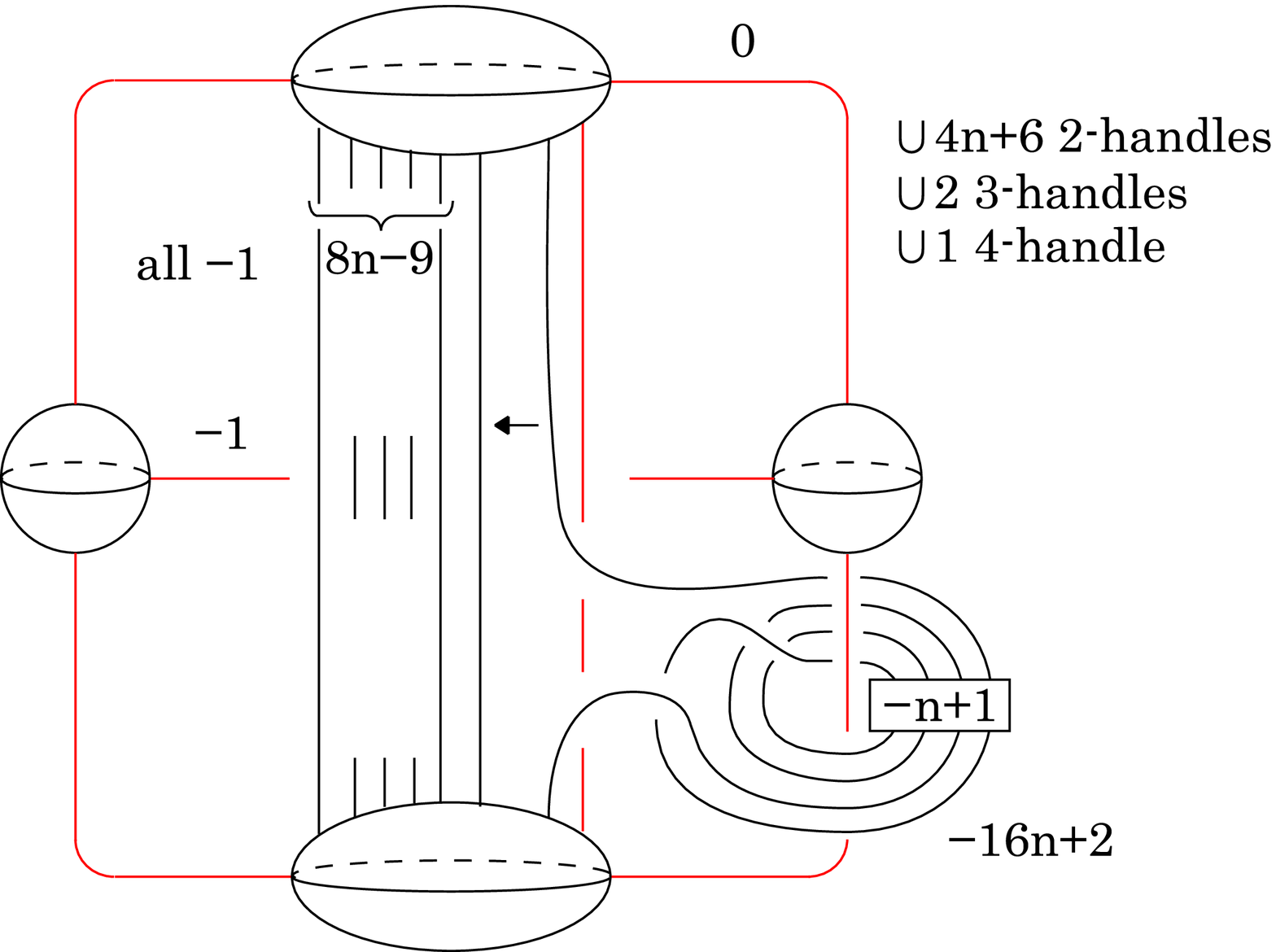}
\caption{$E(n)_4$ $(n\geq 2)$}
\label{fig50}
\end{center}
\bigskip \medskip

\begin{center}
\includegraphics[width=3.1in]{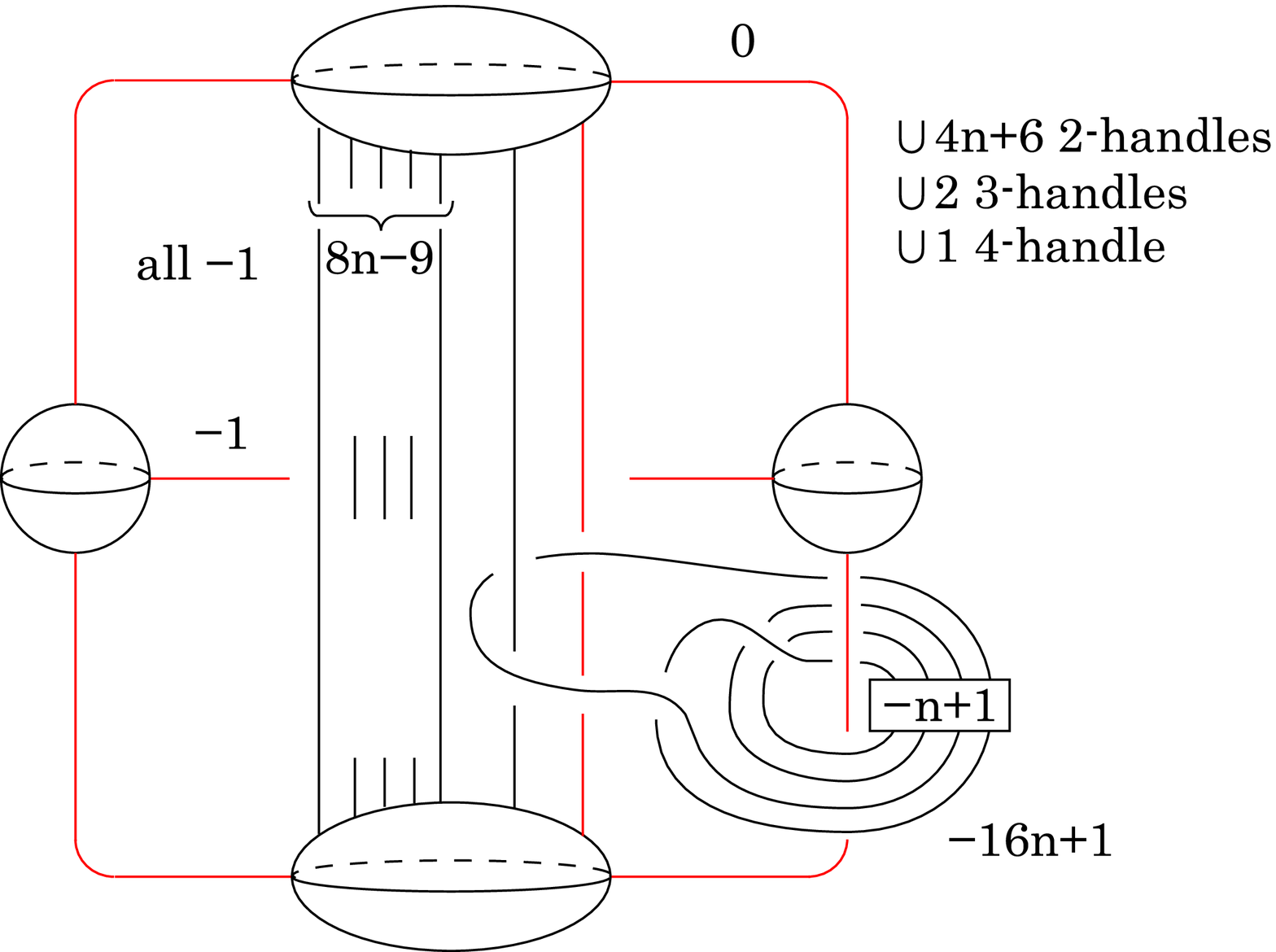}
\caption{$E(n)_4$ $(n\geq 2)$}
\label{fig51}
\end{center}
\bigskip \medskip

\begin{center}
\includegraphics[width=3.1in]{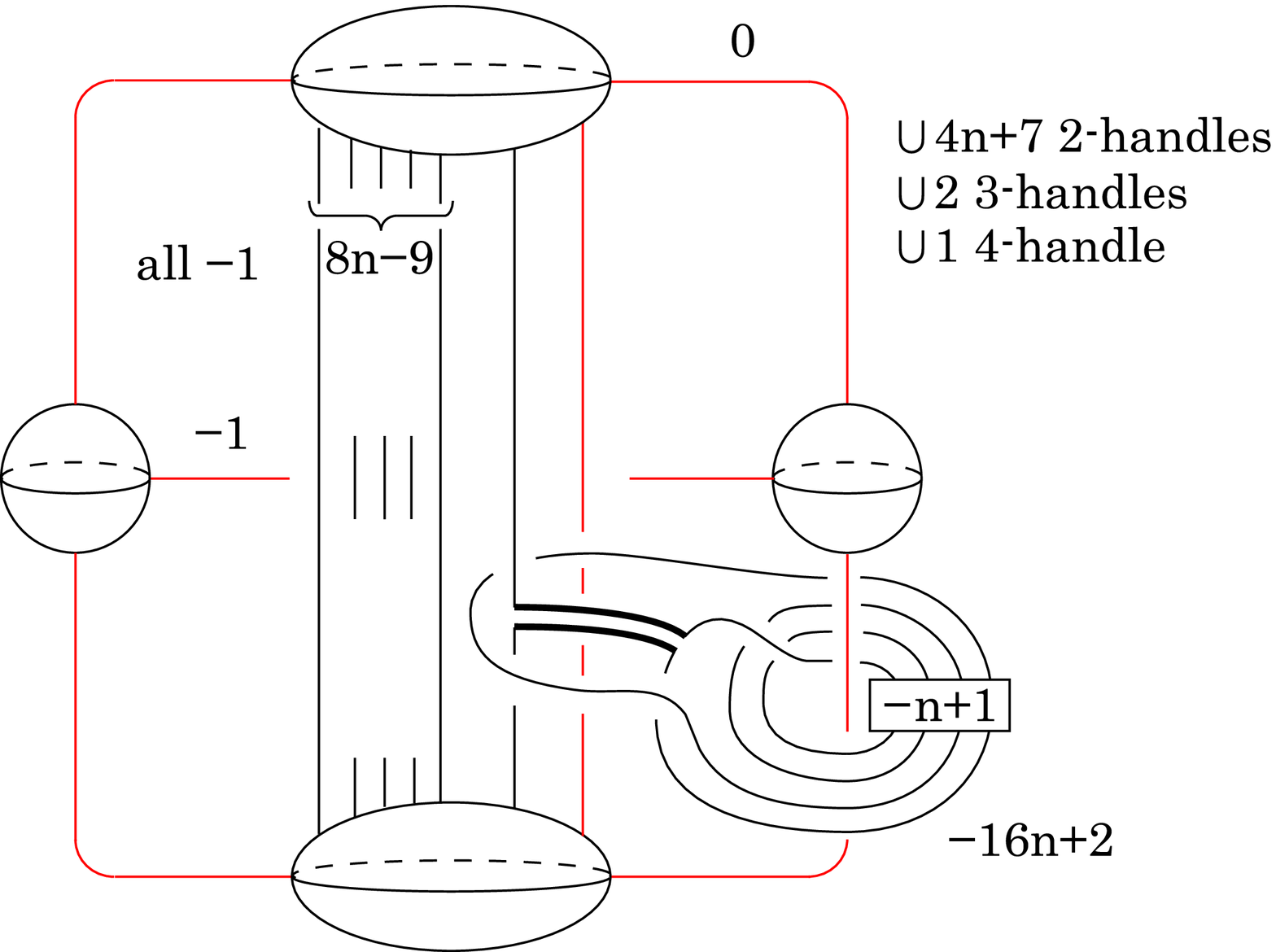}
\caption{$E(n)_4$ $(n\geq 2)$}
\label{fig52}
\end{center}
\end{figure}
\begin{figure}[p]
\begin{center}
\includegraphics[width=3.1in]{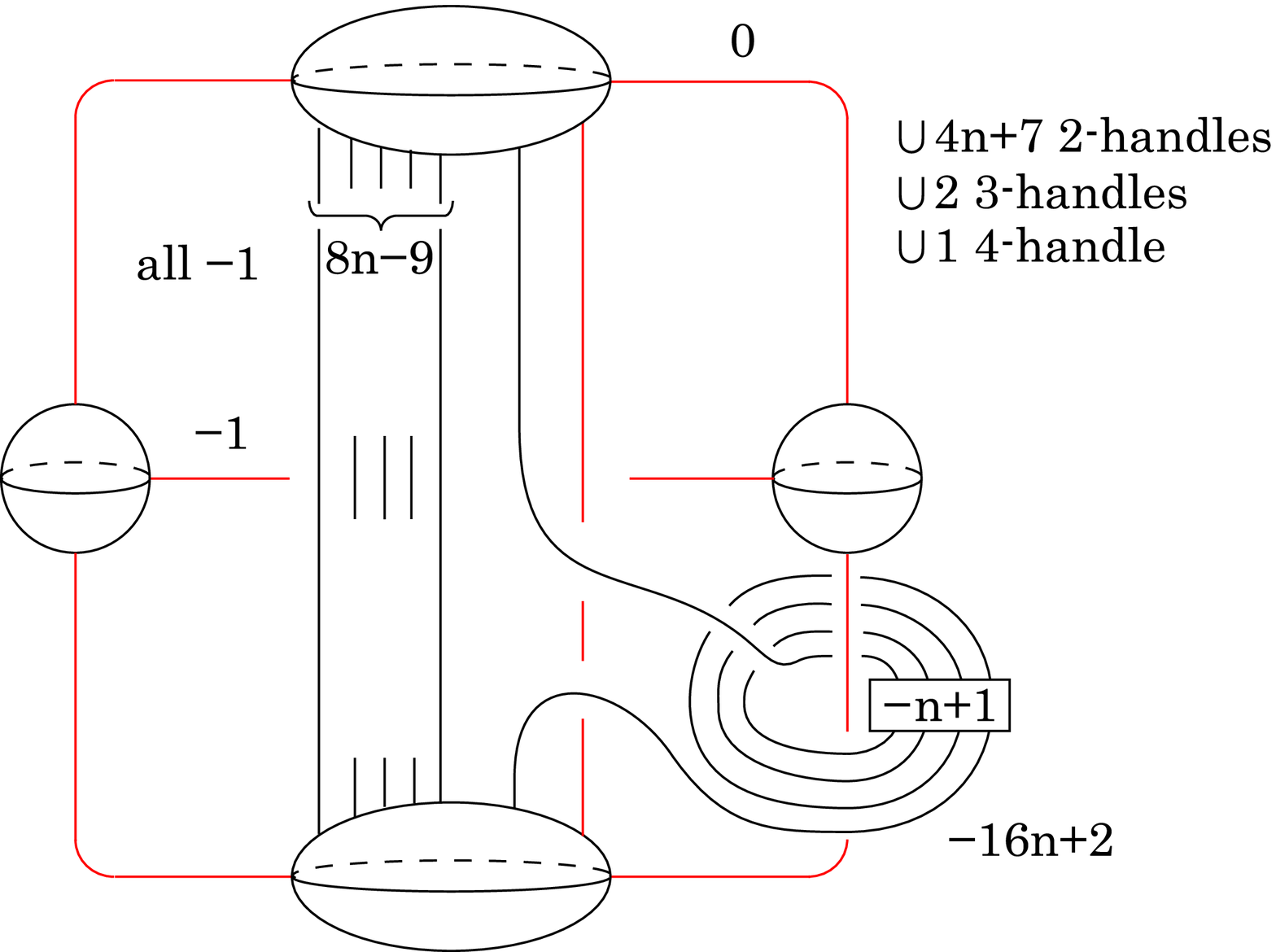}
\caption{$E(n)_4$ $(n\geq 2)$}
\label{fig53}
\end{center}
\bigskip \medskip

\begin{center}
\includegraphics[width=3.1in]{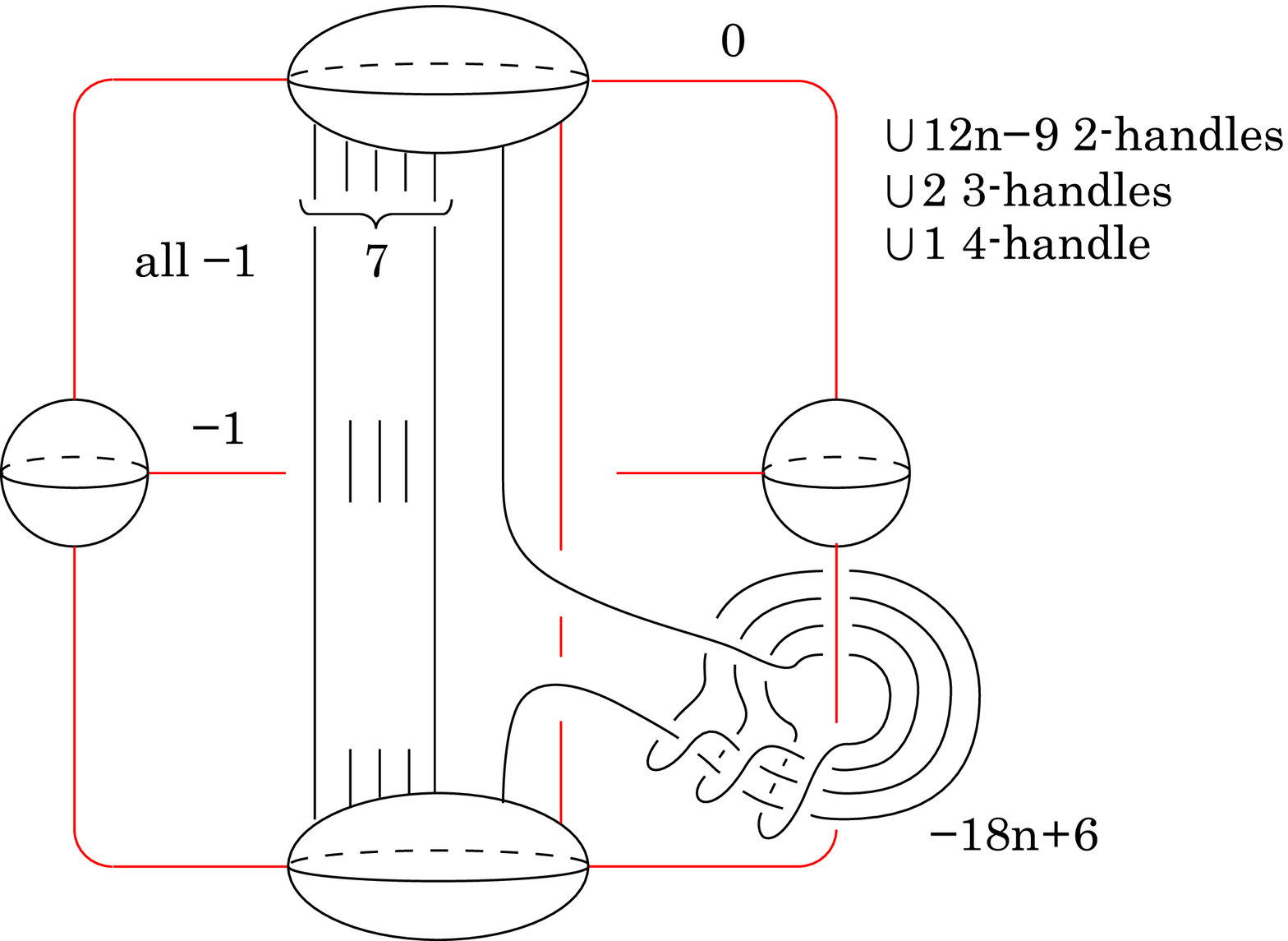}
\caption{$E(n)_4$}
\label{fig54}
\end{center}
\bigskip \medskip

\begin{center}
\includegraphics[width=3.1in]{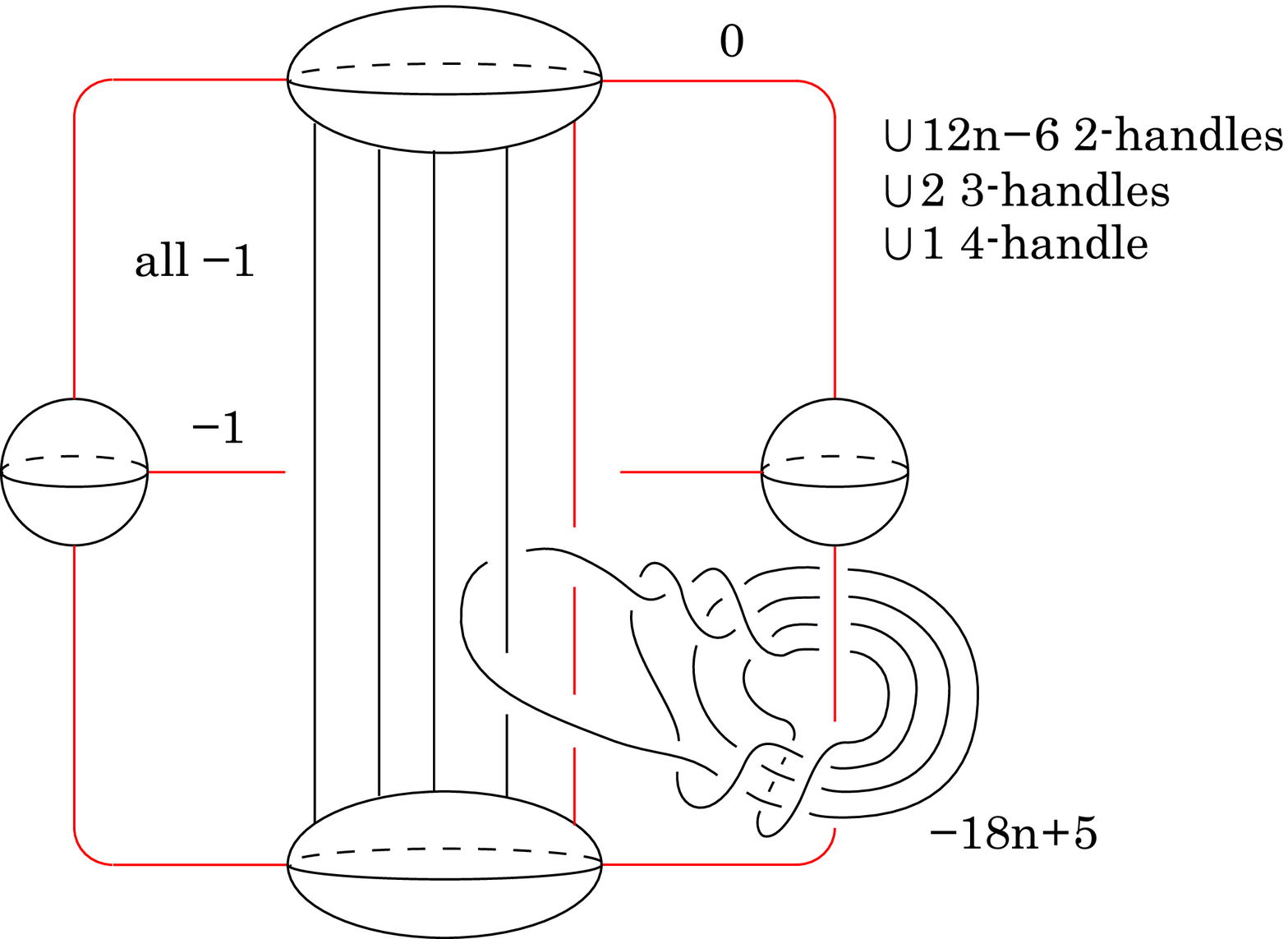}
\caption{$E(n)_4$}
\label{fig55}
\end{center}
\end{figure}
\begin{figure}[p]
\begin{center}
\includegraphics[width=2.9in]{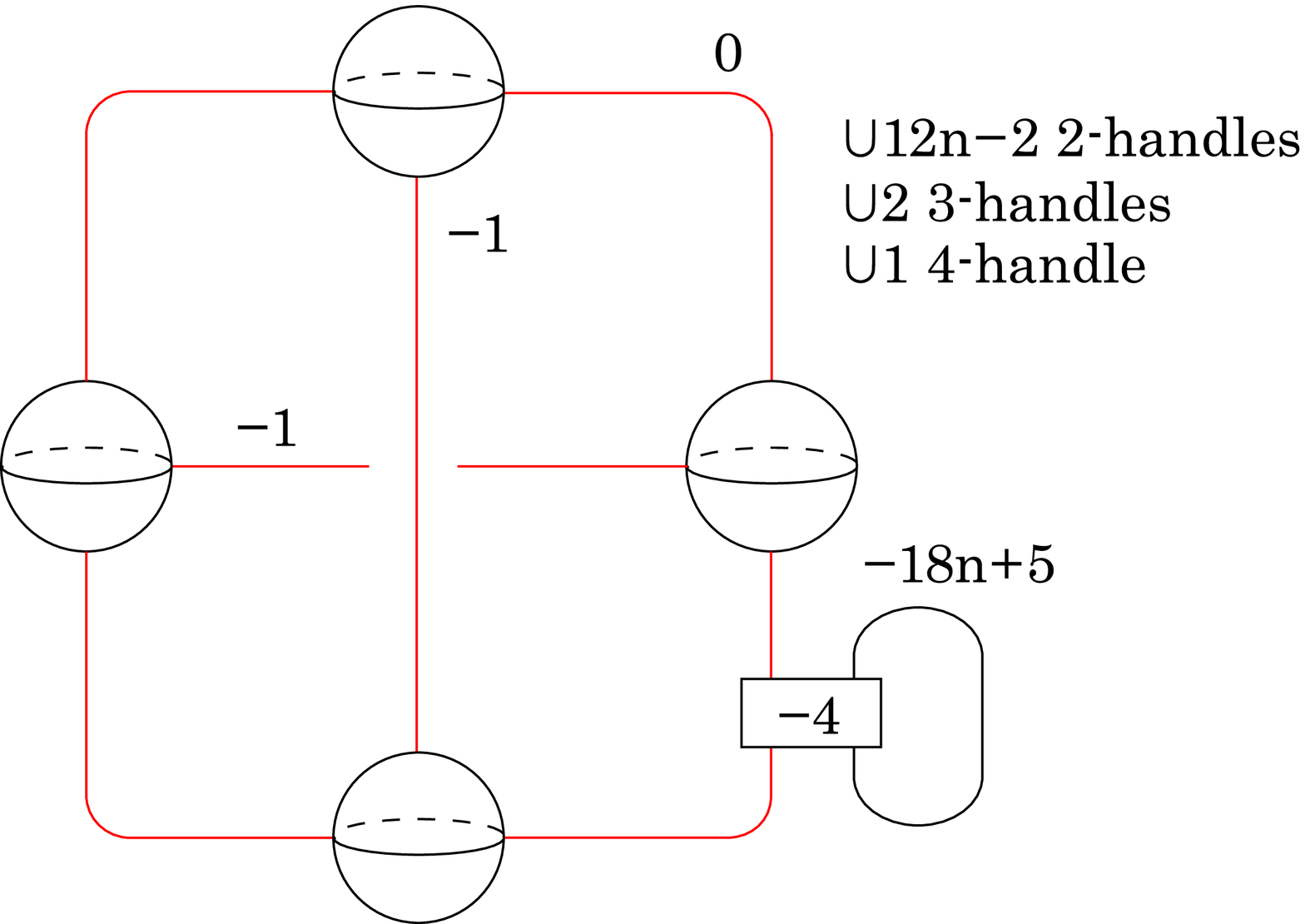}
\caption{$E(n)_4$}
\label{fig56}
\end{center}
\end{figure}
\end{document}